\documentclass[12pt, leqno]{amsart}

\usepackage{graphicx}
\usepackage{amssymb,amsmath,amsthm,amscd}
\usepackage{mathrsfs}
\usepackage{enumerate}
\usepackage{enumitem}
\usepackage{verbatim}
\usepackage[usenames,dvipsnames]{color}
\usepackage[colorlinks=true, pdfstartview=FitV,
linkcolor=blue,citecolor=blue,urlcolor=blue]{hyperref}
\usepackage[all]{xy}
\usepackage{tikz}
\usepackage{tikz-cd}
\usetikzlibrary{matrix, arrows}
\usepackage{stmaryrd} 
\usepackage{dsfont}
\usepackage{bm} %bold font for greek letters
\usepackage{scalerel}[2016/12/29] 
\usepackage{amsthm}
\usepackage{bbm}

\usepackage[linguistics]{forest}
\usepackage{rotating}

\allowdisplaybreaks[4]

\newcommand{\nc}{\newcommand}

\numberwithin{equation}{section}

\nc{\tomaszcomment}[1]{{\color{blue} Comment:\, #1}}

\usepackage[colorinlistoftodos]{todonotes}
\usepackage{mathtools}
\mathtoolsset{showonlyrefs,showmanualtags}
\usepackage{empheq}

\theoremstyle{plain}
\newtheorem{lem}{Lemma}[section]
\newtheorem{pro}[lem]{Proposition}
\newtheorem{thm}[lem]{Theorem}
\newtheorem{cor}[lem]{Corollary}
\newtheorem{defi}[lem]{Definition}

\newcommand{\Pro}{\begin{pro}}
	\newcommand{\enpro}{\end{pro}}
\newcommand{\Lem}{\begin{lem}}
	\newcommand{\enlem}{\end{lem}}
\newcommand{\Thm}{\begin{thm}}
	\newcommand{\enthm}{\end{thm}}
\newcommand{\Cor}{\begin{cor}}
	\newcommand{\encor}{\end{cor}}
\newcommand{\Defi}{\begin{defi}}
	\newcommand{\enDefi}{\end{defi}}
\newcommand{\Proof}{\begin{proof}}
	\newcommand{\enproof}{\end{proof}}

\theoremstyle{definition} %?
\newtheorem{rem}[lem]{Remark}

\newtheorem{Convention}[lem]{Convention}

\newcommand{\Conv}{\begin{Convention}}
	\newcommand{\enconv}{\end{Convention}}
\nc{\Rem}{\begin{rem}}
	\nc{\enrem}{\end{rem}}

\newcommand{\arxiv}[1]{\href{http://arxiv.org/abs/#1}{\tt arXiv:\nolinkurl{#1}}}

\newcommand{\monoto}{\hookrightarrow}
\nc{\epito}{\twoheadrightarrow}
\newcommand{\isoto}[1][]{\mathop{\xrightarrow%
		[{\raisebox{.3ex}[0ex][.3ex]{$\scriptstyle{#1}$}}]%
		{{\raisebox{-.6ex}[0ex][-.6ex]{$\mspace{2mu}\sim\mspace{2mu}$}}}}}

\newcommand{\commentout}[1]{}

\nc{\rmkend}{\hfill$\triangledown$}
\nc{\defend}{\hfill$\triangle$}

\DeclareMathSymbol{\shortminus}{\mathbin}{AMSa}{"39}

\nc{\kor}{\mathbb{C}}%{\Bbbk}
\nc{\fld}{\Bbbk}

\nc{\indx}{\mathbb{I}}
\nc{\Phipm}{\pmb{\phi^\pm}}
\nc{\gqaa}{U_q\widehat{\mathfrak{g}}}
\nc{\qla}{U_qL\mathfrak{g}} 
\nc{\qlae}{U_q^{\mathbb{K}}L\mathfrak{g}}
\nc{\Sym}{\on{Sym}}
\nc{\phipm}{\phi^\pm}
\nc{\uqsl}{U_q(\widehat{\mathfrak{sl}}_2)}
\nc{\ck}{\mathfrak{c}}
\nc{\Ck}{\mathfrak{C}}
\nc{\gOq}{\mathcal{O}_q(\widehat{\g})}
\nc{\gOqi}{\mathcal{O}_q^{[i]}}
\nc{\gOqc}{\mathcal{O}_q^{\mathbf{c}}(\widehat{\g})}
\nc {\gDrOq}{{}^{\mathrm{Dr}}\Oq(\widehat{\g})}

\nc{\degdr}{\on{deg}^{\mathrm{Dr}}}

\nc{\Tbr}{\mathbf{T}}
\nc{\ddA}{\ddot{A}}
\nc{\ide}{J}
\nc{\ddH}{\ddot{H}}

\nc{\Uu}{\widetilde{\mathbf{U}}}

\nc{\CC}{C}
\nc{\cc}{c}
\nc{\sss}{s}

\nc{\CCC}{\mathfrak{C}}
\nc{\ccc}{\mathfrak{c}}

\nc{\mac}[1]{\triangleleft_{#1}}
\nc{\Modfd}{\operatorname{Mod^{fd}}}
\nc{\Modfdc}{\operatorname{Mod}^{\operatorname{fd}}_\mathbf{c}}
\nc{\Modfdw}{\operatorname{Mod^{fdw}_\mathbf{c}}}
\nc{\Oqp}[1]{\Oq^{#1}}

\nc{\fext}[2]{{#1}[\negthinspace[#2]\negthinspace]}

\nc{\ev}[2]{\mathsf{V}_{#1}({#2})}
\nc{\evh}[1]{\operatorname{ev}_{#1}}

\nc{\xics}{\xi_{\mathbf{c}, \mathbf{s}}}
\nc{\cs}{\mathbf{c}, \mathbf{s}}
\nc{\cops}{\overline{\Delta}_{\cs}}
\nc{\tops}{\Delta_{\cs}}
\nc{\etacs}{\eta_{\cs}}
\nc{\etacz}{\eta_{\mathbf{c},0}}

\nc{\Bg}{B}
\nc{\Ag}{A}
\nc{\As}{\pmb{\Ag}(z)}
\nc{\Asmone}{\pmb{\Ag}^{(-1)}(z)}
\nc{\Aps}{\pmb{\Ag}_+(z)}
\nc{\Apsi}{\pmb{\Ag}_{i,+}(z)}
\nc{\Apsone}{\pmb{\Ag}^{(1)}_+(z)}
\nc{\Apsvar}[1]{\pmb{\Ag}^{(#1)}_+(z)}
\nc{\Apsmone}{\pmb{\Ag}^{(-1)}_+(z)}
\nc{\Ams}{\pmb{\Ag}_-(z)}
\nc{\Hg}{H}
\nc{\Thg}{\Theta}
\nc{\Thgs}{\pmb{\Theta}}
\nc{\Thgsr}{\pmb{\grave{\Theta}}}
\nc{\Thgsa}{\pmb{\acute{\Theta}}(z)} 

\nc{\Oqi}{\mathcal{O}_q^{[i]}}

\nc{\bThg}{\overline{\Theta}}
\nc{\gThg}{\grave{\Theta}}
\nc{\lam}{\pmb{\lambda}}
\nc{\poles}{\on{poles}}
\nc{\pol}{\on{pol}}

\DeclareRobustCommand{\qbinom}{\genfrac{[}{]}{0pt}{}}
\nc{\Oqc}{\mathcal{O}_q^{\mathbf{c}}}
\nc{\Oqg}{\mathcal{O}_q^{c_1}}

\nc{\DD}{\pmb{D}}
\nc{\KK}{\mathbb{K}}

\nc{\xgp}{x^+}
\nc{\Xgp}{\pmb{x}^+}
\nc{\xgm}{x^-}
\nc{\Xgm}{\pmb{x}^-}
\nc{\xgpm}{x^{\pm}}
\nc{\psig}{\phi^+}
\nc{\phig}{\phi^-}
\nc{\Psig}{\pmb{\phi}^+}
\nc{\Phig}{\pmb{\phi}^-}
\nc{\hg}{h}

\nc{\Eg}{e}
\nc{\Kg}{K}

\DeclareRobustCommand{\sqbin}{\genfrac{[}{]}{0pt}{}}

\nc{\qq}{(q-q^{-1})^{-1}}
\nc{\factor}{\Omega}

\nc{\Oq}{\mathcal{O}_q}
\nc{\car}{\mathcal{H}}
\nc{\qaa}{U_qL\mathfrak{sl}_2}
\nc{\drqaa}{{}^{\mathrm{Dr}}\qaa}
\nc{\usl}{U_q\mathfrak{sl}_2}
\nc{\uLsl}{U_qL\mathfrak{sl}_2}
\nc{\ugLsl}{U_q^{\fld}L\mathfrak{sl}_2}
\nc{\Serre}{\mathsf{Serre}}
\nc{\UXp}{UX_+}

\nc{\ada}{\ad_{\bhg_{-1}}}
\nc{\adb}{\ad_{\bhg_{1}}}
\nc{\adc}{\ad_{\bHg_{1}^{(2)}}}

\nc{\tAg}{\widetilde{A}}
\nc{\hAg}{\widehat{A}}

\nc{\bHg}{\overline{\Hg}}
\nc{\bhg}{\overline{\hg}}
\nc{\tHg}{\widetilde{H}}

\nc {\DrOq}{{}^{\mathrm{Dr}}\mathcal{O}_q}
\nc{\gup}[1]{^{(#1)}}
\nc{\gupp}{^{(1),+}}
\nc{\mi}{^{-1}}
\nc{\adh}{\operatorname{ad}_{\bar{\hg}_{-1}}}
\nc{\adhp}{\operatorname{ad}_{\bar{\hg}_{1}}}
\nc{\Omg}{\Omega^{-1}}
\nc{\Htwo}{\overline{H}_1^{(2)}}
\nc{\ad}{\operatorname{ad}}

\nc{\cartoq}{\mathcal{H}_{\Oq}}
\nc{\cart}{\mathcal{H}}

\newcommand{\on}{\operatorname}
\nc{\be}{\begin{enumerate}}
	\nc{\ee}{\end{enumerate}}
\newcommand{\eq}{\begin{equation}}
	\newcommand{\eneq}{\end{equation}}
\nc{\bc}{\begin{cases}}
	\nc{\ec}{\end{cases}}
\newcommand{\eqn}{\begin{eqnarray*}}
	\newcommand{\eneqn}{\end{eqnarray*}}
\newcommand{\ba}{\begin{array}}
	\newcommand{\ea}{\end{array}}

\newcommand{\C}{{\mathbb C}}

\newcommand{\Z}{{\mathbb Z}}

\newcommand{\g}{{\mathfrak{g}}}
\nc{\Ad}{\operatorname{Ad}}

\nc{\gr}{\on{gr}} 
\nc{\ourR}{\mathbf{R}}
\nc{\ourQ}{\mathbf{Q}}

\newcommand{\Hom}{\operatorname{Hom}}
\newcommand{\End}{\operatorname{End}}
\nc{\Aut}{\operatorname{Aut}}

\nc{\coker}{\operatorname{coker}}

\nc{\Img}{\on{Im}}
\nc{\res}{\on{res}}

\nc{\modv}[1]{{#1}\operatorname{-mod}}

\nc{\bl}{\bigl(}
\nc{\br}{\bigr)}

\nc{\ie}{{\em i.e.}, }
\nc{\eg}{{\em e.g.}, }
\nc{\ten}{\otimes}
\nc{\id}{\operatorname{id}}
\nc{\wt}{\widetilde}

\newlength{\mylength}
\DeclareRobustCommand{\SkipTocEntry}[5]{}

\AtBeginDocument{%
   \def\MR#1{}
}

\title[Drinfeld rational fractions for quantum symmetric pairs]
{Drinfeld rational fractions for affine Kac--Moody quantum symmetric pairs}

\author[T. Prze\'{z}dziecki]{Tomasz Prze\'{z}dziecki}
\address{School of Mathematics, University of Edinburgh, Peter Guthrie Tait Rd, Edinburgh, EH9 3FD, United Kingdom, 
OrciD: 0000-0001-9700-1007}
\email{\href{mailto:tprzezdz@exseed.ed.ac.uk}{tprzezdz@exseed.ed.ac.uk}}



\thanks{The author was supported by the ERC Starting Grant No. 759967
	{\em Categorified Donaldson--Thomas Theory}, and the EPSRC grant No.\ EP/W022834/1 \emph{Kac--Moody quantum symmetric pairs, KLR algebras and generalized Schur--Weyl duality}.}

\begin{document}

\begin{abstract} 
We formulate a precise connection between the new Drinfeld presentation of a quantum affine algebra $U_q\widehat{\g}$  and the new Drinfeld presentation of affine coideal subalgebras of split type recently discovered by Lu and Wang. 
In particular, we establish a ``factorization formula'', expressing the commuting ``Drinfeld--Cartan''-type operators $\Theta_{i,k}$ in the coideal subalgebra in terms of the corresponding Drinfeld generators of $U_q\widehat{\g}$, modulo the ``Drinfeld positive half" of $U_q\widehat{\g}$. We study the spectra of these operators on finite dimensional representations, and describe them in terms of rational functions with an extra  symmetry. These results can be seen as the starting point of a $q$-character theory for affine Kac--Moody quantum symmetric pairs. Additionally, we prove a compatibility result linking Lusztig's and the Lu-Wang-Zhang braid group actions. 
\end{abstract}

\maketitle

\setcounter{tocdepth}{1}
\tableofcontents

\section{Introduction}

\subsection{A trinity of presentations} 

It is well known that quantum affine algebras~$U_q\widehat{\g}$ admit three distinct presentations: the Kac--Moody (also:\ Drinfeld--Jimbo, Serre) presentation \cite{jimbo-85, drinfeld-quantum-groups-87}, the new Drinfeld (or loop) presentation \cite{drinfeld-dp, beck-94, Damiani12, Damiani15} and the FRT (also:\ RTT) presentation \cite{Res-STS}. The first of these can be seen as a quantization of the usual presentation of an affine Kac--Moody Lie algebra, while the second quantizes its construction as a central extension of a loop Lie algebra. Precise isomorphisms between these  different presentations 
were described in \cite{Ding-Fr, FrMukHopf}. 

It has recently emerged that the same trinity of presentations exists in the world of quantum symmetric pairs. A (classical)  symmetric pair $(\g, \g^\theta)$ consists of a semisimple Lie algebra $\g$, together with a Lie subalgebra fixed under an involution $\theta$. We note that~$\g^\theta$ is reductive. The classification of symmetric pairs is equivalent to the classification of real forms of complex semisimple Lie algebras, and can be formulated combinatorially in terms of Satake diagrams. For example, given an appropriate involution, $(\mathfrak{sl}_n, \mathfrak{so}_n)$ is a symmetric pair. 
However, the situation at the quantum level is more nuanced. 
%Since the fixed point subalgebra $\g^\theta$ is reductive, it possesses a Drinfeld--Jimbo quantization as a Hopf algebra. 
The Drinfeld--Jimbo quantizations of a Lie algebra and its fixed point subalgebra are typically not compatible. The correct quantization of a symmetric pair was proposed by Letzter \cite{letzter-99, letzter-02}. According to her definition, a quantum symmetric pair consists of a quantum group and a coideal subalgebra (rather than a Hopf subalgebra), which specializes to $U(\mathfrak{g}^\theta)$ in the limit $q \mapsto 1$. 

The generalization of quantum symmetric pairs to the Kac--Moody case was carried out by Kolb \cite{kolb-14}. In particular, he provided a description of coideal subalgebras by generators and Serre-like relations. 
The RTT presentation of coideal subalgebras was constructed in special cases by Molev, Ragoucy and Sorba \cite{MRS-03}, and was preceded by earlier work on twisted Yangians \cite{Olshan}. 
On the other hand, the Drinfeld-type presentation of affine quantum symmetric pairs is a much more recent development. 
The first steps in this direction were taken by Baseilhac and Kolb in \cite{bas-kol-20}, where $q$-root vectors and the  relations between them were computed in one of the rank $1$ cases (the $q$-Onsager algebra). 
Their approach was inspired by Damiani's work for $\widehat{\mathfrak{sl}}_2$ \cite{damiani93}, and involved the construction of a new braid group action on the $q$-Onsager algebra, using a computer-based proof. A computer-free proof of the existence of this braid group action was given independently by Terwilliger in \cite{Ter18}. 
A Drinfeld-type presentation for all simply laced split Satake diagrams was later constructed by Lu and Wang \cite{lu-wang-21}, using a more general braid group action stemming from a Hall algebra realization of coideal subalgebras \cite{LuWangKM}. A generalization to non simply-laced split Satake diagrams was obtained by Zhang \cite{ZhangDr}. 

\subsection{New developments} 

The theory of affine Kac-Moody quantum symmetric pairs (aKM QSP) is a very dynamic field, with many significant advances in recent years. 
In order to frame our work in a broader context, we list the most relevant ones below: 
\begin{enumerate}[label=\roman*)]
\item the new Drinfeld presentation for aKM QSP of split type \cite{bas-kol-20, lu-wang-21, ZhangDr}; 
\item relative braid group symmetries on QSP \cite{KolbPel, WangZhang, WangZhang2, ZhangBr}, generalizing Lusztig's braid group action; 
\item the Hall algebra approach to aKM QSP \cite{LuRuanWang, lu-ruan-wang-23, LuWangKM};
\item universal K-matrices for aKM QSP \cite{AppelVlaar2, AppelVlaar}; 
\item generalization of the Kang--Kashiwara--Kim functor, a type of Schur--Weyl duality, to aKM QSP \cite{AppelPrz}, forging a link with Khovanov--Lauda--Rouquier-type algebras. 
\end{enumerate} 

In the present paper we investigate the implications of i) and ii), i.e., the new Drinfeld presentation and braid group action, with view to developing the representation theory of aKM QSP. 
We expect our work will also exhibit interesting connections to iv) and~v) - this is addressed in more detail in \S \ref{subsec: appl future dir}. 

\subsection{Classical motivations}  

In the seminal works of Bao and Wang \cite{bao-wang-18, bao-wang-18b}, a general philosophy was proposed that most fundamental algebraic, geometric, and categorical constructions from quantum group theory should be generalizable to quantum symmetric pairs. 
Around the same time, advances in this direction were also made in the influential papers \cite{balagovic-kolb-19} and \cite{ehrig-stroppel}. 
Accordingly, Lu and Wang's discovery of a Drinfeld-type presentation for aKM QSP can be seen as a major step towards generalizing many classical ideas from the theory of quantum affine algebras. Most importantly, these include: 
\be[label=\roman*)] 
\item the work of Chari and Pressley on \emph{Drinfeld polynomials}, their relation to spectra of certain commuting operators in the Drinfeld presentation, and their role in the classification of finite dimensional representations \cite{chari-pressley-qaa, chari-pressley-94}; 
\item the work of Frenkel, Mukhin and Reshetikhin on $q$\emph{-characters}, associated algorithms, connection to deformed $\mathcal{W}$-algebras and cluster theory \cite{FrenRes, FrenMuk-comb}, 
\item Hernandez's theory of quantum affinizations and the fusion product \cite{hernandez-05, hernandez-07};  
\item Nakajima's geometric construction of finite dimensional representations \cite{Nakajima-01}. 
\ee

The main goal of this paper is to develop a QSP generalization of i), and provide a foundation for ii). A comprehensive study of $q$-characters for quantum symmetric pairs will be presented in a forthcoming paper with Li \cite{li-przez} (see also \cite{LiPrz} and \S \ref{subsec: q char} below). 
Some potential implications of our work for iii) are also discussed in \S \ref{subsec: appl future dir}. 

\subsection{A shift of perspective} 

In the foundational works of Letzter \cite{letzter-99, letzter-02} and Kolb \cite{kolb-14}, the coideal structure is of primary importance. 
However, QSP coideal subalgebras admit a description by generators and relations, and may, accordingly, be regarded as an independent algebraic object (often called an $\imath$quantum group), forgetting the embedding into the quantum group. This point of view is more prevalent in recent literature. It is also rather natural for several reasons. Firstly, $\imath$quantum group constructions are not simple restrictions of the corresponding quantum group constructions. For example, the braid group action from \cite{WangZhang2} is not a restriction of Lusztig's braid group action. Secondly, the proofs are often more technically involved and challenging than their quantum group counterparts, exhibiting many features absent from the latter. Examples include a lack of triangular decompositions and canonical raising and lowering operators, the presence of many rank $1$ cases and parameters, as well as non-obvious integral forms. 

Nevertheless, from the point of view of representation theory, considering an \linebreak $\imath$quantum group as embedded inside a quantum group has two key advantages. Firstly, it gives access to restriction functors, which are a natural source of representations. For example, in the case of the $q$-Onsager algebra, essentially every irreducible finite dimensional representation can be obtained as a restriction of an $\uLsl$-representation \cite{ito-ter-10}. Secondly, the coproduct endows the category of representations of a coideal subalgebra with a monoidal action. The monoidal structure of the category of finite dimensional representations of a quantum affine algebra has been a subject of intense study, with many important applications to cluster theory \cite{kashiwara-02, Hernandez-kr, Hernandez-stp, Her-Lec-10, Her-Lec-15}. The corresponding monoidal actions are likely to be equally interesting. 

\subsection{The main results} 

These considerations lead to the main problem addressed in this paper: 

\emph{What is the relation between the Drinfeld presentation of a quantum affine algebra and Lu and Wang's Drinfeld-type presentation of a QSP coideal subalgebra?} 

Remarkably, we establish that these two presentations are, in fact, mutually compatible in a suitable sense.  
Before stating our results, 
let us make a few remarks about the basic assumptions and notations. Let $\g$ be a simple Lie algebra of type $\mathsf{ADE}$. We consider the corresponding untwisted quantum affine algebra $U_q\widehat{\g}$, and its coideal subalgebras $\gOq$ of split type. The latter are also often referred to as generalized $q$-Onsager algebras. The coideal subalgebra $\gOq$ depends on two sets of parameters $\mathbf{c}$ and $\mathbf{s}$\footnote{In this introduction, in order to simplify notation, we will suppress the dependence on the parameters. In other words, $\gOq$ refers here to a QSP coideal subalgebra with chosen $\mathbf{c}$ and $\mathbf{s}$. Later, in \S \ref{sec: q-ons} and \S \ref{sec: qsp spl}, $\gOq$ will refer to the \emph{universal} version of this algebra, where the parameters $\mathbf{c}$ are replaced by central elements. The coideal subalgebras with fixed parameters will then be denoted as $\Oqp{\mathbf{c},\mathbf{s}}(\widehat{\g})$.}. The case $\mathbf{s} = 0$ is called standard. 
We use the usual notation for the Drinfeld generators of $U_q\widehat{\g}$, 
with ``raising operators'' $x_{i,k}^+$, ``lowering operators'' $x_{i,k}^-$ and ``Cartan operators'' $h_{i,k}, \phi^{\pm}_{i,k}$ 
 (see \S \ref{subsec: qaa gen} for details).   
The coideal subalgebra $\gOq$, in its Drinfeld presentation, has two main kinds of generators: $A_{i,k}$ (somewhat analogous to $x_{i,k}^\pm$ in $U_q\widehat{\g}$ and denoted $B_{i,k}$ in \cite{lu-wang-21}) and $\Theta_{i,k}$ (analogous to $\phi_{i,k}^\pm$ in $U_q\widehat{\g}$). When considering the generators $\Theta_{i,k}$ as a series, we use a slightly different normalization, denoted $\Thgsr_i(z)$, than in \cite{lu-wang-21}. The reason for it is explained in Remark \ref{rem: Thg grave more natural}. 

\subsubsection{Factorization and coproduct theorems} 

Of particular importance in the Drinfeld presentation of $U_q\widehat{\g}$ is the (almost) commutative subalgebra generated by the elements~$\phi_{i,k}^\pm$. The spectra of these operators are used to define Drinfeld polynomials and $q$-characters. In a coideal subalgebra $\gOq$, there exists an analogous commutative subalgebra generated by elements $\Theta_{i,k}$. Our key result, which we call the \emph{factorization theorem}, enables us to express the series $\Thgsr_i(z)$ in terms of the series $\pmb{\phi}^\pm_i(z)$, modulo the ``Drinfeld positive half'' of $U_q\widehat{\g}$, i.e., the subalgebra spanned by elements of positive degree in the Drinfeld gradation, where
$\degdr \xgpm_{i,k} = \pm \alpha_i, \ \degdr h_{i,k} = 0$ 
 (see \S \ref{subsec: qaa gen} for details). 
%, i.e., the subalgebra $U_+$ spanned by elements of positive $Q$-degree (where $Q$ is the finite root lattice). 

\Thm \label{intro main thm} 
Let $\g$ be of type $\mathsf{A}$. Then: 
\eq 
\Thgsr_i(z) \equiv \pmb{\phi}_i^-(z\mi)\pmb{\phi}_i^+(\CC z)  \fext{\mod U_{+}}{z}.  
\eneq 
\enthm 

Here $C$ is a scalar depending on the choice of parameters $\mathbf{c}$. 
For non-standard parameters $\mathbf{s} \neq 0$, the statement of Theorem \ref{intro main thm} is slightly more technical, see Theorem~\ref{cor: group like Theta} for details. The restriction to type $\mathsf{A}$ stems from our use of an explicit combinatorial argument related to braid groups (see \S \ref{sssub: compat}--\ref{sssub: complex}  below). %We are confident that this argument can be generalized to other types, which will be explored in future work. 
Successively more general versions of Theorem \ref{intro main thm} are proven in Theorem \ref{thm: fact} (rank one, standard parameters), Theorem \ref{cor: group like Theta} (rank one, nonstandard parameters) and Corollary \ref{cor: ultimate factorization theorem} (arbitrary rank).

As a consequence of Theorem \ref{intro main thm}, we obtain the following coproduct formula. 

\Thm \label{intro: coproduct} 
The series $\Thgsr_i(z)$ is ``approximately group-like'', i.e., 
\[
\Delta(\Thgsr_i(z)) \equiv \Thgsr_i(z) \otimes \Thgsr_i(z) \quad \mod \fext{(\Oqp{}(\widehat{\g}) \otimes U_{+})}{z}.  
\]
\enthm 

For the proof of Theorem \ref{intro: coproduct}, also see Theorem \ref{cor: group like Theta} and Corollary \ref{cor: ultimate factorization theorem} below. We remark that, since the initial release of this paper, type $\mathsf{BCD}$ versions of the two theorems above have been established in \cite{LiPrz}. 

\subsubsection{Drinfeld rational fractions} 

The Drinfeld presentation of quantum affine algebras is especially well-suited for studying their representation theory. A key role is played by spectra of the series $\pmb{\phi}_i^\pm(z)$ on finite dimensional representations. By the theorems of Chari--Pressley \cite{chari-pressley-qaa, chari-pressley-94} and Frenkel--Reshetikhin \cite{FrenRes}, they can be described as expansions of rational functions which are ratios of $q^2$-shifted polynomials. In the case of the highest weight vector, these polynomials are known as Drinfeld polynomials, and they suffice to classify finite dimensional irreducible representations. To formulate a sensible character theory, the full eigenspace decomposition with respect to the operators $\pmb{\phi}_i^\pm(z)$ and the corresponding eigenvalues are needed. 

As a natural generalization to the world of aKM QSP, we consider the spectra of the series $\Thgsr_i(z)$ on finite dimensional representations of $\gOq$. Theorem \ref{intro main thm} enables us to compute these spectra for all representations $V$ obtained via restriction from $U_q\widehat{\g}$. More precisely, they can be described in terms of rational functions $F(z)$ with a certain symmetry, which we call \emph{Drinfeld rational fractions}.    

\Thm \label{intro thm2}
Consider a generalized eigenspace of $\Thgsr_i(z)$ in $V$. 
The corresponding series of eigenvalues $\pmb{D}(z)$ is the expansion of a rational function of the form 
\eq 
\pmb{D}(z) = \frac{F(q^{-1}z)}{F(qz)}, 
\eneq
where $F(z)$ is a rational function satisfying the $C$-twisted unitarity property: 
\eq 
F(q\CC \mi z\mi) = F(q^{-1}z)\mi. 
\eneq
\enthm

Theorem \ref{intro thm2} (and its close cousins, such as Corollary \ref{cor: FR thm Oq} below) is a natural starting point for a $q$-character theory for aKM QSP. For a more precise definition of, as well as a way to compute, Drinfeld rational fractions, see \eqref{eq: F as ratio P P*}; for the proof of Theorem~\ref{intro thm2}, see Theorem \ref{cor: Drinfeld rational functions} below. 

\subsubsection{Compatibility of braid group actions} 
\label{sssub: compat}

It should be emphasized that the braid group action used to construct the Drinfeld presentation for coideal subalgebras in \cite{bas-kol-20, lu-wang-21}  is \emph{different} from the usual Lusztig braid group action used to construct the Drinfeld presentation for quantum affine algebras. In particular, the former is \emph{not} a restriction of the latter. 
Nevertheless, it turns out that, remarkably, the two actions are, roughly speaking, compatible modulo the ``Drinfeld positive half'' $U_+$ of $U_q\widehat{\g}$. Not only is this a key ingredient allowing us to generalize Theorem \ref{intro main thm} from the $q$-Onsager algebra to higher rank coideal subalgebras, but also a very interesting fact on its own. 
Together with Theorems \ref{intro main thm} and \ref{intro: coproduct}, it suggests a general philosophy that although various structures on coideal subalgebras and quantum affine algebras are distinct, they are in fact compatible modulo $U_+$. 

Let us make the aforementioned compatibility of braid group actions more precise. The Drinfeld presentations of both $U_q\widehat{\g}$ and $\gOq$ rely crucially on braid group operators corresponding to the fundamental weights $\omega_i$, applied to the usual Kac--Moody generators. Such generators have the form $B_i = F_i - c_iE_i K_i^{-1}$ in the coideal case. Let $T_{\omega_i}$ denote the usual Lusztig braid group operator, and $\mathbf{T}_{\omega_i}$ the corresponding Lu--Wang--Zhang braid group action operator. 

\Thm 
\label{intro: br gr}
Let $\g$ be of type $\mathsf{A}$. Then 
\[
\mathbf{T}_{\omega_i'}(B_i) \equiv T_{\omega_i'}(B_i) \mod U_{\neq i,+}. 
\]
Hence the diagram
\eq 
\begin{tikzcd}%[ row sep = 0.2cm]
\Oq \arrow[r] \arrow[d, hookrightarrow] & U_qL\mathfrak{sl}_2 \arrow[d, hookrightarrow]   \\
\mathcal{O}_q(\widehat{\mathfrak{sl}}_N) \arrow[r] & U_qL\mathfrak{sl}_N 
\end{tikzcd}
\eneq
commutes modulo $U_{\neq i,+}$. 
\enthm  

Here $U_{\neq i,+}$ is the subalgebra spanned by elements of positive $Q_i$-degree, where $Q_i$ is the root sublattice of $Q$ spanned by all positive roots except $\alpha_i$. 
Theorem \ref{intro: br gr} is proven in \S \ref{sec: comp br gr} (see Corollary \ref{thm: final br gr act compat}). For the details of the notation, see \S \ref{sec: qsp spl}. 
The generalization of Theorem \ref{intro: br gr} to $\mathsf{BCD}$ types can be found in \cite{LiPrz}.

\subsubsection{A rise in complexity} 
\label{sssub: complex}

Our proofs of Theorems \ref{intro main thm}--\ref{intro: br gr} are purely algebraic. However, they are significantly harder than the proofs of the corresponding results for quantum affine algebras. %For example, the proof of the coproduct formulae for the series $\pmb{\phi}_i^\pm(z)$ (mod $U_+$) relies on relatively straightforward inductions using the Drinfeld presentation. In contrast, the proof for~$\Thgsr_i(z)$, although also inductive in nature, is much more involved. 
The reason for this sharp increase in difficulty can be traced back to the complexity of the Drinfeld presentations themselves. There are two key types of relations which appear both in the original Drinfeld presentation of $U_q\widehat{\g}$ and the Lu--Wang presentation of $\gOq$. 
Let us call elements $x_{i,k}^\pm$ (resp.\ $A_{i,k}$) basic, and $h_{i,k}, \phi_{i,k}^\pm$ (resp.\ $H_{i,k}, \Theta_{i,k}$) diagonal. The first key relation is a commutation relation between basic and diagonal elements, producing further basic elements. In the case of $U_q\widehat{\g}$, the RHS of this relation \eqref{eq: hx rel} involves only one term, while in the relation \eqref{eq: rel2} for $\gOq$ two terms are present. The second key relation is a commutation relation between basic elements, producing diagonal elements. Again, this relation \eqref{eq: sl2 x+x- rel} is much easier for $U_q\widehat{\g}$, involving only a simple commutator on the LHS and typically one term on the RHS. In contrast, the corresponding relation \eqref{eq: rel3} for $\gOq$ involves two $q$-commutators on the LHS and typically two terms on the RHS. 

\subsubsection{Proof outline}

The key result, i.e., the factorization theorem (Theorem \ref{intro main thm}) is proven in three stages. The first stage involves rank one, i.e., the $q$-Onsager algebra, and so-called standard embeddings of the coideal subalgebra, i.e., with parameters $\mathbf{s} = 0$. The corresponding proof is carried out in \S \ref{sec: factorization}. The main idea is to decompose the image of the series $\pmb{A}_+(z)$ inside $U_q\widehat{\g}$ into homogeneous components with respect to the Drinfeld grading. It then suffices to compute the $\pm 1$ degree components. Relations \eqref{eq: rel2}--\eqref{eq: rel3} are used throughout. The second stage of the proof is a generalization to non-standard parameters $\mathbf{s} \neq 0$. Our proof relies on the fact that every non-standard coideal subalgebra can be obtained from a standard one via twisting with a character. Accordingly, it suffices to describe one-dimensional representations of the $q$-Onsager algebra and understand how the coproduct behaves with respect to such representations. This is carried out in \S \ref{sec: odr}-\ref{sec: mult property}. The final stage involves a generalization to quantum symmetric pairs of higher rank. The key role is played here by rank one subalgebras inside $\gOq$ isomorphic to the $q$-Onsager algebra. We study how these subalgebras interact with analogous rank one subalgebras in $U_q\widehat{\g}$. Using explicit extended affine Weyl group combinatorics, we establish the braid group action compatibility from Theorem~\ref{intro: br gr}, which enables us to derive the most general version of the factorization theorem. This is carried out in \S \ref{sec: qsp spl}-\ref{sec: comp br gr}. 

\subsection{Application to $q$-characters} 
\label{subsec: q char}

\nc{\Uh}{U_q(\widetilde{\mathfrak{h}})}
\nc{\Uhz}{\fext{U_q(\widetilde{\mathfrak{h}})}{z}}
\nc{\Rep}{\on{Rep}} 

The results of this paper are sufficient to formulate a notion of $q$-characters for quantum symmetric pairs and show that such characters are compatible, in an appropriate sense, with the usual $q$-characters for quantum affine algebras. A detailed exposition, for all the classical types, can be found in \cite{LiPrz}. Let us summarize here the main ideas. 

It is well known that, for quantum affine algebras, the notion of $q$-characters can be defined in at least three equivalent ways: via the universal $R$-matrix, via the spectrum of Drinfeld--Cartan operators, or via Nakajima's quiver varieties. In the quantum symmetric pairs case, it appears that the geometric and integrable systems methods are not yet advanced enough to yield an analogous theory of `boundary' $q$-characters. %In particular, the unavailability of a Levendorsky--Soibelman-type factorization of the affine universal $K$-matrix is a substantial obstacle. 

Therefore, we propose to define boundary $q$-characters directly via the Lu--Wang presentation. More precisely, we consider the generalized eigenspace decomposition of a finite dimensional representation with respect to the action of the operators $\Theta_{i,m}$, and let boundary $q$-characters encode the multiplicities of such eigenspaces. This is equivalent to taking the trace of a certain operator in the completion $\fext{\gOq \otimes \Uh}{z}$, yielding a map $\chi^\imath_q \colon \Rep \gOq \to \Uhz$. Applying Theorems \eqref{intro main thm}--\eqref{intro thm2}, one can  show that $\chi^\imath_q$ is compatible with the usual $q$-character map (\cite[Corollary 9.3]{LiPrz}). More precisely,  the diagram 
\[
\begin{tikzcd}[ row sep = 0.2cm]
\Rep U_qL\g \arrow[r, "\chi_q"] & \Z[Y_{i,a}^{\pm 1}]_{i \in I, a \in \C^{\times}}  \\
 \curvearrowright  & \curvearrowright  \\
\Rep \gOq \arrow[r, "\chi^\imath_q"] & \Uhz 
\end{tikzcd}
\]
commutes if $\Uhz$ is endowed with an appropriate `twisting' action of $\Z[Y_{i,a}^{\pm 1}]_{i \in I, a \in \C^{\times}}$. 
In particular, this result yields an easy way to compute the boundary $q$-characters of restriction representations.

\subsection{Other applications and future directions} 
\label{subsec: appl future dir}

We will now discuss various additional applications and insights stemming from our main results. 

\subsubsection{Affine Cartan subalgebra and representation theory} 
Our results suggest that the commutative subalgebra generated by the operators $\Theta_i$ should be viewed as a kind of ``affine Cartan subalgebra" in $\gOq$. Moreover, the eigenvalues of these operators are analogous to $\ell$-weights, and the associated eigenspace decomposition provides an analogue of a weight space decomposition. 

Nevertheless, the picture is far from clear. Let us mention the most important points. 
Firstly, since, unlike quantum groups, coideal subalgebras do not possess a triangular decomposition, finding the right replacement for the Cartan subalgebra is quite challenging. This problem was first addressed by Letzter \cite{Letzter-cartan} for quantum symmetric pairs of finite type. 
The notion of a Cartan subalgebra proposed therein was used by Wenzl \cite{Wenzl} to construct Verma modules for the coideal subalgebra quantizing $U(\mathfrak{so}_n)$. The previously mentioned lack of canonical raising and lowering operators makes this construction highly non-trivial. Remarkably, the resulting classification of irreducible modules \cite{Ior-Kl, Wenzl} involves so-called representations of non-classical type, which are not deformations of representations of $U(\mathfrak{so}_n)$. 

It would be interesting to integrate Wenzl's theory with our approach based on the affine Cartan subalgebra. However, this also appears rather non-trivial. In the case of quantum affine algebras, there is a natural compatibility between 
the subalgebra generated by the operators $\phi_{i,k}^\pm$ from the Drinfeld presentation (we will also refer to it as the ``affine Cartan") and the finite Cartan subalgebra, i.e., $\phi_{i,0}^\pm = K_i^{\pm1}$.
In contrast, in affine coideal subalgebras, $\Theta_{i,0}$ is merely a scalar and the higher $\Theta_{i,k}$'s do not appear to commute with the generators of Letzter's proposed Cartan subalgebra. 

A related problem, impeding the construction of a sensible highest weight theory for coideal subalgebras, is the lack of a  canonical ordering on the eigenvalues of the operators $\Theta_{i,k}$. Finally, we would like to remark that, on restriction representations, the eigenspace decomposition with respect to the affine Cartan subalgebra in $\gOq$ is, in general, quite different from (and hence incompatible with) the usual $U_q(\widehat{\g})$-weight decomposition with respect to the operators $K_i$. Put succinctly: using our results, it is possible to express the eigenvalues of the affine Cartan in $\gOq$ in terms of the eigenvalues of the affine Cartan in $U_q(\widehat{\g})$, but it is much harder to compute the corresponding eigenvectors. 

\subsubsection{The symplectic case} 

Remarkably, the situation in the case of the QSP coideal subalgebras of type $\mathsf{AII}$ (also known as symplectic twisted $q$-Yangians) is quite different. Classification results were obtained in \cite{molev-06, gow-molev},  using the RTT presentation. In particular, there exists a classification of irreducibles in terms of Drinfeld polynomials. 
The nicer behaviour of the symplectic case is, to a large extent, due to the fact that these coideal subalgebras contain some of the Cartan elements $K_i$. This facilitates the formulation of a highest weight theory compatible with restriction. 

Presently, a Drinfeld-type presentation is only available for aKM QSP of quasi-split type \cite{bas-kol-20, lu-wang-21, ZhangDr,  WangZhang2, WangZhang3}, corresponding to Satake diagrams with no black dots. %and in type $\mathsf{A}$ leads to coideal subalgebras quantizing $U(\widehat{\mathfrak{so}}_n)$. 
However, it is natural to expect that a Drinfeld presentation should exist for all aKM QSP, even though the technical complexity of the proofs is likely to increase with the complexity of the Satake diagrams. 
We also expect that, once the Drinfeld presentation becomes available in new cases, Theorem~\ref{intro main thm} and methods used to establish it will quite straightforwardly generalize to these new cases. In particular, in the symplectic case, corresponding to the Satake diagram with an alternating colouring, the analogue of our theorem should have more direct representation theoretic consequences, allowing us to recover the classification from \cite{gow-molev} from the point of view of the Drinfeld presentation. 
 
\subsubsection{Tridiagonal pairs} 

The representation theory of coideal subalgebras, outside of the symplectic case, is largely unknown. An exception is the split (i.e.\ orthogonal) case in rank one, which coincides with the $q$-Onsager algebra. This algebra has been extensively studied by Ito and Terwilliger \cite{ito-ter-08, ito-ter-10} using completely different combinatorial methods, involving association schemes and tridiagonal pairs. In particular, they formulated the notion of a Drinfeld polynomial for modules over the $q$-Onsager algebra which is, at least on the surface, completely different from ours. Moreover, they obtained an essentially complete classification of irreducible finite dimensional modules. 

We believe it will prove fruitful to further explore the interplay between Ito and Terwilliger's combinatorial approach and our approach based on the Drinfeld presentation. In particular, it is desirable to understand the relation between the different notions of Drinfeld polynomials, and to find an independent classification argument \'{a} la Chari and Pressley \cite{chari-pressley-qaa}. One may also expect that higher rank quantum symmetric pairs will lead to interesting generalizations of tridiagonal pairs and the associated combinatorics. 

\subsubsection{Drinfeld coproduct} 

According to Theorem \ref{intro: coproduct}, the series $\Thgsr_i(z)$ has a group-like property modulo $\fext{U_+}{z}$. 
We also prove (Proposition \ref{pro: twisted primitive}) that the series $\pmb{A}_i{(z)}$ exhibits a kind of ``twisted primitive'' behaviour. 
This is completely analogous to the behaviour of the series $\pmb{\phi}_i^{\pm}(z)$ and $\pmb{x}_i^{\pm}(z)$ in the quantum affine algebra. In the latter case, there exists another coproduct, called the (deformed) Drinfeld coproduct \cite{hernandez-05, hernandez-07}, with respect to which such group-like and twisted primitive properties hold exactly. 
Our results provide further evidence that the same should be true in the world of quantum symmetric pairs. In particular, the theory of quantum affinizations should admit a QSP generalization. 

\subsubsection{Factorization of the universal K-matrix} 

Originally, the notion of a $q$-character for quantum affine algebras was formulated in terms of the universal R-matrix \cite{FrenRes}. Its connection to spectra of the operators $\phi_{i,k}^{\pm}$ was established using an infinite product 
 factorization formula \cite{khor-tol, lev-sob-st, damiani-r}. 
In the QSP case, we are proceeding in the opposite direction, i.e., we propose to define $q$-characters in terms of the spectra of the operators $\Theta_{i,k}$. We expect that an equivalent formulation in terms of the universal K-matrix from \cite{AppelVlaar} is also possible. In particular, there should exist an analogue of the Levendorsky--Soibelman formula for the universal K-matrix. In the finite case, similar factorization results were already established in \cite{dobson-kolb-19, WangZhang}.

\subsubsection{Schur--Weyl duality and KLR algebras} 

It has been known for some time that quantum affine algebras and KLR algebras categorify the same objects. The extra grading on KLR algebras can be seen as explaining the various $q$-deformations appearing in the representation theory of quantum affine algebras, such as $q$-characters and deformations of Grothendieck rings. This link was made precise by Kang, Kashiwara and Kim \cite{kang-kashiwara-kim-18}, who constructed ``generalized Schur--Weyl duality'' functors from categories of modules over KLR algebras to those over quantum affine algebras. Their work was generalized to the QSP case in \cite{AppelPrz}. Accordingly, we expect that $q$-characters for aKM QSP should also have an interesting interpretation in terms of \emph{orientifold} KLR algebras and the shuffle modules from \cite{przez-oklr}. 

\subsection{Organization}

Sections \ref{sec: q-ons} and \ref{sec: qsp spl} consist mainly of recollections of existing material. In Section \ref{sec: rat}, we establish basic rationality properties of generating series, as well as the important ``$C$-symmetry'' property. Section \ref{sec: factorization} is the technical core of the paper, and contains the proof of the first version of the factorization theorem (Theorem \ref{intro main thm}) - the rank one case with parameters $\mathbf{s} = 0$. In Section \ref{sec: drf}, we derive various consequences for spectra of the operators $\Theta_k$. In particular, we formulate the notion of Drinfeld rational fractions and prove Theorem \ref{intro thm2}. Sections \ref{sec: odr} and \ref{sec: mult property} extend the factorization theorem to the case of non-standard parameters $\mathbf{s} \neq 0$. We also prove the coproduct theorem (Theorem \ref{intro: coproduct}). In Section \ref{sec: comp br gr}, we prove Theorem \ref{intro: br gr} and generalize the factorization theorem to higher ranks. 

\addtocontents{toc}{\SkipTocEntry}
	
\section*{Acknowledgements} 
I would like to express my gratitude for all the support I have received. I am particularly indebted to Andrea Appel, who has pointed out to me the importance of rationality properties and one dimensional representations, as well as to Jian-Rong Li for the careful reading of the manuscript and our ongoing collaboration on a related project. 
I would also like to thank Paul Terwilliger for insightful discussions about the $q$-Onsager algebra, Alexander Molev for correspondence about twisted Yangians, Daniel Orr and Mark Shimozono for providing a reference on the combinatorics of the extended affine Weyl group, as well as Pascal Baseilhac, Weiqiang Wang and Weinan Zhang for comments on the first version of this paper. 

\section{The $q$-Onsager algebra} 
\label{sec: q-ons}

Throughout this paper, we assume that $q \in \C^\times$ is not a root of unity. We use the standard notation for quantum numbers, i.e., 
\[
[k] = \frac{q^k - q^{-k}}{q-q^{-1}}, \qquad [k]! = [k][k-1]\cdots[1], \qquad \qbinom{k}{l} = \frac{[k]!}{[k-l]!l!}, 
\]
for $k,l \in \Z_+$, and $[-k] = -[k]$. We also set $[a,b]_q = ab - qba$, as usual.

\subsection{Quantum affine $\mathfrak{sl}_2$}

The quantum affine algebra $U_q\widehat{\mathfrak{sl}}_2$ is the associative algebra over $\kor$ with generators $\Eg_i^\pm, \Kg_i^{\pm 1}$, where $i \in \{0,1\}$, and relations:
\begin{alignat*}{3}
\Kg_i\Kg_i^{-1} &= \Kg_i^{-1}\Kg_i = 1, \\ 
\Kg_0\Kg_1 &= \Kg_1 \Kg_0, \\
\Kg_i\Eg_i^{\pm} &= q^{\pm2} \Eg_i^\pm\Kg_i, \\
\Kg_i\Eg_j^{\pm} &= q^{\mp2} \Eg_j^\pm\Kg_i  && \quad (i \neq j), \\
[\Eg_i^{\pm}, \Eg_j^{\mp}] &= \pm\delta_{ij} \frac{\Kg_i - \Kg_i^{-1}}{q - q^{-1}}, \\
\Serre(\Eg_i^{\pm}, \Eg_j^\pm) &= 0 && \quad (i \neq j),
\end{alignat*}
where $\Serre(x,y)$ denotes the usual $q$-Serre relations, i.e., 
\[
\Serre(x,y) = \sum_{r=0}^{3} (-1)^r \qbinom{3}{r} x^{3-r}yx^r. 
\] 
The algebra $U_q\widehat{\mathfrak{sl}}_2$ is a Hopf algebra. We will use the coproduct
\eq \label{eq: coprod on U}
\Delta(\Eg_i^+) = \Eg_i^+ \otimes 1 + \Kg_i \otimes \Eg_i^+,  \quad \Delta(\Eg_i^-) = \Eg_i^- \otimes \Kg_i^{-1} + 1 \otimes \Eg_i^-, \quad \Delta(\Kg_i^{\pm 1}) = \Kg_i^{\pm 1} \otimes \Kg_i^{\pm 1}. 
\eneq
This choice of coproduct is consistent with \cite{kolb-14}, and differs from the coproduct in \cite{chari-pressley-qaa} by a flip. The counit is given by 
\[
\varepsilon(\Eg_i^{\pm}) = 0, \quad \varepsilon(\Kg_i^{\pm 1}) = 1. 
\]

\subsection{The quantum loop algebra} 
The quantum loop algebra $\uLsl$ is the quotient of $U_q\widehat{\mathfrak{sl}}_2$ by the relation 
\eq
\Kg_1 = \Kg = \Kg_0\mi. 
\eneq 
Let us recall the Drinfeld presentation of the quantum loop algebra $\uLsl$. By \cite{drinfeld-dp, beck-94}, $\uLsl$ is isomorphic to the algebra generated by $\xgpm_k, \hg_l, \Kg^{\pm 1}$, where $k \in \Z, \ l \in \Z - \{0\}$, subject to the following relations:
\begin{align}
\Kg\Kg^{-1} =& \ \Kg^{-1}\Kg = 1, \\
[\hg_k, \hg_l] =& \ 0, \\
\Kg\hg_k =& \ \hg_k\Kg, \\ 
\Kg \xgpm_k =& \ q^{\pm 2}\xgpm_k\Kg, \label{eq: Kx rel} \\
[\hg_k,\xgpm_l] =& \pm \textstyle \frac{[2k]}{k}\xgpm_{k+l} \label{eq: hx rel},\\
\label{eq: sl2 4x rel} \xgpm_{k+1}\xgpm_l - q^{\pm 2} \xgpm_l\xgpm_{k+1} =& \ q^{\pm2}\xgpm_k\xgpm_{l+1} - \xgpm_{l+1}\xgpm_k, \\
\label{eq: sl2 x+x- rel}
[\xgp_k, \xgm_l] =& \ \textstyle\frac{1}{q-q^{-1}}(\psig_{k+l} - \phig_{k+l}),
\end{align}
where 
\begin{align}
\Psig(z) = \sum_{k=0}^\infty \psig_k z^k =& \ \Kg \exp\left( (q-q^{-1}) \sum_{k=1}^\infty \hg_k z^k \right), \\ 
\Phig(z) = \sum_{k=0}^\infty \phig_{-k} z^{-k} =& \ \Kg^{-1} \exp\left( -(q-q^{-1}) \sum_{k=1}^\infty \hg_{-k} z^{-k} \right).  
\end{align}

For later convenience, we note that, in particular, 
\eq
\psig_0 = \Kg, \quad \phig_0 = \Kg\mi, \quad \psig_1 = (q-q\mi)\Kg\hg_1, \quad \phig_{-1} = - (q-q\mi)\Kg\mi\hg_{-1}. 
\eneq
We also write
\eq
\Xgp(z) = \sum_{k=-\infty}^\infty \xgp_k z^k, \quad
\Xgm(z) = \sum_{k=-\infty}^\infty \xgm_k z^k. 
\eneq

The Serre and Drinfeld presentations are related by the following isomorphism (\cite[Proposition 3.7]{beck-94}):
\[
\Kg \mapsto \Kg, \quad \Eg_1^\pm \mapsto \xgpm_0, \quad \Eg_0^+ \mapsto -\Kg^{-1}\xgm_1, \quad \Eg_0^- \mapsto -\xgp_{-1}\Kg.   
\]

We consider $\uLsl$ as a graded algebra with 
\eq \label{eq: qOns Dr gr}
\deg \xgpm_k = \pm1, \quad \deg h_k = 0. 
\eneq
Given $Y \in \uLsl$, let $Y\gup{n}$ denote its homogeneous component of degree $n$. 
For $i \geq 1$, let $U_{\geq i}$ be the subalgebra of $\uLsl$ spanned by elements of degree at least $i$. If $i=1$, we write $U_+ = U_{\geq 1}$. 

\subsection{The $q$-Onsager algebra}

The defining relations of the $q$-Onsager algebra, known as the $q$-Dolan--Grady relations, were first considered by Terwilliger \cite{Ter1}. The $q$-Onsager algebra itself was first defined in \cite{Ter2}. In mathematical physics, it was also studied by Baseilhac in \cite{Bas1}. 
In the framework of \cite{kolb-14}, 
the $q$-Onsager algebra is a quantum symmetric pair coideal subalgebra of $U_q\widehat{\mathfrak{sl}}_2$,  corresponding to the Satake diagram of type $A_1^{(1)}$ with no black vertices and the trivial automorphism of the underlying affine Dynkin diagram. It is also commonly referred to as the split affine $\imath$quantum group of rank one (e.g., \cite{lu-wang-21}). 

As in \cite{lu-wang-21}, we consider a central extension of the $q$-Onsager algebra called the \emph{universal} $q$-Onsager algebra. 
It has the following definition by generators and relations. 

\Defi
The universal $q$-Onsager algebra $\Oq$ is the algebra generated by $B_0, B_1$ and invertible central elements $\ccc_0, \ccc_1$ subject to relations:
\eq
\label{eq: Ons rel}
\sum_{r=0}^3 (-1)^r \sqbin{3}{r} \Bg_i^{3-r}\Bg_j\Bg_i^r = - q \ccc_i [2]^2[\Bg_i, \Bg_j] \qquad (0 \leq i \neq j \leq 1).
\eneq
\enDefi

It follows from \cite[Theorem 8.3]{kolb-14} that $Z(\Oq) = \C[\ccc_0^{\pm1}, \ccc_1^{\pm 1}]$.  
Given $\mathbf{c} = (c_0, c_1) \in (\C^\times)^{2}$, we have a central character 
\eq
\label{eq: central character} 
\chi_{\mathbf{c}} \colon Z(\Oq) \to \C^\times, \qquad \ccc_i \mapsto c_i. 
\eneq
Let $\Oqc$ be the corresponding central reduction, i.e., the quotient of $\Oq$ by the ideal generated by $\ccc_i - c_i$. 
We make no distinction in the notation between elements of $\Oq$ and their images in $\Oqc$. 
It is known that algebras $\Oqc$ and $\mathcal{O}_q^{\mathbf{c}'}$ are isomorphic for any choices of parameters $\mathbf{c}, \mathbf{c}'$ (see, e.g., \cite[Proposition 9.2]{kolb-14} or \cite[Lemma 2.5.1]{Wat21}). 

\Rem 
Let us briefly comment on our choice of conventions. 
In $\Oqc$, 
the relations \eqref{eq: Ons rel} specialize to the defining relations of \cite[Example 7.6]{kolb-14}. 
On the other hand, the dictionary between our conventions and those of \cite{lu-wang-21} is given by setting 
\[
v = q, \qquad \mathbb{K}_i = q^2 \ccc_i.
\] 
\enrem

\subsection{The Lu--Wang presentation}

Recently, a Drinfeld-type presentation of $\Oq$ was found by Lu and Wang \cite{lu-wang-21}. 

\Defi 
Let $\DrOq$ be the algebra generated by $\Hg_m$ and $\Ag_r$, where $m\geq1$, $r\in\Z$, and invertible central elements $\ccc_1, \CCC$, subject to the following relations:
\begin{align}
[\Hg_m,\Hg_n] &=0, \label{eq: rel1} \\
[\Hg_m, \Ag_{r}] &= \textstyle \frac{[2m]}{m} (\Ag_{r+m}- \Ag_{r-m}\CCC^m), \label{eq: rel2} \\
\label{eq: rel3}
[\Ag_{r}, \Ag_{s+1}]_{q^{-2}}  -q^{-2} [\Ag_{r+1}, \Ag_{s}]_{q^{2}}
&= \ccc_1\CCC^r\Thg_{s-r+1} - q^{-2}\ccc_1\CCC^{r+1}\Thg_{s-r-1}  \\ \notag
&\quad  + \ccc_1\CCC^s\Thg_{r-s+1}  -q^{-2}\ccc_1 \CCC^{s+1}\Thg_{r-s-1}, 
\end{align}
where $m,n\geq1$; $r,s\in \Z$; 
and 
\[
1+ \sum_{m=1}^\infty (q-q^{-1})\Thg_{m} z^m  =  \exp \left( (q-q^{-1}) \sum_{m=1}^\infty \Hg_m z^m \right).
\]
In particular, $\Thg_1 = \Hg_1$. 
By convention, $\Thg_0 = (q-q^{-1})^{-1}$ and $\Thg_m = 0$ for $m\leq -1$. 
\enDefi 

By \cite[Theorem 2.16]{lu-wang-21}, there is an algebra isomorphism 
\eq
\label{eq: droq->oq}
\DrOq \isoto{} \Oq, \quad \quad \Ag_0 \mapsto \Bg_1, \quad \Ag_{-1} \mapsto q^{-2}\ccc_0^{-1}\Bg_0, 
\quad \ccc_1 \mapsto \ccc_1, \quad \CCC \mapsto q^4\ccc_0\ccc_1. 
\eneq

\nc{\UqLsl}[1]{U_qL\mathfrak{sl}_{#1}}

\subsection{Coideal structures} \label{subsec: coideal str sl2}

Following \cite[Example 7.6]{kolb-14}, 
for any $\mathbf{s} = (\sss_0, \sss_1) \in \kor^2$, there 
is an injective algebra homomorphism 
\eq
\label{eq: Kolb emb}
\eta_{\mathbf{c},\mathbf{s}} \colon 
\Oqc \monoto \uLsl, \quad \Bg_i \mapsto \Eg_i^- - \cc_i \Eg_i^+ \Kg_i^{-1} + \sss_i \Kg_i^{-1}, %\quad \ccc_i \mapsto \cc_i \cdot 1 \qquad (i=0,1),
\eneq
giving $\Oqc$ the structure of a right coideal subalgebra of $\uLsl$. 
Let 
\[
\tops = \Delta \circ \etacs \colon \ \Oqc \to \Oqc \otimes \uLsl 
\]
be the associated coproduct. 
In terms of the Drinfeld presentation, the map \eqref{eq: Kolb emb} is given by: 
\eq 
\label{eq: droq->sl}
\Ag_0 \mapsto \xgm_0 - \cc_1 q^2 \Kg^{-1}\xgp_0 + \sss_1 \Kg^{-1}, \quad \Ag_{-1} \mapsto  -q^{-4}\cc_0^{-1} \Kg\xgp_{-1} + \xgm_1 +  q^{-2}\cc_0^{-1}\sss_0 \Kg. 
\eneq

%restriction functor
The choice of a coideal subalgebra structure, or, equivalently, the embedding $\eta_{\mathbf{c},\mathbf{s}}$, gives rise to a right monoidal action 
\eq \label{eq: monoidal action}
\mac{\mathbf{c},\mathbf{s}}\colon \Modfdc(\Oq) \boxtimes \Modfd(\UqLsl{2}) \to \Modfdc(\Oq),
\eneq
where 
$\Modfdc(\Oq)$ is the category of finite dimensional $\Oq$-modules with a fixed central character \eqref{eq: central character}. If $\mathbf{s} = (0,0)$, we call the monoidal action 
\emph{standard}. 

\subsection{One-dimensional modules} 

Restricting the trivial representation of $\uLsl$ via $\eta_{\mathbf{c},\mathbf{s}}$ gives rise to a continuous family of one-dimensional $\Oqc$-modules $\kor_{\mathbf{c},\mathbf{s}} = \eta_{\mathbf{c},\mathbf{s}}^*(\kor)$, parametrised by $\mathbf{s} \in \C^2$. 
These modules correspond to the characters 
\eq \label{eq: xi s char}
\xi_{\mathbf{c},\mathbf{s}} = \varepsilon\circ\eta_{\mathbf{c},\mathbf{s}}\colon \ \Oqc \ \to \ \kor,
\eneq
explicitly given by 
\[
\xi_{\mathbf{c},\mathbf{s}}(A_0)=s_1, \quad \xi_{\mathbf{c},\mathbf{s}}(A_{-1}) = q^{-2} \cc_0\mi s_0. 
\]

Let $\Oqp{\mathbf{c},\mathbf{s}} = \eta_{\mathbf{c},\mathbf{s}}(\Oqc)$. 
By \cite[Sec.~3.5]{dobson-kolb-19}, the character $\xi_{\mathbf{c},\mathbf{s}}$ induces an algebra isomorphism 
\[
\phi_{\mathbf{c},\mathbf{s}}:\Oqp{\mathbf{c},0} \to \Oqp{\mathbf{c},\mathbf{s}}
\]
such that
$\eta_{\mathbf{c},\mathbf{s}} = \phi_{\mathbf{c},\mathbf{s}} \circ \eta_{\mathbf{c},0}$. Explicitly, $\phi_{\mathbf{c},\mathbf{s}}$ is
given by $\phi_{\mathbf{c},\mathbf{s}} = (\wt{\xi}_{\mathbf{c},\mathbf{s}}\ten\id) \circ \Delta|_{\Oqp{\mathbf{c},0}}$, where
$\wt{\xi}_{\mathbf{c},\mathbf{s}} = \xi_{\mathbf{c},\mathbf{s}} \circ (\eta'_{\mathbf{c},0})^{-1}$, and $\eta'_{\mathbf{c},0}$ denotes $\eta_{\mathbf{c},0}$ with codomain restricted to its image $\Oqp{\mathbf{c},0}$. \\

In the following, we will extensively study $\Oq$-modules which arise from restricting $\uLsl$-modules via $\eta_{\mathbf{c},\mathbf{s}}$. The following lemma shows that any such module can be realized via a \emph{standard} action on a one-dimensional representation. 

\begin{lem} \label{rem: s-reduction}
Let $V \in \Modfd(\UqLsl{2})$. Then, there is an isomorphism of $\Oqc$-modules
\[ 
\eta_{\mathbf{c},\mathbf{s}}^*(V) \cong \kor_{\mathbf{c},\mathbf{s}}\mac{\mathbf{c},0} V, 
\]
given by the canonical identification of vector spaces $V \cong \kor\ten V$.
\end{lem}

\begin{proof}
It is enough to observe that $\eta_{\mathbf{c},\mathbf{s}}=\phi_{\cs}\circ\eta_{\mathbf{c},0}$, \ie
\eq \label{eq: eta cs delta}
\eta_{\mathbf{c},\mathbf{s}}(Y) = (\wt{\xi}_{\mathbf{c},\mathbf{s}} \ten \id) \circ \Delta_{\mathbf{c},0}(Y)
\eneq
for any $Y \in\Oqc$.
\end{proof}

\subsection{$\CC$-Duality} \label{subsec: C-duality}

Let $\Modfdw(\Oq)$ be the full subcategory of $\Oq$-modules $M$ such that: 
\begin{enumerate}[label=(\alph*)]
\item $M$ is finitely generated, 
\item $M$ has central character $\chi_{\mathbf{c}}$, 
\item there is a joint generalized eigenspace decomposition
\eq \label{eq: gen eig dec cat}
M = \bigoplus_{\bm\lambda = (\lambda_s)_{s\in \Z_{\geq 0}}} M_{\bm\lambda}, \qquad M_{\bm\lambda} = \{ v \in M \mid \exists p \ \forall s \in \Z_{\geq 0}: (\Thg_{s} - \lambda_{s})^p \cdot v = 0 \}, 
\eneq
with $M_{\bm\lambda}$ finite dimensional. 
\item For each $\bm\lambda$, the vector space generated by $A_r \cdot M_{\bm\lambda}$ $(r \in \Z)$ is finite-dimensional. 
\end{enumerate}
If $M_{\bm\lambda} \neq \{0\}$, we call $\bm\lambda$ a \emph{weight} of $M$. In the context of quantum affine algebras, this is usually referred to as an $\ell$-weight or a pseudo-weight (see, \emph{e.g.}, \cite{hernandez-05, her-jim-12}). 
If $\bm \lambda$ and $\bm \mu$ are weights, we write $\bm\lambda \to \bm \mu$ if there exists $v \in M_{\bm\lambda}$ and $r\in \Z$ such that $A_r\cdot v$ has a non-trivial projection onto $M_{\bm \mu}$.

\Lem
There is an algebra anti-automorphism $\tau \colon \Oq \to \Oq$ sending
\[
\Ag_r \mapsto \CCC^r \Ag_{-r}, \quad \Hg_m \mapsto \Hg_m, \quad \ccc_i \mapsto \ccc_i, 
\]
inducing a duality
\[
\Modfdw(\Oq) \to \Modfdw(\Oq), \quad M \mapsto M^*, 
\]
where $M^*$ is the restricted dual with respect to the decomposition \eqref{eq: gen eig dec cat}. 
\enlem

\Proof
One directly checks that the defining relations are preserved up to switching the order of multiplication. This holds for \eqref{eq: rel2} because
\begin{align}
\textstyle\frac{[2m]}{m}\tau( \Ag_{r+m}-\Ag_{r-m}\CCC^m) =& \ 
\textstyle \frac{[2m]}{m}\CCC^r (\Ag_{-r-m}\CCC^{m} - \Ag_{-r+m}) \\
 =& \ 
-[\Hg_m, \CCC^r\Ag_{-r}] = \tau([\Ag_{r}, \Hg_m]).  
\end{align} 
To see that \eqref{eq: rel3} is also preserved, observe that switching the order of multiplication on the LHS amounts to swapping $r$ and $s$, and the RHS is invariant under such a swap.
\enproof

We will refer to $M^*$ as the $\CC$\emph{-dual} of $M$. 

\section{Rationality}
\label{sec: rat}

\nc{\Asp}{\mathbf{A}_+}
\nc{\Asm}{\mathbf{A}_-}
\nc{\Aspm}{\pmb{A}_\pm}
\nc{\Thr}{{\bm\vartheta}}
\nc{\Ar}{\mathcal{A}}
\nc{\RepfdOq}{\operatorname{Rep_{fd}}(\Oq)}
\nc{\RepfdUqLsl}[1]{\operatorname{Rep_{fd}}(\UqLsl{#1})}

The rationality property of the generating series of quantum loop algebras 
is due to Beck–Kac \cite{kac-beck-96} and Hernandez \cite{hernandez-07} 
(see also \cite{gautam-tl-16}). 	
Below we establish an analogous property for the generating series of the $q$-Onsager algebra. 
Set 
\eq
\As = \sum_{r= - \infty}^\infty \Ag_r z^r, \quad  
\Aps = \sum_{r \geq 0} \Ag_r z^r, \quad \Ams = - \sum_{r \leq -1} \Ag_r z^r. 
\eneq
\eq \label{eq: Thg acute def}
\Thgs(z) = \sum_{r \geq 0} \Thg_r z^r, \quad \Thgsa = \frac{1-q^{-2}\CCC z^2}{1-\CCC z^2} \Thgs(z). 
\eneq
For convenience, we also abbreviate
\[
\bHg_1 = [2]^{-1} \Hg_1, \qquad \bhg_{\pm1} = [2]^{-1} \hg_{\pm1}, \qquad \ad_{\bHg_1}(-) = [\bHg_1, -]. 
\] 
\begin{pro} \label{pro: rationality main pro}
Let $V$ be a module in $\Modfdw(\Oq)$. 
Set 
\[
\CC = q^4\cc_0 \cc_1. 
\]
Then, for every weight $\bm\lambda$ of $V$: 
\be
\item The generating series
	$
	\Aspm{(z)}\in\fext{\Hom(V_{\bm\lambda}, \bigoplus_{\bm\lambda \to \bm\mu}V_{\bm\mu})}{z^{\pm1}}
	$
	are the expansions of the \emph{same} rational function\footnote{In the later sections we will typically use the same notation for rational functions and their series expansions where no confusion can arise.} 
	\begin{align}\label{eq:Ar} 
		\Ar(z)=(1-\ad_{\bHg_1}z-\CC z^2)^{-1}(A_0+ \CC z A_{-1})
	\end{align}
	at $z=0$ and $z=\infty$, respectively. 
\item Moreover, the series 
	$
	\Thgsa\in\fext{\End(V_{\bm\lambda})}{z} 
	$
	is the expansion of the rational function
	\begin{align}
	\label{eq:Thr}
	\Thr(z)=\frac{(1-\CC z^2)\Thg_0+\cc_1^{-1}\CC z^2\left(z\mi[A_{-1}, \Ar(z)]_{q^{-2}}-q^{-2}[A_0,\Ar(z)]_{q^2}\right)}{1-\CC z^2}
	\end{align}
	at $z=0$. 
\item The rational function $\Thr(z)$ satisfies the “$\CC$\emph{-symmetry}" property:
\eq \label{eq: C-symmetry}
\Thr(z)=\Thr(\CC\mi z^{-1}).
\eneq
\ee
\end{pro}

\Proof
\emph{(1)}
Relation \eqref{eq: rel2}, with $m=1$, implies that 
\begin{align}
[\bHg_1, \Aspm(z)] =& \ z\mi \Aspm(z) - z\CC \Aspm(z) -z\mi \Ag_{0} - \CC \Ag_{-1}.
%-[\bHg_1, \Ams] =& \ z\mi \Ams - z\CC \Ams + z\mi\Ag_{0} + \CC \Ag_{-1}
\end{align}
Multiplying both sides by $z$ and rearranging the terms, we get 
\eq
\label{eq: Aplus series}
(1-\ad_{\bHg_1}z-\CC z^2)\Aspm(z) = \Ag_0 + \CC \Ag_{-1}z.
\eneq

\emph{(2)}
Setting $r=-1$ in \eqref{eq: rel3}, we get, for $s \geq 0$, 
\eq \label{eq: As vs Ths eq}
[\Ag_{-1}, \Ag_{s+1}]_{q^{-2}}  -q^{-2} [\Ag_{0}, \Ag_{s}]_{q^{2}}
= \cc_1(C^{-1}\Thg_{s+2} - q^{-2}\Thg_{s} + \delta_{s=0}\Thg_0). 
\eneq
Also note that, setting $r = s = -1$, relation \eqref{eq: rel3} simplifies to
\eq \label{eq: H1 simple}
[\Ag_{-1},\Ag_0]_{q^{-2}} = q^{-4}\cc_0^{-1}\Hg_1 = \cc_1\CC\mi  \Thg_1. 
\eneq
Hence, summing \eqref{eq: As vs Ths eq} over all $s \geq 0$, we get 
\begin{align}
-z \cc_1\CC^{-1}\Thg_1 \ +& \ z^2\left(z\mi[\Ag_{-1}, \Aps]_{q^{-2}} 
- q^{-2} [\Ag_{0}, \Aps]_{q^{2}}\right) = \\ 
=& \ \cc_1\left((C^{-1}- q^{-2}z^2)\Thgs(z) 
 - \CC\mi \Thg_1 z  + (z^2 - \CC\mi)\Thg_0\right).
\end{align}
This simplifies to
\begin{align} \label{eq: Theta series in A terms}
(1 - q^{-2}\CC z^2)\Thgs(z) =& \ (1-\CC z^2)\Thg_0 \\
& +  z^2\cc_1\mi\CC\left(z\mi[\Ag_{-1}, \Aps]_{q^{-2}}  -q^{-2} [\Ag_{0}, \Aps]_{q^{2}}\right). 
\end{align}
Dividing both sides by $(1-\CC z^2)$, we get \eqref{eq:Thr}. 

\emph{(3)} On the other hand, for $s \leq -1$ and $r=-1$, relation \eqref{eq: rel3} simplifies to
\eq
\label{eq: As vs Ths eq neg}
[\Ag_{-1}, \Ag_{s+1}]_{q^{-2}}  -q^{-2} [\Ag_{0}, \Ag_{s}]_{q^{2}}
= \cc_1(C^{s}\Thg_{-s} - q^{-2}\CC^{s+1}\Thg_{-s-2} +\CC\mi( \delta_{s=-1}\Thg_1 + \delta_{s=-2}\Thg_0)). 
\eneq
Summing \eqref{eq: As vs Ths eq} over all $s \leq -1$, we get 
\begin{align}
z\mi \cc_1\CC^{-1}\Thg_1 \ +& \ z\mi[\Ag_{-1}, \Ams]_{q^{-2}} 
- q^{-2} [\Ag_{0}, \Ams]_{q^{2}} = \\ 
=& \ \cc_1\left((1- q^{-2}\CC\mi z^{-2})\Thgs(\CC\mi z\mi) 
 +\CC\mi \Thg_1 z\mi  + (z^{-2}\CC\mi - 1)\Thg_0\right).
\end{align}
This simplifies to
\eq
(1 - q^{-2}\CC\mi z^{-2})\Thgs(\CC\mi z\mi) = (1-\CC\mi z^{-2})\Thg_0 + \cc_1\mi\left( z\mi[\Ag_{-1}, \Ams]_{q^{-2}}  -q^{-2} [\Ag_{0}, \Ams]_{q^{2}}\right). 
\eneq
Dividing both sides by $(1-\CC\mi z^{-2})$, we get 
\begin{align}
\Thr(\CC\mi z\mi) =& \ \frac{(1-\CC\mi z^{-2})\Thg_0 + \cc_1\mi\left( z\mi[\Ag_{-1}, \Ams]_{q^{-2}}  -q^{-2} [\Ag_{0}, \Ams]_{q^{2}}\right)}{1-\CC\mi z^{-2}} \\
=& \ \frac{(1-\CC z^{2})\Thg_0 - \cc_1\mi\CC z^2\left( z\mi[\Ag_{-1}, \Ams]_{q^{-2}}  -q^{-2} [\Ag_{0}, \Ams]_{q^{2}}\right)}{1-\CC z^{2}}. 
\end{align} 
Since $\Aps$ and $-\Ams$ are both expansions of the same rational function $\Ar(z)$, 
the identity $\Thr(z)=\Thr(\CC\mi z^{-1})$ follows.
\enproof

\begin{rem} \label{rem: Theta vs acute Theta}
The $\CC$-symmetry property does not hold for $\Thgs(z)$. This is the first indication that the generators $\acute{\Thg}_s$ are more natural from a representation theoretic point of view. We will see further evidence for this later (c.f.\ Remark \ref{rem: Thg grave more natural}).  
\end{rem}

%%%%%%%%%%%%%%%%%%%%%%%%%%%%%%%%
%%%%%%%%%%%%%%%%%%%%%%%%%%%%%%%%

\section{Factorization of $\Thgs(z)$} \label{sec: factorization}

The goal of this section is to prove our first main result, namely a factorization theorem for $\Thgs(z)$ (Theorem \ref{thm: fact}), and interpret it in terms of Drinfeld rational functions (Theorem \ref{cor: Drinfeld rational functions}). The factorization theorem establishes a connection between the Drinfeld presentation of the $q$-Onsager algebra and the Drinfeld presentation of $\uLsl$. It is key to understanding restricted representations. 

More precisely, we are able to express the image of the series $\Thgs(z)$ in $\uLsl$ as a product of shifted series $\pmb{\phi}^+(z)$ and $\pmb{\phi}^-(z)$, modulo $U_+$. The proof relies on an explicit computation of the homogeneous components of low degree of the images of $\Oq$-generators.  
The main step in the argument is the calculation of $\Apsone$ in \S \ref{sec: A+1 calc}. 

For the entirety of \S \ref{sec: factorization}, we fix a central character \eqref{eq: central character} and the   
embedding \eqref{eq: Kolb emb}, assuming that 
\eq \label{s-assumption}
\sss_0 = 0 = \sss_1.  
\eneq
This assumption is justified by Lemma \ref{rem: s-reduction}, and will be lifted in \S \ref{sec: mult property}. For the sake of readability, we make no distinction in the notation between an element of $\Oq$ and its image under \eqref{eq: Kolb emb}. In other words, throughout \S \ref{sec: factorization}, $\Ag_i$ stands for $\eta_{\mathbf{c},0}(\Ag_i)$, etc. 

\subsection{Auxiliary relations} 

We begin by deducing several auxiliary relations in $\uLsl$. 
By \cite[(1.2.7)]{Nakajima-01}, the following relation holds:
\eq
(z - q^{-2}w)\Phig(z)\Xgp(w) = (q^{-2}z - w)\Xgp(w)\Phig(z). 
\eneq
Collecting homogeneous terms, we get
\eq
\label{eq: Nakajima hom}
\phig_{-r}\xgp_s + \xgp_{s-1}\phig_{-r+1} = q^{-2}(\xgp_s\phig_{-r} + \phig_{-r+1}\xgp_{s-1})
\eneq
for $r \geq 1$. 

\Lem
\label{lem: xphi}
We have
\eq
\xgp_k \phig_{-r} = q^2 \phig_{-r} \xgp_k + (q^4 - 1) \sum_{s = 1}^r q^{2(s-1)} \phig_{s-r} \xgp_{k-s}. 
\eneq
\enlem

\Proof
The proof is by induction on $r$. The base case $r=0$ is just the relation \eqref{eq: Kx rel}. For $r >0$, we get 
\begin{align}
\xgp_k \phig_{-r} =& \ q^2 \phig_{-r} \xgp_k - \phig_{-r+1}\xgp_{k-1} + q^2 \xgp_{k-1}\phig_{-r+1} \\
 =& \ q^2 \phig_{-r} \xgp_k - \phig_{-r+1}\xgp_{k-1} + q^2 \left(q^2 \phig_{-r+1} \xgp_{k-1} + (q^4 - 1) \sum_{s = 1}^{r-1} q^{2(s-1)} \phig_{s-r+1} \xgp_{k-s-1}\right) \\
=& \ q^2 \phig_{-r} \xgp_k + (q^4 - 1) \sum_{s = 1}^r q^{2(s-1)} \phig_{s-r} \xgp_{k-s}, 
\end{align}
where the first equality follows from \eqref{eq: Nakajima hom} and the second by induction. 
\enproof

\Lem \label{lem: xgp phi comm}
We have
\[
\xgp_{k}\Phig(z\mi) = q^2 \Phig(z\mi) \xgp_{k} + %(q^2 - q^{-2})(-q^4)^{k}z^k \Phig(-q^{-2}z\mi)\Xgp_{\leq k-1}(-q^{-4}z^{-1})
(q^4 - 1)z \Phig(z\mi) (1-q^2 \ada z)\mi \xgp_{k-1}. 
\]
Hence
\[
[\xgp_{k}, \Phig(z\mi)] = \Phig(z\mi) (1-q^2 \ada z)\mi (q^2 - 1) (1 + \ada z) \xgp_k. 
\]
\enlem

\Proof 
By Lemma \ref{lem: xphi},  
we have 
\begin{align}
\xgp_k \Phig(z\mi) =& \ q^2 \Phig(z\mi) \xgp_{k} + (q^4-1) \sum_{r \geq 0} \sum_{s=1}^r q^{2(s-1)} \phig_{s-r} \xgp_{k-s} z^r \\
=& \ q^2 \Phig(z\mi) \xgp_{k} + (q^2-q^{-2})q^{2k}z^k \sum_{r \geq 0} \sum_{s=1}^r \phig_{s-r} q^{2(s-k)} \xgp_{k-s} z^{r-k} \\ 
=& \ q^2 \Phig(z\mi) \xgp_{k}  + (q^2 - q^{-2})q^{2k}z^k \Phig(z\mi)\Xgp_{\leq k-1}(q^{-2}z^{-1}) \\
=& \ q^2 \Phig(z\mi) \xgp_{k} + 
(q^4 - 1)z \Phig(z\mi) (1-q^2 \ada z)\mi \xgp_{k-1}.  
%=& \ \Phig(z\mi) (1+q^4 \ada z)\mi \big( 
\end{align}
Therefore, 
\begin{align}
[\xgp_{k}, \Phig(z\mi)] =& \  (q^2-1) \Phig(z\mi) \xgp_{k} + (q^4 - 1)z \Phig(z\mi) (1-q^2 \ada z)\mi \xgp_{k-1} \\
=& \ \Phig(z\mi) (1-q^2 \ada z)\mi (q^2 - 1) (\xgp_k + \xgp_{k-1} z ) \\
=& \ \Phig(z\mi) (1-q^2 \ada z)\mi (q^2 - 1) (1 + \ada z) \xgp_k, 
\end{align}
as required. 
\enproof 

\Lem \label{lem: nak comm phi x}
We have 
\[
z\mi[\xgm_1, \Phig(z\mi)]_{q^{-2}} -q^{-2} [\xgm_{0}, \Phig(z\mi)]_{q^{2}} = 0. 
\]
\enlem

\Proof
By \cite[(1.2.7)]{Nakajima-01}, 
\[
(z - q^{-2}w) \Phig(z\mi) \Xgm(w\mi) = (q^{-2}z - w) \Xgm(w\mi)\Phig(z\mi). 
\]
Separating homogeneous terms of degree zero in $w$ and setting $w=1$, we get  
\[
\Phig(z\mi) (z \xgm_0 - q^{-2} \xgm_1) = q^{-2} (z\xgm_0 - \xgm_1) \Phig(z\mi).
\]
Rearranging the equation above and dividing by $z$ yields the required formula. 
\enproof

\subsection{Basic implications of the Lu--Wang relations} 

We first use the relations \eqref{eq: rel1}--\eqref{eq: rel3} to express $\Hg_1$ in terms of the Drinfeld generators of $\qaa$. 

\Lem
We have 
\begin{align}
\label{eq: H1}
\Hg_1 =& \ 
(\CC\hg_1 - \hg_{-1}) - q^{-2} \cc_1 [\xgp_0,\xgp_{-1}]_{q^6}. 
\end{align}
\enlem

\Proof
Recall from \eqref{eq: H1 simple} that 
\eq
[\Ag_{-1},\Ag_0]_{q^{-2}} = q^{-4}\cc_0^{-1}\Hg_1. 
\eneq
Applying \eqref{eq: droq->sl} to the LHS, the relation becomes 
\begin{align}
\label{eq: H-lemma calc} [\Ag_{-1},\Ag_0]_{q^{-2}} =& \ [\xgm_1,\xgm_0]_{q^{-2}} +q^{-2} \cc_0\mi\cc_1 [\Kg\xgp_{-1},\Kg^{-1}\xgp_0]_{q^{-2}} \\
\quad& \ - q^{2}\cc_1[\xgm_1,\Kg^{-1}\xgp_0]_{q^{-2}} - q^{-4}\cc_0\mi[\Kg\xgp_{-1},\xgm_0]_{q^{-2}}\\
=& \ -(q^{-6}\cc_0\mi\cc_1 [\xgp_0, \xgp_{-1}]_{q^6} 
+ \cc_1\Kg\mi[\xgm_1,\xgp_0] + q^{-4}\cc_0\mi\Kg[\xgp_{-1},\xgm_0]) \\
=& \ -q^{-6}\cc_0\mi\cc_1 [\xgp_0,\xgp_{-1}]_{q^6} + \cc_1\hg_1 - q^{-4}\cc_0\mi \hg_{-1}. 
\end{align}
Note that the first term on the RHS of \eqref{eq: H-lemma calc} vanishes thanks to \eqref{eq: sl2 4x rel}. 
\enproof

Recall the notation for homogeneous components in the Drinfeld gradation introduced just after \eqref{eq: qOns Dr gr}. Our next task is to compute $\Ag_r^{(-1)}$, i.e., the component of degree~$-1$.

Relation \eqref{eq: rel2} (with $m=1$) implies that 
\eq
\label{eq: Ar}
\Ag_{r} = [\bHg_1,\Ag_{r-1}] + \CC \Ag_{r-2}.
\eneq

\Lem 
\label{lem: Ar-1}
We have: 
\begin{alignat}{3}
\label{eq: Armin} 
\Ag_r\gup{l} =& \ 0 &\quad \quad &(r \in \Z, \ l \leq -2 \mbox{ or } l \in 2\Z), \\
\label{eq: Armin1} 
\Ag_r\gup{-1} =& \ \xgm_{-r} &\quad \quad &(r \in \Z). 
\end{alignat}
\enlem

\Proof
Recall from \eqref{eq: droq->sl} that, assuming \eqref{s-assumption}, generators $\Ag_{-1}$ and $\Ag_0$ only have terms of degree $\pm 1$. 
The first equality \eqref{eq: Armin} then follows from \eqref{eq: Ar}, by downward and upward inductions, using the fact that $\Hg_1$ has only terms of degree $0$ and $2$ (see \eqref{eq: H1}). 

The second equality holds for $r = -1,0$ by \eqref{eq: droq->sl}. 
By induction and Lemma \ref{eq: H1}, we get 
\begin{align}
\Ag_r\gup{-1} =& \ [\bHg_1\gup{0},\Ag_{r-1}\gup{-1}] + \CC\Ag_{r-2}\gup{-1} \\ 
=& \  [\CC\bhg_1 - \bhg_{-1},  \xgm_{-r+1} ] + \CC  \xgm_{-r+2} \\
=& \  -\CC \xgm_{-r+2} + \xgm_{-r} + \CC \xgm_{-r+2} \\ 
=& \  \xgm_{-r} 
\end{align}
for $r \geq 1$. A similar argument using the relation $\Ag_{r} = \CC^{-1}(\Ag_{r+2} - [\bHg_1, \Ag_{r+1}])$ and downward induction establishes \eqref{eq: Armin1} for $r \leq -2$. 
\enproof

In terms of generating series, \eqref{eq: Armin1} can be rephrased as
\eq
\Asmone = \Xgm(z^{-1}). 
\eneq

\subsection{$\Apsone$ calculation} \label{sec: A+1 calc}

Abbreviate
\eq \label{eq: Omega factor}
\factor = (1+\ada z)(1 - \CC \adb z). 
\eneq

\Lem
\label{lem: A1}
We have 
\begin{align}
\Apsone =& \ 
 -\cc_1 \Phig(z^{-1})\Omega\mi\big((q^2 - q^{-2})\Xgp_{\leq -2}(q^{-2}z^{-1}) + q^2\xgp_0 + q^4\xgp_{-1}z\big). 
\end{align}
\enlem 

\Proof
First, observe that
\begin{align}
\left((1-\ad_{\bHg_1}z-\CC z^2)\Aps\right)^{(1)} =& \ (1-\ad_{\bHg_1^{(0)}}z-\CC z^2)\Apsone - z\ad_{\bHg_1^{(2)}}\Apsmone \\
=& \ \Omega\Apsone - z \adc \Apsmone, \\
(\Ag_0 + \CC \Ag_{-1}z)^{(1)} =& \ -\cc_1(\Kg \xgp_{-1}z + q^2 \Kg^{-1}\xgp_0). 
\end{align}
Hence, by \eqref{eq: Aplus series}, we get
\eq \label{eq: Omega A K formula}
\Omega\Apsone - z \adc \Apsmone = 
-\cc_1(\Kg \xgp_{-1}z + q^2 \Kg^{-1}\xgp_0). 
\eneq

We now calculate the second term on the LHS. We have 
\begin{align} \label{eq: ad H A sum}
\adc \Apsmone =& \  -q^{-2} \frac{\cc_1}{[2]} \sum_{r \geq 0} \big[[\xgp_0,\xgp_{-1}]_{q^6}, \xgm_{-r}\big]z^r. 
\end{align}
Let us compute the terms in the sum on the RHS. We get 
\begin{align}
\xgp_0\xgp_{-1} \xgm_{-r} =& \ \xgp_0 (\xgm_{-r}\xgp_{-1} - \qq \phig_{-r-1}) \\
=& \ - \qq((\phig_{-r} - \delta_{r=0} \psig_0)\xgp_{-1} + \xgp_0\phig_{-r-1}) + \xgm_{-r}\xgp_0\xgp_{-1}, \\
\xgp_{-1}\xgp_{0} \xgm_{-r} =& \ \xgp_{-1} (\xgm_{-r}\xgp_{0} - \qq (\phig_{-r} - \delta_{r=0}\psig_{0})) \\
=& \ - \qq(\phig_{-r-1}\xgp_{0} + \xgp_{-1}(\phig_{-r}- \delta_{r=0} \psig_0)) + \xgm_{-r}\xgp_{-1}\xgp_{0} \\ 
=& \ - q^{-2}\qq((\phig_{-r} - \delta_{r=0} \psig_0)\xgp_{-1} + \xgp_0\phig_{-r-1}) + \xgm_{-r}\xgp_{-1}\xgp_{0}, 
\end{align}
where the last equality follows from \eqref{eq: Nakajima hom}. 
Hence, using Lemma \ref{lem: xphi}, we get 
\begin{align}
\big[[\xgp_0,\xgp_{-1}]_{q^6}, \xgm_{-r}\big] =& \ q^2[2](\phig_{-r}\xgp_{-1} + \xgp_0\phig_{-r-1} - \delta_{r=0}\Kg\xgp_{-1}) \\
=& \ [2] \big( q^4\phig_{-r-1}\xgp_0 + q^6\phig_{-r}\xgp_{-1} + (q^4 - 1) \sum_{s = 2}^{r+1} q^{2s} \phig_{s-r-1} \xgp_{-s} - \delta_{r=0} q^2\Kg \xgp_{-1} \big). 
\end{align}
Combining the equality above with \eqref{eq: ad H A sum}, we obtain  
\begin{align}
\adc \Apsmone =& \  \cc_1\Kg \xgp_{-1} -q^{-2} \cc_1\sum_{r \geq 0} ( q^4\phig_{-r-1}\xgp_0 + q^6\phig_{-r}\xgp_{-1})z^r \\
\quad& \ - \cc_1(q^2 - q^{-2}) \sum_{r \geq 0}\sum_{s = 2}^{r+1} q^{2s} \phig_{s-r-1} \xgp_{-s} z^r. \\
=& \ \cc_1 \Kg \xgp_{-1} + \cc_1q^{2}\Kg\mi \xgp_0 z^{-1} - \cc_1q^2\Phig(z\mi)\xgp_{0}z\mi \\
\quad& \  - \cc_1q^4\Phig(z\mi)\xgp_{-1} 
- \cc_1(q^2-q^{-2})\Phig(z^{-1})\Xgp_{\leq -2}(q^{-2}z^{-1})z\mi. 
\end{align}

Finally, combining the equality above with \eqref{eq: Omega A K formula}, we get 
\begin{align}
\Omega \Apsone =& \ 
- \cc_1 \Phig(z^{-1})\big((q^2 - q^{-2})\Xgp_{\leq -2}(q^{-2}z^{-1}) + q^2\xgp_0 + q^4\xgp_{-1}z\big).
\end{align}
Dividing both sides by $\Omega$ and using the fact that $\Omega$ commutes with $\Phig(z^{-1})$ completes the proof. 
\enproof

In light of Lemma \ref{lem: A1}, we define
\eq \label{eq: mmathbb A}
\mathbb{A} =  -\cc_1 \Omega\mi\big((q^2 - q^{-2})\Xgp_{\leq -2}(q^{-2}z^{-1}) + q^2\xgp_0 + q^4\xgp_{-1}z\big), 
\eneq
so that $\Apsone = \Phig(z^{-1}) \mathbb{A}$. 

\subsection{Proof of the factorization theorem}

Throughout this subsection the symbol ``$\equiv$" denotes equalities holding modulo $U_+$. %define 
\Lem \label{lem: Theta degrees} \label{lem: Theta 2 parts}
We have 
\begin{align} 
(1 - q^{-2}\CC z^2)\Thgs(z) \equiv& \  
\frac{(\Kg - \CC \Kg^{-1} z^2)\Phig(z^{-1}) }{q-q\mi} \\
 \quad& \ -q^2 \cc_0 z \Phig(z^{-1}) [\mathbb{A}, (1+q^2\ada z)x_1^-] \quad 
\mod \fext{U_+}{z}. 
\end{align} 
\enlem

\Proof
Recall from \eqref{eq: Theta series in A terms} that 
\[
(1 - q^{-2}\CC z^2)\Thgs(z) = (1-\CC z^2)\Thg_0 + z^2\cc_1\mi\CC\left(z\mi[\Ag_{-1}, \Aps]_{q^{-2}}  -q^{-2} [\Ag_{0}, \Aps]_{q^{2}}\right). 
\]
The elements $A_{-1}$ and $A_0$ only contain terms of degree $\pm 1$, whereas $\Aps$ has additional terms of higher degree. Since we are working modulo $U_+$, we can discard all the homogeneous components $\Apsvar{l}$ with $l \geq 2$. 
Therefore, 
\begin{align}
z\mi[\Ag_{-1}, \Aps]_{q^{-2}}  -q^{-2} [\Ag_{0}, \Aps]_{q^{2}} \equiv& \ 
(z\mi \Ag_{-1} - q^{-2}\Ag_{0} )(\Apsmone + \Apsone) \\
& \ + (\Apsmone + \Apsone) (- q^{-2}z\mi \Ag_{-1} + \Ag_0). 
\end{align} 

Clearly, we can ignore terms of degree $2$ above. 
Next, we show that terms of degree $-2$ cancel out, i.e., 
\eq \label{eq: LHS Ax}
(z\mi \Ag_{-1}\gup{-1} - q^{-2}\Ag_{0}\gup{-1})\Apsmone  
 + \Apsmone(- q^{-2}z\mi \Ag_{-1}\gup{-1} + \Ag_0\gup{-1}) = 0. 
\eneq 
Indeed, computing the LHS of \eqref{eq: LHS Ax}, we get 
\begin{align}
\mbox{LHS} =
 (z\mi \xgm_{1} -q^{-2} \xgm_0)\Xgm_{\leq 0}(z\mi) 
+ \Xgm_{\leq 0}(z\mi) (q^{-2}z\mi \xgm_1 + \xgm_0). 
\end{align}
The $z^s$-term in the resulting series is equal to 
\[
(\xgm_{1}\xgm_{-s-1} -q^{-2} \xgm_0\xgm_{-s}) + (-q^{-2}\xgm_{-s-1}\xgm_1 + \xgm_{-s} \xgm_0) 
\] 
if $s \geq 0$ and $[\xgm_1,\xgm_0]_{q^{-2}}$ if $s=-1$. In both cases, the expression vanishes thanks to relation \eqref{eq: sl2 4x rel}, with $k=0$ and $l = -s-1$. 

It follows that 
\begin{align} \label{eq: new A-A equiv}
\ \ z\mi[\Ag_{-1}, \Aps]_{q^{-2}}  -q^{-2} [\Ag_{0}, \Aps]_{q^{2}} \equiv& \ 
(z\mi \Ag_{-1}\gup{1} - q^{-2}\Ag_{0}\gup{1})\Apsmone \\
& \ + (z\mi \Ag_{-1}\gup{-1} - q^{-2}\Ag_{0}\gup{-1})\Apsone \\
& \ + \Apsmone(- q^{-2}z\mi \Ag_{-1}\gup{1} + \Ag_0\gup{1}) \\
& \ + \Apsone(- q^{-2}z\mi \Ag_{-1}\gup{-1} + \Ag_0\gup{-1}).  
\end{align} 
Our next step is to compute the four summands on the RHS above. 

For $s \geq 0$, we have 
\begin{align}
(\Ag_{-1}\gup{1}\xgm_{-s-1} - q^{-2} \Ag_{0}\gup{1}\xgm_{-s}) =& \ 
(-q^{-4}\cc_0^{-1}\Kg\xgp_{-1}\xgm_{-s-1} + \cc_1\Kg\mi \xgp_0 \xgm_{-s}),  \\ 
(-q^{-2}\xgm_{-s-1}\Ag_{-1}\gup{1} + \xgm_{-s}\Ag_{0}\gup{1}) =& \ 
(q^{-4}\cc_0^{-1}\Kg\xgm_{-s-1}\xgp_{-1} - \cc_1\Kg\mi \xgm_{-s}\xgp_0 ),  \\ 
\end{align}
and summing the two terms above yields 
\[
\frac{  \Kg \phig_{-s-2} - \CC \Kg^{-1} (\phig_{-s} - \delta_{s=0} \psig_0) }{\cc_1\mi \CC(q-q\mi)}.
\]
We also have 
\[
\Ag_{-1}\gup{1}\xgm_{0} -q^{-2}\xgm_0\Ag_{-1}\gup{1} = 
-q^{-4}\cc_0\mi \Kg [\xgp_{-1}, \xgm_0] = (q-q\mi)\mi \cc_1 \CC\mi \Kg \phig_{-1}. 
\] 
Putting the equalities above into a series, we get 
\begin{align} \label{eq: Ax->phi 1}
%\big( \Ag_{-1}\gup{1} \Xgm_{\leq 0}(-q^{-2}z\mi)z\mi -& \ q^{-2}\Ag_0\gup{1}\Xgm_{\leq 0}(-q^{-2}z\mi)\big) 
(z\mi \Ag_{-1}\gup{1} - q^{-2}\Ag_{0}\gup{1})\Apsmone \ +& \ \Apsmone(- q^{-2}z\mi \Ag_{-1}\gup{1} + \Ag_0\gup{1}) 
 \equiv \\
\equiv& \ \frac{(\Kg - \CC \Kg^{-1} z^2)\pmb{\phi}^-(z^{-1}) + (\CC z^2 - 1)}{\cc_1\mi \CC z^2(q-q\mi)}.
\end{align} 

Moreover, by Lemma \ref{lem: A1}, we have 
\begin{align}
(z\mi \Ag_{-1}\gup{-1} - q^{-2}\Ag_{0}\gup{-1})\Apsone =& \ (z\mi \Ag_{-1}\gup{-1} - q^{-2}\Ag_{0}\gup{-1})\Phig(z^{-1}) \mathbb{A} \\ 
=& \ (x_1^- z\mi -q^{-2} x_0^-) \Phig(z^{-1}) \mathbb{A} \\ 
=& \ \Phig(z^{-1})(q^{-2}x_1^- z\mi - x_0^-)  \mathbb{A}, 
\end{align}
where $\mathbb{A}$ is as in \eqref{eq: mmathbb A}, and the third equality follows from Lemma \ref{lem: nak comm phi x}. 
On the other hand, 
\[
\Apsone(- q^{-2}z\mi \Ag_{-1}\gup{-1} + \Ag_0\gup{-1}) = 
\Phig(z^{-1}) \mathbb{A} (-q^{-2}\xgm_1z\mi + \xgm_0). 
\] 
Hence 
\begin{align} \label{eq: Ax->phi 2}
(z\mi \Ag_{-1}\gup{-1} - q^{-2}\Ag_{0}\gup{-1})\Apsone \ +& \
\Apsone(- q^{-2}z\mi \Ag_{-1}\gup{-1} + \Ag_0\gup{-1}) = \\
=& \ -q^{-2}z\mi\Phig(z^{-1}) [\mathbb{A}, (1+q^2\ada z)x_1^-]. 
\end{align} 
Substituting \eqref{eq: Ax->phi 1} and \eqref{eq: Ax->phi 2} back into \eqref{eq: new A-A equiv} and \eqref{eq: Theta series in A terms} completes the proof. 
\enproof 

\Lem \label{lem: ad A one formula}
We have
\begin{align}
(1-q^2\ada z)\mathbb{A} =
- \cc_1q^2 (1-\ada z)\Xgp_{\geq0}(\CC z). 
\end{align}
\enlem

\Proof
We begin with an auxiliary calculation. 
Abbreviate
\begin{align}
%M_1 =& \ \cc_1q^2(1-z\mi)\Omega\mi (\Kg \xgp_{-1}z - \Kg^{-1}\xgp_0) \\
M =& \ -[(q^2 - q^{-2})\Xgp_{\leq -2}(q^{-2}z^{-1}) + q^2\xgp_0 + q^4\xgp_{-1}z]. 
\end{align}
We have 
\begin{align}
- (1-q^2\ada z)(q^2 - q^{-2})\Xgp_{\leq -2}(q^{-2}z^{-1}) =& \ (q^2 - q^{6})\xgp_{-2} z^2 \\
-(1-q^2\ada z)(q^2\xgp_0 + q^4\xgp_{-1}z) =& \ q^{6}\xgp_{-2}z^2 -q^2\xgp_0. 
\end{align}
Hence
\eq
(1-q^2\ada z)M = q^2\xgp_{-2} z^2 - q^2 \xgp_0 = (1 - \ada^2z^2)(-q^2\xgp_0)
\eneq
and so 
\eq \label{eq: Omega M}
(1-q^2\ada z)\Omega\mi M = -q^2\frac{1-\ada z}{1 - \CC \adb z}\xgp_0. 
\eneq

Recall that, according to Lemma \ref{lem: A1},  
\eq \label{eq: Apsone Omega M}
\mathbb{A} = \cc_1 \Omega\mi M. 
\eneq
Therefore,
\begin{align}
(1-q^2\ada z)\mathbb{A} =& \ \cc_1 (1-q^2\ada z)\Omega\mi M \\
=& \ -q^2\cc_1 \frac{1-\ada z}{1 -\CC \adb z}\xgp_0 \\
=& \ - \cc_1q^2 (1-\ada z)\Xgp_{\geq0}(\CC z), 
\end{align}
yielding the desired formula. The first equality follows from \eqref{eq: Apsone Omega M}, the second from \eqref{eq: Omega M} and the third from \eqref{eq: hx rel}.
\enproof

We now have all the ingredients to prove the main result of this section. 

\Thm[Factorization theorem] \label{thm: fact}
Assuming that $\sss_0 = \sss_1 = 0$, 
we have 
\eq
\Thgs(z) \equiv \frac{1-\CC z^2}{(q-q\mi)(1-q^{-2}\CC z^2)}\pmb{\phi}^-(z\mi)\pmb{\phi}^+(\CC z)  \mod \fext{U_+}{z}.  
\eneq
\enthm

\Proof 
We have the following chain of equalities: 
\begin{align}
-q^2\cc_0 z \Phig(z^{-1}) [\mathbb{A}, (1 \ +& \ q^2\ada z)x_1^-] = 
-q^2 \cc_0 z \Phig(z^{-1}) [(1-q^2\ada z)\mathbb{A}, x_1^-] \\
=& \ \CC  z \Phig(z^{-1})[(1-\ada z)\Xgp_{\geq0}(\CC z), \xgm_1] \\
=& \ \CC  z \Phig(z^{-1})[\Xgp_{\geq0}(\CC z), (1+\ada z)\xgm_1] \\ 
=& \ \frac{\CC z\Phig(z^{-1})}{q-q\mi} \big( (\CC z)\mi\Psig_{\geq1}(\CC z)  - z \Psig(\CC z) + z\Kg\mi  \big) \\
=& \ \frac{\Phig(z^{-1})}{q-q\mi} \big((1-\CC z^2) \Psig(\CC z) - \Kg + \CC z^2 \Kg\mi \big). 
\end{align} 
The first and third equalities follow from the fact that $\ada$ is a derivation acting trivially on the generators $\psi_i, \phi_i$. 
The second equality follows from Lemma \ref{lem: ad A one formula}. 
The theorem now follows directly from Lemma \ref{lem: Theta 2 parts}. 
\enproof

\Rem \label{rem: Thg grave more natural}
The choice of the imaginary root vectors $\Thg_s$ in new Drinfeld presentation of the $q$-Onsager algebra in \cite[\S 2.7]{lu-wang-21} was suggested by the construction of the $\imath$Hall algebra of the projective line from \cite{lu-ruan-wang-23}. %using $1$-periodic complexes of coherent sheaves on the projective line. 
Earlier, a different set of root vectors $\acute{\Thg}_s$ (see \eqref{eq: Thg acute def}) was found in \cite{bas-kol-20}. The corresponding Drinfeld presentation can be found in \cite[\S 5.3]{lu-wang-21}. 
Theorem \ref{thm: fact} shows that the generators $\acute{\Thg}_s$, further renormalized by $(q-q\mi)$, are indeed more natural from the point of view of representation theory, as already observed in Remark \ref{rem: Theta vs acute Theta}. 
Define 
\eq \label{eq: grave vs acute theta}
\grave{\Thg}_s = (q-q\mi) \acute{\Thg}_s.
\eneq
Then Theorem \ref{thm: fact} takes the following more elegant form. 
\enrem

\Cor \label{cor: fact}
Assuming $\sss_0 = \sss_1 = 0$, we have 
\eq \label{eq: fact thm renormalized} 
\Thgsr(z) \equiv \pmb{\phi}^-(z\mi)\pmb{\phi}^+(\CC z)  \fext{\mod U_+}{z}.  
\eneq
\encor

%%%%%%%%%%%%%%%%%%%%%%%%%%%%%
%%%%%%%%%%%%%%%%%%%%%%%%%%%%%

\section{Drinfeld rational fractions} 

\label{sec: drf}

\subsection{Spectra of restricted representations} 
\label{subsec: spectra} 

Let $\cart$ (resp.\ $\cartoq$) be the commutative subalgebra of $\uLsl$ (resp.\ $\Oq$) generated by elements $\phi^{\pm}_i$ (resp.\ $\Theta_i$). 
Recall that every irreducible finite dimensional $\uLsl$-module $V$ can be endowed with a $\Z_{\leq 0}$-grading (compatible with the $\Z$-grading on $\uLsl$) such that the highest weight vector is in degree zero. Indeed, this is clear, in the case of evaluation representations, from their explicit description, and follows in the general case from the fact that every finite dimensional irreducible module is obtained as a tensor product of evaluation representations \cite[Theorem 4.11]{chari-pressley-qaa}. %references 
This grading also induces a (vector space) filtration with $F_i(V) = \{ v \in V \mid \deg v \geq i \}$. Since this filtration is preserved by the action of $\cart$, we can consider the associated graded $\gr (V)$ as a $\cart$-module. 

The $\Oq$-module $\eta_{\mathbf{c},0}^*(V)$ has a decomposition into generalized $\cartoq$-eigenspaces. 
Let $w$ be generalized $\cartoq$-eigenvector and $\overline{w}$ its image in $\gr(V)$. Corollary \ref{cor: fact} implies that $\overline{w}$ is a generalized $\cart$-eigenvector, and that the $\cartoq$-eigenvalues of $w$ can be computed using the $\cart$-eigenvalues of $\overline{w}$. More precisely, in terms of generating series, the eigenvalue of~$\Thgsr(z)$ on $w$ is the same as the $\pmb{\phi}^-(z\mi)\pmb{\phi}^+(\CC z)$-eigenvalue on $\overline{w}$. This fact has many interesting consequences. 

\subsection{Frenkel--Reshetikhin-type theorem} 

Let $w \in \eta_{\mathbf{c},0}^*(V)$ be a generalized $\cartoq$-eigenvector. Let $\pmb{d}_w^\pm(z) \in \fext{\kor}{z^{\pm 1}}$ and $\pmb{D}_w(z) \in \fext{\kor}{z}$ be defined by 
\eq
\pmb{\phi}^+(z) \cdot \overline{w} = \pmb{d}_w^+(z)  \overline{w}, \qquad 
\pmb{\phi}^-(z) \cdot \overline{w} = \pmb{d}_w^-(z) \overline{w}, \qquad
\pmb{\grave{\Thg}}(z) \cdot w = \pmb{D}_w(z) w.
\eneq 
By \cite[Proposition 1]{FrenRes}, there exist polynomials $Q(z), R(z) \in \C[z]$ with constant term $1$ such that 
\eq
\label{eq: FR formula}
\pmb{d}_w^\pm(z)= q^{\deg Q - \deg R} \frac{Q(q^{-1}z)R(qz)}{Q(qz)R(q^{-1}z)}. 
\eneq 

Given a polynomial $P(z) \in \C[z]$ with constant term $1$, let $P^\dag(z)$ be the polynomial with constant term $1$ whose roots are obtained from those of $P(z)$ via the transformation $a \mapsto \CC\mi a\mi$; and let $P^*(z)$ be the polynomial with constant term $1$ whose roots are the inverses of the roots of $P(z)$. Observe that 
\eq \label{eq: P* formula} 
P(z) = \gamma_P z^{\deg P} P^\dag(\CC\mi z\mi), \qquad 
P(z) = \gamma_{P} z^{\deg P} P^*(z^{-1}),
\eneq
where $\gamma_P$ is the product of the roots of $P^*(z)$. 

\Cor
\label{cor: FR thm Oq}
The eigenvalues of $\pmb{\grave{\Thg}}(z)$ on $w$ are of the form: 
\eq \label{eq: FR formula Oq}
\pmb{D}_w(z) = \frac{\ourQ(q^{-1}z)}{\ourQ(qz)} \frac{\ourQ^\dag(qz)}{\ourQ^\dag(q^{-1}z)}, \footnote{I thank Jian-Rong Li for calling my attention to, and helping me verify this formula.}
\eneq
as elements of $ \fext{\kor}{z}$, where $\ourQ(z)$ is a polynomial with constant term $1$. Explicitly, 
\[
\ourQ(z) = {Q(\CC z)R^*(z)}, \quad  \ourQ^\dag(z) = {R(\CC z)Q^*(z)}. 
\]
\encor 

\Proof
Corollary \ref{cor: fact}, together with the argument from \S \ref{subsec: spectra}, implies that 
\[ 
\pmb{D}_w(z) = \pmb{d}_w^+(\CC  z) \pmb{d}_w^-(z\mi).  
\]
Using Frenkel and Reshetikhin's formula \eqref{eq: FR formula}, as well as \eqref{eq: P* formula}, we get 
\begin{align} \label{eq: Dw long formula}
\pmb{D}_w (z) =& \    \frac{Q^*(q z)R^*(q\mi z)Q(\CC q^{-1}z)R(\CC qz)}{Q(\CC q z)R(\CC q^{-1} z)Q^*(q\mi z)R^*(qz)} \\
=& \   \left( \frac{Q(\CC q^{-1}z)R^*(q\mi z)}{Q(\CC q z)R^*(q z)} \right) \left( \frac{R(\CC q z)Q^*(q z)}{R(\CC q^{-1}z)Q^*(q\mi z)}\right), 
\end{align} 
as claimed. 
\enproof 

\Rem
It is interesting to compare Corollary \ref{cor: FR thm Oq} with the corresponding results for quantum affine algebras. 
 
Firstly, we observe that formula \eqref{eq: FR formula Oq} has an extra symmetry that is absent in \eqref{eq: FR formula}. Namely, the polynomials $Q(z)$ and $R(z)$ are not related in any special way, unlike $\ourQ(z)$ and $\ourQ^\dag(z)$. 

Secondly, for quantum affine algebras, the case of the highest weight vector $v_0$ is special because $R(z) = 1$ in \eqref{eq: FR formula}. This special case was proven by Chari and Pressley in \cite[Theorem 3.4]{chari-pressley-qaa}, and the polynomial $Q(z)$ is called the Drinfeld polynomial. In contrast, in our case, both $\ourQ(z)$ and $\ourQ^\dag(z)$ still appear in \eqref{eq: FR formula Oq} (with opposite $q^2$-shifts) for $v_0$. 

This fact admits the following interpretation. The notion of a highest weight, familiar from the representation theory of quantum affine algebras, is not directly applicable to modules restricted to the $q$-Onsager algebra. The main obstacle is the lack of canonical annihilation operators analogous to $x_k^+$. In fact, \eqref{eq: Armin1} implies that, in general, none of the Drinfeld generators $A_k$ of $\Oq$ annihilate the highest weight vector $v_0$. 
\enrem 

\Defi
We will refer to the polynomial $\ourQ(z)$ from Corollary \ref{cor: FR thm Oq} as the \emph{Drinfeld polynomial}  associated to an $\cartoq$-eigenspace in $\eta_{\mathbf{c},0}^*(V)$.  
\enDefi

\subsection{Rational fractions} 
\label{ss: DP and DRF}

We will now present another interpretation of the eigenvalues series $\pmb{D}_w(z)$, namely, as a ratio of rational functions with an additional symmetry. 

\Defi
The \emph{Drinfeld rational fraction}\footnote{Note that the Drinfeld rational fraction is defined up to sign. This ambiguity can be resolved by choosing a branch of the square root function or considering $\Z/2\Z$-orbits of such functions.} (DRF) associated to an $\cartoq$-eigenspace in $\eta^*_{\mathbf{c},0}(V)$ is 
\begin{align} \label{eq: F as ratio P P*}
F(z) =& \ \gamma_{\ourQ}\mi\CC^{\frac{1}{2}\deg \ourQ} \frac{\ourQ(z)}{\ourQ^\dag(z)} \in \kor(z), 
\end{align} 
\enDefi 

\Thm \label{cor: Drinfeld rational functions}
The eigenvalues of $\pmb{\grave{\Thg}}(z)$ on $w$ can be expressed as the following ratio of Drinfeld rational fractions: 
\eq \label{eq: D ratio rat func F}
\pmb{D}_w(z) = \frac{F(q^{-1}z)}{F(qz)}. 
\eneq
Moreover, the rational function $F(z)$:
\be
\item does not have a pole or zero at $z=0$ or $z = \infty$, 
\item satisfies the $\CC$-\emph{twisted unitarity property}: 
\eq \label{eq: def twisted unitarity property}
F(q\CC \mi z\mi) = F(q^{-1}z)\mi. 
\eneq
\ee
\enthm

\Proof 
The identity \eqref{eq: D ratio rat func F} follows directly from the definition \eqref{eq: F as ratio P P*} and \eqref{eq: Dw long formula}. 
Since, by \cite[Proposition 1]{FrenRes}, neither $P_V(z)$ nor $P^*_V(z)$ has a zero at $z=0$, part~(1) follows. 
Next, using \eqref{eq: P* formula} , we compute
\begin{align} 
F(q\CC \mi z\mi) =& \ \gamma_{\ourQ}\mi \CC^{\frac{1}{2}\deg \ourQ}  \frac{\ourQ(q\CC\mi z\mi)}{\ourQ^\dag(q\CC\mi z\mi)} \\ 
=& \ \gamma_{\ourQ}\mi \CC^{\frac{1}{2}\deg \ourQ} \frac{\gamma_{\ourQ}(q\CC\mi z\mi)^{\deg \ourQ}\ourQ^\dag(q\mi z)}{\ourQ^\dag(q\CC\mi z\mi)} \\ 
=& \ (q\CC^{-\frac{1}{2}}z\mi)^{\deg \ourQ}  \frac{\ourQ^\dag(q\mi z)}{\ourQ^\dag(q\CC\mi z\mi)}, \\ 
F(q^{-1}z)\mi =& \ \gamma_{\ourQ} \CC^{-\frac{1}{2}\deg \ourQ} \frac{\ourQ^\dag(q\mi z)}{\ourQ(q\mi z)} \\ 
=&  \ \gamma_{\ourQ} \CC^{-\frac{1}{2}\deg \ourQ} \frac{\ourQ^\dag(q\mi z)}{\gamma_{\ourQ}(q\mi z)^{\deg \ourQ} \ourQ^\dag(q\CC\mi z\mi)} \\
=& \ (q\CC^{-\frac{1}{2}}z\mi)^{\deg \ourQ}  \frac{\ourQ^\dag(q\mi z)}{\ourQ^\dag(q\CC\mi z\mi)}, 
\end{align} 
which yields (2). 
\enproof

\Rem 
Equation \eqref{eq: D ratio rat func F} has the form of a ratio involving a single $q^{-2}$-shift, familiar from the theory of Drinfeld polynomials. However, we emphasize that this kind of formula only holds in our case when polynomials are replaced by rational functions. 

The twisted unitarity property can be construed as an analogue of the condition $u^{\deg P}P(u\mi) = q^{-\deg P}P(uq^2)$ in \cite[Theorem 4.6]{gow-molev}, which appeared in the setting of the symplectic twisted $q$-Yangian. 

Note that, by \cite[Theorem 11.7]{kolb-14}, the $q$-Onsager algebra corresponds to the \emph{orthogonal} twisted $q$-Yangian of rank one. Therefore, the results of \cite{gow-molev} are not applicable in our case. There is also a major difference in methodology as \emph{op.\ cit.} uses the FRT presentation. 
\enrem

\Rem
The notion of a Drinfeld polynomial, in the context of the $q$-Onsager algebra, has already been introduced by Ito--Terwilliger in \cite{ito-ter-08, ito-ter-10}. 
However, the definition in \emph{op.\ cit.} is completely different from our definitions of Drinfeld polynomials and rational   fractions. In particular, \emph{a priori}, it has nothing to do with Lu--Wang's new Drinfeld presentation, introduced only twelve years later. 

On the other hand, the two concepts share many similarities. They both appear to play an important role in the classification of finite dimensional modules and behave nicely with respect to tensor products. It would be interesting to understand the precise connection between the two notions. 
\enrem 

\Rem
The idea of replacing the RHS of $\pmb{d}_V^+(z) = \frac{P_V(q^{-2}z)}{P_V(z)}$ with more general rational functions has already appeared in \cite{her-jim-12} in the study of category $\mathcal{O}$ for $U_q\widehat{\mathfrak{b}}$, and in \cite{her-shifted} in the study of shifted quantum affine algebras. 
\enrem

%%%%%%%%%%%%%%%%%%%%%%%%%%%%%%%%%
%%%%%%%%%%%%%%%%%%%%%%%%%%%%%%%%%

\section{One-dimensional representations}
\label{sec: odr}

The action of $\Oqc$ on a one-dimensional representation $V$ is determined by the two scalars by which the generators $\Ag_{-1}$ and $\Ag_0$ act on $V$. Therefore, given any pair of such scalars $t_{-1}, t_0$, we can construct $V$ by restricting the trivial $\uLsl$-representation via \eqref{eq: Kolb emb}, with $s_0 = q^2 \cc_0 t_{-1}$ and $s_1 = t_0$. 
Let $\kor_{\cs} = \eta_{\cs}^* (\kor)$ denote this restriction, and let 
$\xi_{\cs} \colon \Oq \to \kor$ be the character corresponding to the action of $\Oq$ on $\kor_{\cs}$, as in \eqref{eq: xi s char}. 
Abbreviate
\[
\DD_{\cs}(z) = \DD_{\kor_{\cs}}(z) = \xi_{\cs}(\Thgsr(z)). 
\]

\Lem \label{lem: theta char formula}
We have
\eq \label{eq: 1 dim D formula}
\DD_{\cs}(z) = %\cc_1\mi\frac{(1-\CC z^2)^2 + \alpha \CC z^2  +  \beta \CC z(1+\CC z^2)}{(1-\CC z^2)^2},
\frac{q\mi(q-q\mi)^2\cc_1\mi\CC(\alpha z + \beta(1+\CC z^2))z + (1-\CC z^2)^2}{(1-\CC z^2)^2}, 
\eneq
where
\[ 
\alpha = \CC (q^{-2}\cc_0\mi \sss_0)^2 + \sss_1^2, \quad
\beta = q^{-2}\cc_0\mi \sss_0 \sss_1.  \] 
\enlem 

\Proof
Applying $\xics$ to \eqref{eq:Thr} and multiplying by $q-q\mi$, we get 
\eq \label{eq: one dim Theta gen formula} 
\DD_{\cs}(z) = 1 + (1-\CC z^2)\mi \cc_1\mi \CC z q\mi(q-q^{-1})^2 (\xics(\Ag_{-1}) + \xics(\Ag_0)z) \xics(\Aps). 
\eneq
Relation \eqref{eq: rel2} implies that 
\[
\xi_{\cs}(\Ag_{2r}) = \CC^r\sss_1, \quad \xi_{\cs}(\Ag_{2r-1}) = \CC^r q^{-2}\cc_0\mi \sss_0,
\] 
so 
\eq \label{eq: Aps xics formula}
\xics(\Aps) = (1-\CC z^2)\mi(\sss_1 +\CC q^{-2}\cc_0\mi \sss_0z). 
\eneq
Substituting \eqref{eq: Aps xics formula} into \eqref{eq: one dim Theta gen formula} yields the required formula. 
\enproof 

\Pro
There exists a rational function $F_{\cs}(z) \in \kor(z)$ such that
\[
\DD_{\cs}(z) = \left(F_{\cs}(z)F_{\cs}(\CC\mi z\mi)\right)\mi. 
\]
This function is unique up to sign and transforming its poles under the map $a \mapsto \CC\mi a\mi$. 
Moreover: 
\be
\item $F_{\cs}(z) = \pm 1$ if and only if $\sss_0 = \sss_1 = 0$, 
\item the numerator and denominator of $F_{\cs}(z)$ have degree one if and only if $\sss_0 \neq 0 \neq \sss_1$ and $\frac{\sss_0}{\sss_1} = \pm \sqrt{\frac{\cc_0}{\cc_1}}$.  
\ee 
\enpro

\begin{proof}
Let 
\[
G(z) = q\mi(q-q\mi)^2\cc_1\mi\CC(\alpha z + \beta(1+\CC z^2))z + (1-\CC z^2)^2 = \CC^2 \prod_{i=1}^4 (z - \gamma_i).
\]
Then, by Lemma \ref{lem: theta char formula}, 
\[
\DD_{\cs}(z) = \frac{-\CC\mi G(z)z^{-2}}{(1-\CC z^2)(1-\CC\mi z^{-2})}. 
\]
By the $C$-symmetry of $\DD_{c,s}(z)$, we have 
\[ 
G(z)z^{-2} = \CC^2G(\CC\mi z\mi)z^2.
\]
Hence $\gamma_3 = \CC\mi \gamma_1\mi$, $\gamma_4 = \CC\mi \gamma_2\mi$ and 
\[
G(z) = (\gamma_1\gamma_2)\mi \CC^2(z-\gamma_1)(z-\gamma_2)(\CC\mi z\mi - \gamma_1)(\CC\mi z\mi - \gamma_1)z^2. 
\]
Therefore
\eq \label{eq: one dim F definition}
F_{\cs}(z) = \pm i \frac{(\gamma_k\gamma_l\CC)^{1/2}(z^2 - \CC\mi)}{(z-\gamma_k)(z-\gamma_l)}, %choice of signs
\eneq
where $k \in \{1, 3\}$ and $l \in \{2, 4\}$. This proves the existence and uniqueness of $F_{\cs}(z)$, completing the proof of the main statement. 

We now deduce the three additional statements. We will use the following equalities, obtained by 
considering the coefficients of the polynomial $G(z)$ on $z^3$ and $z^2$, respectively: 
\begin{align}
\label{eq: alpha beta 1}
\tilde{\beta} = q\mi(q-q\mi)^2\cc_1\mi \beta =& \ \gamma_1 + \gamma_2 + \CC\mi(\gamma_1\mi + \gamma_2\mi), \\ 
\label{eq: alpha beta 2}
\tilde{\alpha} = q\mi(q-q\mi)^2\cc_1\mi \CC \alpha =& \ (\gamma_1\CC + \gamma_1\mi) (\gamma_2\CC + \gamma_2\mi) + 4\CC . 
\end{align}

It follows directly from the definition \eqref{eq: one dim F definition} that 
$F_{\cs}(z) = \pm 1$ if and only if $\gamma_1 = \CC^{-1/2}$ and $\gamma_2 = -\CC^{-1/2}$. By \eqref{eq: alpha beta 1}--\eqref{eq: alpha beta 2}, this is the case if and only if $\alpha = \beta = 0$, which is equivalent to $s_0 = s_1 = 0$. This proves (1). 

Similarly, it follows from \eqref{eq: one dim F definition} that $F_{\cs}(z)$ is a ratio of polynomials of degree one if and only if $\gamma_1 = \pm \CC^{-1/2}$ and $\gamma_2 \neq \mp \CC^{-1/2}$. Without loss of generality, take $\gamma_1 = \CC^{-1/2}$. Then we obtain the quadratic equations:
\[
\gamma_2^2 + (2\CC^{-1/2}-\tilde{\beta}) + \CC\mi = 0, \quad \gamma_2^2 + (2\CC^{-1/2}-\textstyle\frac{1}{2}\CC^{-3/2}\tilde{\alpha}) + \CC\mi = 0, 
\]
which imply that 
$\alpha = 2\CC^{1/2}\beta$, i.e., 
\[
\CC (q^{-2}\cc_0\mi \sss_0)^2 + \sss_1^2 = 2\CC^{1/2}q^{-2}\cc_0\mi \sss_0 \sss_1. 
\]
After rearranging, we get   
\[
(\cc_1^{1/2}\sss_0 - \cc_0^{1/2}\sss_1)^2 = 0.  
\]
This proves (2). 
\end{proof} 

\Rem
Note that, unlike the Drinfeld rational functions from \S \ref{ss: DP and DRF}, the functions $F_{\cs}(z)$ are, in general, \emph{not} twisted unitary (in the sense of \eqref{eq: def twisted unitarity property}). 
\enrem

%%%%%%%%%%%%%%%%%%%%%%%%%%%%%%%%%
%%%%%%%%%%%%%%%%%%%%%%%%%%%%%%%%%

\section{Multiplicative property} 
\label{sec: mult property}

The goal of this section is to prove a more general version of the factorization theorem, 
where we drop the assumption that $s_0 = s_1 = 0$. To do this, we prove a multiplicative property for the series $\Thgsr(z)$ with respect to actions on one-dimensional representations. The key ingredients are formulas for the coproduct of $\Aps$ and $\Thgsr(z)$, modulo certain ideals, established in \S \ref{subsec: cop Aps}--\ref{subsec: cop Thgs} below. 
Once again, we emphasize that the parameters $s_0, s_1$ can be \emph{arbitrary} throughout \S \ref{sec: mult property}.

\subsection{Coproduct formula for $\Aps$} 
\label{subsec: cop Aps}

It follows directly from \eqref{eq: coprod on U} and \eqref{eq: droq->sl} that 
\[
\tops(\Ag_{-1}) = 1 \otimes \etacz(\Ag_{-1}) + \etacs(\Ag_{-1}) \otimes \Kg, \quad 
\tops(\Ag_{0}) = 1 \otimes \etacz(\Ag_{0}) + \etacs(\Ag_{0}) \otimes \Kg\mi. 
\] 
For the sake of readability, we abbreviate $Y \otimes 1 = \etacs(Y) \otimes 1$ and $1 \otimes Y = 1 \otimes \etacz(Y)$, for $Y \in \Oqc$. 
Relation \eqref{eq: H1 simple} implies that 
\[
\Delta(\bHg_1) = 1 \otimes \bHg_1 + \bHg_1 \otimes 1 - (q-q\mi)\Ag_0 \otimes \xgp_{-1} - \CC (q-q\mi) \Ag_{-1} \otimes \xgp_0.  
\]
We will denote the four summands on the RHS above as $Q_1, \cdots, Q_4$, respectively. 
We also abbreviate
\begin{align}
S^r_{1} =& \ 1 \otimes \Ag_r, \\
S^r_{2} =& \  \Ag_{0} \otimes \phig_{-r}, \\
S^r_{3} =& \ \sum_{i=0}^{r-1} \Ag_{r-i} \otimes  \phig_{-i}. %\\
%S^r_{2+3} =& \ \sum_{i=0}^{r} \Ag_{r-i} \otimes (-q^2)^i \phig_{-i}
\end{align}
We write sums of the terms above as, \emph{e.g.}, $Q_{1+2}, S^r_{2+3}, S^r_{1+2+3}$, etc. 
Throughout this subsection the symbol ``$\equiv$" denotes equalities holding modulo $1 \otimes U_{\geq2}$. 

\Lem \label{lem: S ind 1}
We have
\[
S^{r+1}_{1+2+3} \equiv  [Q_{1+2}, S^r_{1+2+3}] + \CC S^{r-1}_{1+2+3} + [Q_{3+4}, S^r_1]    \mod 1 \otimes U_{\geq2}. 
\]
\enlem 

\Proof
We use relation \eqref{eq: Ar}, i.e., 
\eq \label{eq: Ar vs H C spec}
\Ag_{r+1} = [\bHg_1,\Ag_{r}] + \CC \Ag_{r-1},
\eneq
throughout. 

Now assume $r \geq 2$. 
Using \eqref{eq: Ar vs H C spec} and Lemma \ref{lem: Ar-1}, we get  
\begin{align}
[ Q_1, S^r_1] + \CC S^{r-1}_1 =& \ S^{r+1}_1, \\ 
[ Q_3, S^r_1 ] \equiv& \ \Ag_0 \otimes (q-q\mi) [\xgp_{-1},  \xgm_{-r}] \\
	=& \ S^{r+1}_2. 
\end{align} 

Next, using \eqref{eq: Ar vs H C spec} repeatedly, we see that 
\[
[Q_2, S^r_{3}] + \CC S^{r-1}_{2+3} = \sum_{i=0}^{r-1} \Ag_{r+1-i} \otimes \phig_{-i}.  
\]
Moreover,
\begin{align}
[ Q_4, S^r_1] + [Q_2, S^r_{2}] = & \  - \CC \Ag_{-1} \otimes (q-q\mi)[\xgp_0,  \xgm_{-r}] 
+ [\bHg_1, \Ag_0] \otimes \phig_{-r} \\
=& \ ([\bHg_1, \Ag_0] + \CC \Ag_{-1}) \otimes  \phig_{-r} = \Ag_1 \otimes  \phig_{-r}. 
\end{align} 
Hence
\[
[Q_2, S^r_{2+3}] + \CC S^{r-1}_{2+3} + [ Q_4, S^r_1] = S^{r+1}_3. 
\]

Next, observe that 
\[
[Q_1, S^r_{2+3}] \equiv [Q_2, S^r_1] \equiv  0. 
\]
Summing up our calculations so far, we, therefore, obtain
\[
[Q_{1+2}, S^r_{1+2+3}] + \CC S^{r-1}_{1+2+3} + [Q_{3+4}, S^r_1] \equiv S^{r+1}_{1+2+3}, 
\]
as claimed. 
\enproof

Write 
\[
\cops = (\wt{\xi}_{\mathbf{c}, \mathbf{s}} \otimes \id) \circ  \Delta_{\mathbf{c},0},
\]
where $\wt{\xi}_{\mathbf{c}, \mathbf{s}}$ is as in \eqref{eq: xi s char}. 
From now on all the elements on the LHS of the tensor symbol are understood as images under $\xi_{\mathbf{c}, \mathbf{s}}$. For example, $\Ag_{0} \otimes \phig_{-r} = \xi_{\mathbf{c}, \mathbf{s}}(\Ag_{0}) \otimes  \phig_{-r}$, etc. In particular, any two elements on the LHS of the tensor symbol commute. 
Abbreviate 
\[
\kappa = -(q-q\mi)\big(\xi_{\mathbf{c}, \mathbf{s}}(\Ag_0)+\CC  \xi_{\mathbf{c}, \mathbf{s}}(\Ag_{-1}) \adb\big)
\]
so that 
\[
Q_{3+4} = 1 \otimes \kappa \xgp_{-1}. 
\]

Also abbreviate
\begin{align}
S^r_4 =  \sum_{\substack{i+j+k=r; \\ i\geq 1; j,k\geq0}}  \Ag_{k}  \otimes \kappa\mathbbm{h}_{i-1}[\xgp_{-1},  \phi_{-j}^-], 
\end{align}
where $\mathbbm{h}_i$ denotes the $i$-th degree complete homogeneous symmetric polynomial in $-\ada$ and $\CC \adb$,  
so that 
$
\sum_{i=0}^\infty \mathbbm{h}_i z^{i} = \factor\mi, 
$
with $\factor$ as in \eqref{eq: Omega factor}. 

\Lem  \label{lem: S ind 2}
We have 
\[
S^{r+1}_4 \equiv 
[Q_{1+2+3+4}, S^r_4] + \CC S^{r-1}_4 + [Q_{3+4}, S^r_{2+3}]  \mod 1 \otimes U_{\geq2}.  
\]
\enlem

\Proof
Clearly 
\[
[ Q_{2+3+4},  S_4^r] \equiv 0.  
\]
Next, observe that
\begin{align}
[\bHg_1, \mathbbm{h}_{i-1}] + \CC \mathbbm{h}_{i-2} \equiv& \ (\CC \adb - \ada) \mathbbm{h}_{i-1} + \CC \mathbbm{h}_{i-2} - \CC \adb\ada \mathbbm{h}_{i-2} = \mathbbm{h}_{i}. 
\end{align}
Hence  
\begin{align}
[ Q_1, S^r_4 ] + \CC S^{r-1}_4 + [Q_{3+4}, S^r_{2+3}] \equiv S^{r+1}_4,
\end{align}
as claimed. 
\enproof

\Pro \label{pro: twisted primitive}
For each $r \geq 0$, 
\[
\cops(\Ag_r) \equiv S^r_{1+2+3+4}  \mod 1 \otimes U_{\geq2}. 
\]
Hence 
\begin{align} \label{eq: cop A three terms}
\cops(\Aps) \equiv& \ 1 \otimes \Aps + \Aps \otimes \Phig(z\mi) + 
\Aps \otimes \kappa \Gamma,
\end{align}
where
\[
\Gamma = (1-q^2\ada z)\mi (1-\CC \adb z)\mi (q^2-1) \Phig(z\mi) \xgp_{-1}z. 
\] 
\enpro

\Proof
The first statement follows by induction, using \eqref{eq: Ar vs H C spec} and Lemmas \ref{lem: S ind 1}--\ref{lem: S ind 2}. It also follows directly from the definition of $S^r_k$ that 
\[
\sum_{r \geq 0} S^r_1 = 1 \otimes \Aps, \qquad \sum_{r \geq 0} S^r_{2+3} = \Aps \otimes \Phig(z\mi). 
\] 
Moreover, 
\begin{align}
\sum_{r \geq 1} S^r_4 z^r =& \ \Aps \otimes z\kappa \factor\mi [\xgp_{-1}, \Phig(z\mi)] 
= \Aps \otimes \kappa \Gamma, 
\end{align}
where the last equality follows from Lemma \ref{lem: xgp phi comm}. 
\enproof 

We will also use the notation
\[
\widetilde{\Gamma} = (1-q^2\ada z)\mi (1-\CC \adb z)\mi (q^2-1) \xgp_{-1}z 
\]
so that $\Gamma = \Phig(z\mi) \widetilde{\Gamma}$.

\subsection{Coproduct formula for $\Thgsr(z)$}
\label{subsec: cop Thgs}

Throughout this subsection the symbol ``$\equiv$" denotes equalities holding modulo $1 \otimes U_+$. 

\Pro \label{pro: s coproduct Theta} 
We have 
\[
\cops(\Thgsr(z)) \equiv \xics(\Thgsr(z)) \cdot \eta_{\mathbf{c},0}(\Thgsr(z)) \mod \fext{1 \otimes U_+}{z}. 
\]
\enpro

\Proof
We use \eqref{eq: Theta series in A terms} and \eqref{eq: cop A three terms} throughout. 
Let us denote the three summands on the RHS of \eqref{eq: cop A three terms} as $R_1, R_2, R_3$, respectively. 
Then 
\begin{align}
(1-\CC z^2)\cops(\Thg_0) \ + z^2\cc_1\mi\CC\big(z\mi[1\otimes\Ag_{-1}, R_1]_{q^{-2}}  -q^{-2} [&1\otimes\Ag_{0}, R_1]_{q^{2}}\big) = \\
 =& \ 1 \otimes (1 - q^{-2} \CC z^2) \Thgs(z), 
\end{align}
and 
\begin{align}
 z^2\cc_1\mi\CC\big(&z\mi[\Ag_{-1}\otimes \Kg, R_1]_{q^{-2}} -  q^{-2} [\Ag_{0}\otimes \Kg\mi , R_1]_{q^{2}}\big) \equiv  \\
\equiv& \ z^2\cc_1\mi\CC\left(z\mi \Ag_{-1}\otimes [\Kg,  \Xgm(z^{-1})]_{q^{-2}}  - q^{-2} \Ag_{0}\otimes [\Kg\mi ,  \Xgm(z^{-1})]_{q^{2}}\right) \\
\equiv& \ 0,  
\end{align}
because the terms on the RHS of the tensor symbols vanish. This implies that 
\begin{align}
(1-\CC z^2)\cops(\Thg_0) \ +  z^2\cc_1\mi\CC\big(z\mi[\cops(\Ag_{-1}), R_1]_{q^{-2}} -& q^{-2} [\cops(\Ag_{0}), R_1]_{q^{2}}\big) \equiv \\
\equiv& \ 1 \otimes (1 - q^{-2} \CC z^2) \Thgs(z). \label{eq: R1 Teil}
\end{align}

Next, it follows from Lemma \ref{lem: nak comm phi x} that 
\begin{align}
z^2\cc_1\mi\CC\big(z\mi[&1\otimes\Ag_{-1}, R_2]_{q^{-2}} - q^{-2} [1\otimes\Ag_{0}, R_2]_{q^{2}}\big) \equiv \\
\equiv& \ z^2\cc_1\mi\CC\Aps\otimes \left( z\mi[\xgm_1, \Phig(z\mi)]_{q^{-2}} -q^{-2} [\xgm_{0}, \Phig(z\mi)]_{q^{2}}\right) \\
\equiv& \ 0, 
\end{align} 
and it follows directly from the definitions that 
\begin{align}
 z^2\cc_1\mi\CC\big(z\mi[&\Ag_{-1}\otimes \Kg, R_2]_{q^{-2}} - q^{-2} [\Ag_{0}\otimes \Kg\mi , R_2]_{q^{2}}\big) \equiv  \\
\equiv& \ z^2\cc_1\mi\CC(1-q^{-2}) \Aps \otimes \left(z\mi \xics(\Ag_{-1}) \Kg + \xics(\Ag_{0}) \Kg\mi \right)\Phig(z^{-1}). 
\end{align} 
This implies that 
\begin{align}
z^2\cc_1\mi\CC\big(z\mi[&\cops(\Ag_{-1}), R_2]_{q^{-2}} - q^{-2} [\cops(\Ag_{0}), R_2]_{q^{2}}\big) \equiv \\
\equiv& \ z^2\cc_1\mi\CC(1-q^{-2}) \Aps \otimes \left(z\mi \xics(\Ag_{-1}) \Kg + \xics(\Ag_{0}) \Kg\mi \right)\Phig(z^{-1}). \label{eq: R2 Teil}
\end{align}

Moreover, 
\begin{align}
z^2\cc_1\mi\CC\big( z\mi[\cops(\Ag_{-1}),& \ R_3]_{q^{-2}} - q^{-2} [\cops(\Ag_{0}), R_3]_{q^{2}}\big) \equiv \\ 
\equiv& \ z^2\cc_1\mi\CC \Aps \otimes \big(z\mi[\xgm_1, \kappa \Gamma]_{q^{-2}}  -q^{-2} [\xgm_0, \kappa \Gamma]_{q^{2}}\big) \\ 
\equiv& \ -q^{-2}z\cc_1\mi\CC \Aps \otimes \Phig(z\mi) [\kappa \widetilde{\Gamma}, (1 + q^2 \ada z) \xgm_1] \\ 
\equiv& \ -q^{-2}z\cc_1\mi\CC \Aps \otimes \Phig(z\mi) [(1 - q^2 \ada z) \kappa \widetilde{\Gamma}, \xgm_1] \\ 
\equiv& \ q\mi z^2(q-q\mi)^2 \cc_1\mi\CC \Aps \otimes
\Phig(z\mi)\cdot \\
& \cdot \big[(1 - \CC \adb z)\mi  (\xics(\Ag_0)\xgp_{-1} + \CC \xics(\Ag_{-1})\xgp_0), \xgm_1\big],
\end{align}
where Lemma \ref{lem: nak comm phi x} was used in the third equality. 
Let us calculate the term in the last line above. Observe that 
\begin{align}
(1 - \CC \adb z)\mi  (\xics(\Ag_0)\xgp_{-1} \ +& \ \CC \xics(\Ag_{-1})\xgp_0) = \\
=& \ \CC (\xics(\Ag_{-1}) +  \xics(\Ag_{0})z) \Xgp_{\geq0}(\CC z) + \xics(\Ag_0)\xgp_{-1}.  
\end{align}
Hence 
\begin{align}
(q-q\mi)\big[&(1 - \CC \adb z)\mi (\xics(\Ag_0)\xgp_{-1} + \CC \xics(\Ag_{-1})\xgp_0), \xgm_1 \big]  = \\
=& \ 
z\mi (\xics(\Ag_{-1}) +  \xics(\Ag_{0})z) \Psig_{\geq 1}(\CC z)  + \xics(\Ag_0) (\Kg - \Kg\mi) \\
=& \  (\xics(\Ag_{-1})z\mi + \xics(\Ag_{0})) \Psig(\CC z)  - \xics(\Ag_{-1})z\mi \Kg - \xics(\Ag_0)\Kg\mi. 
\end{align} 
It follows that 
\begin{align} \label{eq: R3 Teil} 
z^2\cc_1\mi\CC\big( z\mi[\cops(\Ag_{-1}), R_3]_{q^{-2}} \ -& \ q^{-2} [\cops(\Ag_{0}), R_3]_{q^{2}}\big) \equiv \\ 
\equiv& \ q\mi z^2(q-q\mi) \cc_1\mi\CC \Aps \otimes
\Phig(z\mi)\cdot \\ 
& \{(\xics(\Ag_{-1})z\mi + \xics(\Ag_{0})) \Psig(\CC z) \\
& - \xics(\Ag_{-1})z\mi \Kg - \xics(\Ag_0)\Kg\mi\}. 
\end{align}
The summand corresponding to the last line above cancels out \eqref{eq: R2 Teil}. Therefore, collecting 
\eqref{eq: R1 Teil}--\eqref{eq: R3 Teil} yields 
\begin{align}
(1 - q^{-2} \CC z^2)\cops(\Thgs(z)) \equiv& \ 1 \otimes (1 - q^{-2} \CC z^2) \Thgs(z) \\
& + z (1-q^{-2})\cc_1\mi \CC (\Ag_{-1} + \Ag_{0}z) \Aps \otimes \pmb{\phi}^-(z\mi)\pmb{\phi}^+(\CC z). 
\end{align}
Using Corollary \ref{cor: fact}, this simplifies to 
\begin{align}
(1 - q^{-2} \CC z^2)\cops(\Thgs(z)) \equiv& \ 1 \otimes (1 - q^{-2} \CC z^2) \Thgs(z) \\
& + z (1-q^{-2})\cc_1\mi \CC (\Ag_{-1} + \Ag_{0}z) \Aps \otimes \Thgsr(z). 
\end{align} 
The proposition now follows directly from \eqref{eq: one dim Theta gen formula}. 
\enproof 

\subsection{General factorization theorem}

We can now prove the main result of this section. 

\Thm \label{cor: group like Theta}
For any $\mathbf{s} \in \C^2$, the following factorization and coproduct formulae hold: 
\be
\item generalized factorization: 
\[
\eta_{\cs}(\Thgsr(z)) \equiv \xics(\Thgsr(z)) \cdot \eta_{\mathbf{c},0}(\Thgsr(z)) \quad  \mod \fext{(1 \otimes U_+)}{z},  
\]
\item coproduct formula - the series $\Thgsr(z)$ is ``group-like", i.e., 
\[
\Delta_{\cs}(\Thgsr(z)) \equiv \eta_{\cs}(\Thgsr(z)) \otimes \eta_{\mathbf{c},0}(\Thgsr(z)) \quad  \mod \fext{(\Oqp{\mathbf{c},\mathbf{s}} \otimes U_{+})}{z}.  
\]
\ee
\enthm 

\Proof 
The first part is a direct consequence of Proposition \ref{pro: s coproduct Theta} and \eqref{eq: eta cs delta}. 
For the second part, we have 
\begin{align}
\tops(\Thgsr(z)) \equiv& \  \Delta \circ \etacs (\Thgsr(z)) \\ %(\wt{\xi}_{\mathbf{c},\mathbf{s}}\ten\id) \circ \Delta \circ  \eta_{\mathbf{c},0} \\ 
\equiv& \ \xics(\Thgsr(z)) \cdot \Delta \circ \eta_{\mathbf{c},0}(\Thgsr(z)) \\
\equiv& \ \xics(\Thgsr(z)) \cdot \Delta_{\mathbf{c},0}(\Thgsr(z)) \\
\equiv& \ \xics(\Thgsr(z)) \cdot \etacz(\Thgsr(z)) \otimes \etacz(\Thgsr(z)) \\
\equiv& \ \eta_{\cs}(\Thgsr(z)) \otimes \eta_{\mathbf{c},0}(\Thgsr(z)), 
\end{align}
where %the first equality follows from \eqref{eq: eta cs delta}, 
the second and fifth equalities follow from part (1) of the present theorem. 
%the second from Proposition \ref{pro: s coproduct Theta}, 

Let us prove the fourth equality, i.e., the step 
\begin{equation} \label{eq: prove}
\Delta_{\mathbf{c},0}(\Thgsr(z)) \equiv  \etacz(\Thgsr(z)) \otimes \etacz(\Thgsr(z)) \qquad \mbox{mod} \ \fext{(\mathcal{O}_q^{\mathbf{c},0} \otimes U_+)}{z}. 
\end{equation} 
By Corollary \ref{cor: fact}, we can write
\[
\Thgsr = \pmb{\phi}^+(C z)\pmb{\phi}^-(z^{-1}) + R, \qquad R \in \fext{U_+}{z}. 
\]
Therefore, 
\begin{align*}
\Delta(\Thgsr(z)) &= \Delta(\pmb{\phi}^+(Cz)) \Delta(\pmb{\phi}^-(z^{-1})) + \Delta(R) \\
&= (\pmb{\phi}^+(Cz) \otimes \pmb{\phi}^+(Cz) + R_1)(\pmb{\phi}^-(z^{-1}) \otimes \pmb{\phi}^-(z^{-1}) + R_2) + \Delta(R), 
\end{align*}
with $R_1, R_2 \in \fext{(U_- \otimes U_+)}{z}$. The first equality follows from the fact that the coproduct is an algebra homomorphism, and the second equality from the group-like (mod $U_- \otimes U_+$) property of $\pmb{\phi}^+(z), \pmb{\phi}^-(z)$. One should remark here that the latter is stated in a weaker form in some places in the literature. For example, \cite[Proposition 4.4.(iii)]{chari-pressley-qaa} gives a weaker statement modulo $UX_+ \otimes UX_- + UX_- \otimes UX_+$, where $UX_{\pm}$ are ideals generated by $x_{k}^{\pm}$. In fact, the proof in \emph{loc. cit.} already gives the stronger statement modulo $UX_- \otimes UX_+$ (given our conventions - they use a flipped coproduct). The version of the statement modulo $U_- \otimes U_+$ can be found in \cite[Proposition 7.1]{damiani-r}. Note that the notational conventions therein are somewhat different from ours, the positive (resp. negative) part in the Drinfeld grading is denoted by $\infty$ (resp. $-\infty$) for the elements $E_{(m\delta,i)}$, and vice versa for $F_{(m\delta,i)}$. 

It follows that 
\begin{align}
\Delta(\Thgsr(z)) &= \pmb{\phi}^+(Cz) \pmb{\phi}^-(z^{-1}) \otimes \pmb{\phi}^+(Cz) \pmb{\phi}^-(z^{-1}) + \Delta(R) + Z  \\
&= (\Thgsr(z)-R) \otimes \pmb{\phi}^+(Cz) \pmb{\phi}^-(z^{-1}) + \Delta(R) + Z, \label{eq: appendix}
\end{align}
where
\[
Z = (\pmb{\phi}^+(Cz) \otimes \pmb{\phi}^+(Cz)) R_2 + R_1 (\pmb{\phi}^-(z^{-1}) \otimes \pmb{\phi}^-(z^{-1})) \in \fext{(U_- \otimes U_+)}{z}.
\]
Since $\mathcal{O}_q^{\mathbf{c},0}$ is a coideal subalgebra, we must have $\Delta(\Thgsr(z)) \in \fext{(\mathcal{O}_q^{\mathbf{c},0} \otimes U)}{z}$. In particular, it follows from \eqref{eq: appendix} that 
\begin{equation} \label{eq: appendix2}
- R \otimes \pmb{\phi}^+(Cz) \pmb{\phi}^-(z^{-1}) + \Delta(R) + Z \in  \fext{(\mathcal{O}_q^{\mathbf{c},0} \otimes U)} {z}. 
\end{equation} 
Since $R \in \fext{U_+}{z}$, we can write 
\[
\Delta(R) = a + b, \qquad a \in \fext{(U \otimes U_+)}{z}, \ \ b \in \fext{(U_+ \otimes U_{\leq 0})}{z}. 
\] 
Decomposing \eqref{eq: appendix2} as a sum of homogeneous terms with respect to the right tensor factor, we deduce that 
\[
Z + a \in \fext{(\mathcal{O}_q^{\mathbf{c},0} \otimes U_+)}{z}, \qquad b - (R \otimes \pmb{\phi}^+(Cz) \pmb{\phi}^-(z^{-1})) \in \fext{((\mathcal{O}_q^{\mathbf{c},0}\cap U_+) \otimes U_{\leq 0})}{z}. 
\]
We claim that $\mathcal{O}_q^{\mathbf{c},0}\cap U_+ = \{0\}$. Assuming the claim, \eqref{eq: appendix} implies that 
\[
\Delta(\Thgsr(z)) = \Thgsr(z) \otimes \pmb{\phi}^+(Cz) \pmb{\phi}^-(z^{-1}) + Z + a = \Thgsr(z) \otimes \Thgsr(z) + (Z + a - \Thgsr(z) \otimes R),
\]
with $Z + a - (\Thgsr(z) \otimes R) \in \fext{\mathcal{O}_q^{\mathbf{c},0} \otimes U_+}{z}$, which implies \eqref{eq: prove}.

Let us finally prove the claim.\footnote{I thank Stefan Kolb for suggesting the following argument to me.} 
Let $G^+$ be the subalgebra of $U$ generated by $G_i=E_iK_i^{-1}$ (for $i=0,1$), and $G^-$ the subalgebra generated by $F_0, F_1$. Choose a subset $\mathbf{I}$ of multi-indices  $I=(i_1, i_2,\cdots,i_k)$ such that the monomials $F_I=F_{i_1} F_{i_2} \cdots F_{i_k}$ form a basis of $G^-$. Assume that 
$y = \sum_{I \in \mathbf{I}} a_I B_I$, where the $a_I$ are nonzero scalar coefficients, and $B_I = B_{i_1} B_{i_2} \cdots B_{i_k}$ (we know the $B_I$ form a basis of $\mathcal{O}_q^{\mathbf{c},0}$), 
is a non-zero element in $U_+$. Then $y$ has some leading terms $a_I F_I$ in $G^-$ where $|I|=k$ is maximal with $a_I$ nonzero. These terms have some Drinfeld degree $n$ (which can be zero if $I$ contains $0$ and $1$ the same number of times). However, $y$ also has  'opposite' leading terms $a_I G_I$, where $G_I=G_{i_1} G_{i_2} \cdots G_{i_k}$. Now, these terms are linearly independent from each other and from all the other homogeneous components of $y$. Moreover, if $a_I F_I$ has degree $n$, then $a_I G_I$ has degree $-n$. Hence, $y$ contains either a nonzero term of degree $0$ if $n$ happens to be zero, or there two terms, one of degree $n$ and one of degree $-n$ if $n \neq0$. In each case, the element $y$ does not lie in  $U_+$.  
\enproof

As a consequence, we also obtain a generalization of Theorem \ref{cor: Drinfeld rational functions}. 
Let $V$ be an irreducible finite dimensional $\uLsl$-module and $\eta_{\cs}^*(V)$ its restriction to $\Oqc$ via the embedding \eqref{eq: Kolb emb}, with $\sss_0, \sss_1$ arbitrary. 
Let $w \in \eta_{\mathbf{c},\mathbf{s}}^*(V)$ be a generalized $\cartoq$-eigenvector. Let $\DD_{\eta_{\cs}^*(V),w}(z) \in \fext{\kor}{z}$ be defined by 
$\pmb{\grave{\Thg}}(z) \cdot w = \DD_{\eta_{\cs}^*(V),w}(z) w$. 

The following result says that the series of the eigenvalues of the operators $\grave{\Thg}_r$ on $w$ factorizes into the the corresponding series for the \emph{standard} restriction $\eta_{\mathbf{c},0}^*(V)$ and the one-dimensional representation $\C_{\cs}$. 

\Cor
The series $\DD_{\eta_{\cs}^*(V)}(z)$ can be expressed as the product:
\[
\DD_{\eta_{\cs}^*(V),w}(z) = \DD_{\cs}(z) \cdot \DD_{\eta_{\mathbf{c},0}^*(V),w}(z).  
\]
\encor

\Proof
The statement follows directly from Theorem~\ref{cor: group like Theta}. 
\enproof 

%%%%%%%%%%%%%%%%%%%%%%%%%%%%
%%%%%%%%%%%%%%%%%%%%%%%%%%%%

\section{Affine QSP of split type}
\label{sec: qsp spl}

Our next goal is to generalize the factorization and coproduct formulae (Theorem~\ref{cor: group like Theta}) to QSP coideal subalgebras of higher rank. This will be achieved in \S \ref{sec: comp br gr}. Below we recall the definitions of quantum affine algebras and generalized $q$-Onsager algebras, as well as certain facts about braid group actions. 

\subsection{Quantum affine algebras} 
\label{subsec: qaa gen} 

Let $\indx_0 = \{ 1,\cdots, N \}$ and $\indx = \indx_0 \cup \{0\}$. Let $\g$ be a simple Lie algebra of type $\mathsf{ADE}$ with Cartan matrix $(a_{ij})_{i,j \in \indx_0}$, and $\widehat{\g}$ the corresponding untwisted affine Lie algebra with affine Cartan matrix $(a_{ij})_{i,j \in \indx}$. 
We use standard notations and conventions regarding root systems, weight lattices, Weyl groups, etc., as in, e.g., \cite[\S 3.1]{lu-wang-21}. In particular, we 
let $\alpha_i$ ($i \in \indx$) denote the simple roots of $\widehat{\g}$; let $P$ and $Q$ denote the weight and root  lattices of $\g$, respectively; and let $\omega_i \in P$ ($i \in \indx_0$) be the fundamental weights of $\g$. 

The \emph{quantum affine algebra} $\gqaa$ is the algebra with generators $\Eg_i^\pm, \Kg_i^{\pm 1}$ $(i \in \indx)$ and relations: 
\begin{align}
\Kg_i\Kg_i^{-1} &= \Kg_i^{-1}\Kg_i = 1, \\ 
\Kg_i\Kg_j &= \Kg_j \Kg_i, \\
\Kg_i\Eg_j^{\pm} &= q^{\pm a_{ij}} \Eg_j^\pm\Kg_i, \\
%\Kg_i\Eg_j^{\pm} &= q^{\mp2} \Eg_j^\pm\Kg_i  && \quad (i \neq j), \\
[\Eg_i^{+}, \Eg_j^{-}] &= \delta_{ij} \frac{\Kg_i - \Kg_i^{-1}}{q - q^{-1}}, \\
\Serre_{ij}(\Eg_i^{\pm}, \Eg_j^\pm) &= 0 \qquad (i \neq j),
\end{align}
where
\[
\Serre_{ij}(x,y) = \sum_{r=0}^{1-a_{ij}} (-1)^r \qbinom{1-a_{ij}}{r} x^{1-a_{ij}-r}yx^r. 
\]
We also write
\[ E_i = e_i^+, \quad F_i = e_i^-.\] 
We refer to the presentation above as the Kac--Moody presentation of $\gqaa$. 
The \emph{quantum loop algebra} $\qla$ is the quotient of $\gqaa$ by the ideal generated by the central element $\Kg_\delta - 1$, where $\Kg_\delta$ is defined as in, e.g., \cite[\S 1.3]{beck-94}. 

Let us recall the ``new" Drinfeld presentation of the quantum loop algebra $\qla$. By \cite{drinfeld-dp, beck-94}, $\qla$ is isomorphic to the algebra generated by $\xgpm_{i,k}, \hg_{i,l}, \Kg^{\pm 1}_{i}$, where $k \in \Z$, $l \in \Z - \{0\}$ and $i \in \indx_0$, subject to the following relations:
\begin{align}
\Kg_i\Kg_i^{-1} =& \ \Kg_i^{-1}\Kg_i = 1, \\
\Kg_i \Kg_j =& \ \Kg_j \Kg_i, \\ 
[\hg_{i,k}, \hg_{j,l}] =& \ 0, \\
\Kg_i\hg_{j,k} =& \ \hg_{j,k}\Kg_i, \\ 
\Kg_i \xgpm_{j,k} =& \ q^{\pm a_{ij}} \xgpm_{j,k}\Kg_i, \label{eq: g Kx rel} \\
[\hg_{i,k},\xgpm_{j,l}] =& \pm \textstyle \frac{[k\cdot a_{ij}]}{k}\xgpm_{j,k+l} \label{eq: g hx rel},\\
\xgpm_{i,k+1}\xgpm_{j,l} - q^{\pm a_{ij}} \xgpm_{j,l}\xgpm_{i,r+1} =& \ q^{\pm a_{ij}}\xgpm_{i,k}\xgpm_{j,l+1} - \xgpm_{j,l+1}\xgpm_{i,k}, \label{eq: g sl2 4x rel} \\
[\xgp_{i,k}, \xgm_{j,l}] =& \ \delta_{ij} \textstyle\frac{1}{q-q^{-1}}(\psig_{i, k+l} - \phig_{i, k+l}), \label{eq: g x+-}
\end{align}
\vspace{-0.25cm}
\[
\Sym_{k_1,\cdots,k_r} \sum_{s=0}^r (-1)^t \qbinom{r}{s} \xgpm_{i,k_1} \cdots \xgpm_{i,k_s}\xgpm_{j,l}\xgpm_{i,k_s+1} \cdots \xgpm_{i,k_r} = 0  \quad \mbox{for } r=1-a_{ij} \quad (i \neq j), 
\]
where $\Sym_{k_1,\cdots,k_r} $ denotes symmetrization with respect to the indices $k_1, \cdots, k_r$ and 
\[
\pmb{\phi}^{\pm}_i(z) = \sum_{k=0}^\infty \phipm_{i,\pm k} z^{\pm k} = \Kg_i^{\pm 1} \exp\left( \pm (q-q^{-1}) \sum_{k=1}^\infty \hg_{i,\pm k} z^{\pm k} \right). 
\] 
%index was missing

We will also consider the following version of the quantum loop algebra
\[ \qlae = \qla \otimes_{\C} \C[\KK_j^{\pm 1} \mid j \in \indx] \]
with extended scalars, matching the central elements in $\gOq$ (see \eqref{eq: central K}).  
We consider $\qlae$ as a $Q$-graded algebra with 
\[
\degdr \xgpm_{i,k} = \pm \alpha_i, \quad \degdr h_{i,k} = \degdr \KK_j = 0 \qquad (i \in \indx_0, \ j \in \indx). 
\] 
We will write $\degdr_i$ for degree along the simple root $\alpha_i$. Let $Q_i$ be the sublattice of~$Q$ spanned by all the simple roots except $\alpha_i$. 
Let $U_{d_i=r,+}$ be the subalgebra of $\qlae$ spanned by elements which are both: (i) of degree $r$ along $\alpha_i$, and (ii) positive $Q_i$-degree, i.e., non-negative along all the simple roots and positive along at least one simple root in $Q_i$. That is, 
\begin{align}
U_{d_i=r,+} =& \ \langle y \in \qlae \mid \degdr_i y = r,\ \forall j \neq i \in \indx_0\colon \degdr_j y \geq 0,\ \exists j \neq i\colon \degdr_j y > 0 \rangle, \\
U_{d_i \geq r,+} =& \ \langle y \in \qlae \mid \degdr_i y \geq r,\ \forall j \neq i \in \indx_0\colon \degdr_j y \geq 0,\ \exists j \neq i\colon \degdr_j y > 0 \rangle. 
\end{align}
We will also need the following subalgebras of $\qlae$: 
\begin{align}
U_+ =& \ \langle y \in \qlae \mid \forall j \in \indx_0\colon \degdr_j y \geq 0,\ \exists j \colon \degdr_j y > 0 \rangle, \\ 
%U_- =& \ \langle y \in \qlae \mid \forall j \in \indx_0\colon \degdr_j y \leq 0,\ \exists j \colon \degdr_j y < 0 \rangle, \\ 
U_{\geq 0} =& \ \langle y \in U \mid \forall j \in \indx_0\colon \degdr_j y \geq 0 \rangle, \\
U_{\neq i,+} =& \ \langle y \in \qlae \mid \forall j \neq i \in \indx_0\colon \degdr_j y \geq 0,\ \exists j \neq i\colon \degdr_j y > 0 \rangle, \\ 
%U_{\neq i,-} =& \ \langle y \in \qlae \mid \forall j \neq i \in \indx_0\colon \degdr_j y \leq 0,\ \exists j \neq i\colon \degdr_j y < 0 \rangle, \\ 
U_{i+,0} =& \ \langle y \in \qlae \mid \degdr_i y >0,\ \forall j \neq i \in \indx_0\colon \degdr_j y = 0 \rangle. 
%U_{i-,0} =& \ \langle y \in \qlae \mid \degdr_i y <0,\ \forall j \neq i \in \indx_0\colon \degdr_j y = 0 \rangle. 
\end{align} 
Finally, let $U^i_{+}$ be the subalgebra of $\qlae$ generated by $U_{\neq i,+}$ and $U_{i+,0}$. It is the span of all elements which are: either (a) of positive degree along $\alpha_i$ and non-negative degree along all the other simple roots, or (b) of arbitrary degree along $\alpha_i$ and positive $Q_i$-degree. %Similarly, let $U^i_{-}$ be the subalgebra of $\qlae$ generated by $U_{\neq i,-}$ and $U_{i-,0}$. 

Sometimes we will also consider $\qlae$ as a $\widehat{Q}$-graded algebra with 
\[
\deg E_i =1, \quad \deg F_i = -1, \quad \deg K_i = \deg \KK_i = 0 \qquad (i \in \indx),
\]
where $\widehat{Q}$ is the root lattice of $\widehat{\g}$.

\subsection{The braid group action} 
\label{subsec: br action Lusztig}

The Weyl group $W_0$ of $\g$ is generated by the simple reflections $s_i$ ($i \in \indx_0$), and acts on $P$ in the usual way. The extended affine Weyl group $\widetilde{W}$ is the semi-direct product $W_0 \ltimes P$. 
%where we denote $t_\omega = (1,\omega) \in \widetilde{W}$ for $\omega \in P$. %? 
It contains the affine Weyl group $W = W_0 \ltimes Q = \langle s_i \mid i \in \indx \rangle$ as a subgroup, and can also be realized as the semi-direct product 
\eq \label{eq: ext aff weyl 2 rel}
\widetilde{W} \cong \Lambda \ltimes W,
\eneq
where $\Lambda = P/Q$ is a finite group of automorphisms of the Dynkin diagram of $\widehat{\g}$. Let $\widetilde{\mathcal{B}}$ be the braid group associated to $\widetilde{W}$. 

%There is a braid group action given by %on...
It is well known \cite{lusztig-94} that,
for each $i \in \indx$, there exists an automorphism $T_i$ of $\qla$ such that $T_i(K_\mu) = K_{s_i\mu}$ and 
\begin{alignat}{3}
T_i(E_i) =& \ -F_iK_i, \quad& T_i(E_j) =& \ \sum_{r+s = -a_{ij}} (-1)^r q^{-r} E_i^{(s)}E_jE_i^{(r)}, \\ 
T_i(F_i) =& \ -K_i\mi E_i, \quad& T_i(F_j) =& \ \sum_{r+s = -a_{ij}} (-1)^r q^{r} F_i^{(r)}F_jF_i^{(s)}, 
\end{alignat} 
for $\mu \in P$ and $i \neq j$. 
In particular, if $a_{ij} = -1$ then
\eq \label{eq: br act Fi Ei ord}
T_i(E_j) = E_iE_j - q\mi E_j E_i, \qquad 
T_i(F_j) = F_jF_i - qF_iF_j. 
%-q(F_iF_j - q\mi F_jF_i). %
\eneq
For each $\lambda \in \Lambda$, there is also an automorphism $T_\lambda$ such that $T_\lambda(E_i) = E_{\lambda(i)}$, $T_\lambda(F_i) = F_{\lambda(i)}$ and $T_\lambda(K_i) = K_{\lambda(i)}$.   
These automorphisms satisfy the relations of the braid group $\widetilde{\mathcal{B}}$. 
We extend this action to an action on $\qlae$ 
%\[ \qlae = \qla \otimes_{\C} \C[\KK_j^{\pm 1} \mid j \in \indx] \]
by setting 
$T_i(\KK_\mu) = \KK_{s_i\mu}$ and $T_\lambda(\KK_i) = \KK_{\lambda(i)}$. 
%K_i mathbb, ordinary K too
%details...lusztig, LW 

\subsection{Generalized $q$-Onsager algebras}

%The $q$-Onsager algebra is a quantum symmetric pair coideal subalgebra of $\uqsl$. It has been studied extensively from the point of view of mathematical physics, especially in the context of quantum integrable systems with boundary \cite{}.
Baseilhac and Belliard introduced in \cite{BasBel} \emph{generalized} $q$-Onsager algebras, which are coideal subalgebras of $\gqaa$. These algebras are also known in the literature as affine $\imath$quantum groups of split type \cite{lu-wang-21}. We first recall the original ``Serre-type'' presentation of these algebras. 

%The defining relations of the $q$-Onsager algebra, the so-called $q$-Dolan-Grady relations, were found in \cite{Bas} to admit a higher rank generalization. 

\Defi
The universal generalized $q$-Onsager algebra $\gOq$ is the algebra generated by $\Bg_i$ and invertible central elements $\ck_i$ $(i \in \indx)$ subject to relations 
\begin{alignat}{3}
[\Bg_i, \Bg_j] =& \ 0 \qquad && \mbox{if } a_{ij} = 0, \\
\Bg_i^2 \Bg_j - [2]\Bg_i \Bg_j \Bg_i + \Bg_j \Bg_i^2 =& \ -q \ck_i \Bg_j \qquad && \mbox{if } a_{ij} = -1, \\
\sum_{r=0}^3 (-1)^r \sqbin{3}{r} \Bg_i^{3-r}\Bg_j\Bg_i^r =& \ - q \ck_i [2]^2[\Bg_i, \Bg_j] \qquad && \mbox{if } a_{ij} = -2. 
\end{alignat} 
\enDefi

It follows from \cite[Theorem 8.3]{kolb-14} that $Z(\gOq) = \C[ \ck_i^{\pm 1} \mid i \in \indx ]$.  

Recently, Lu and Wang \cite{lu-wang-21} found a ``Drinfeld-type" presentation of $\gOq$. We will use the following shorthand notations: 
\begin{align}
S(r_1, r_2 \mid s; i,j) =& \ \Ag_{i,r_1} \Ag_{i,r_2} \Ag_{j,s} - [2] \Ag_{i,r_1} \Ag_{j,s} \Ag_{i,r_2} + \Ag_{j,s}\Ag_{i,r_1} \Ag_{i,r_2}, \\
\mathbb{S}(r_1, r_2 \mid s; i,j) =& \ S(r_1, r_2 \mid s; i,j) + S(r_2, r_1 \mid s; i,j), \\ 
R(r_1, r_2 \mid s; i,j) =& \ q^2\ck_i\Ck^{r_1}\big( - \sum_{p \geq 0} q^{2p} [2] [\Thg_{i, r_2-r_1-2p-1}, \Ag_{j,s-1}]_{q^{-2}}  \Ck^{p+1} \\
& \ - \sum_{p \geq 1} q^{2p-1} [2] [\Ag_{j,s}, \Thg_{i, r_2-r_1-2p}]_{q^{-2}} \Ck^{p} - [\Ag_{j,s}, \Thg_{i,r_2-r_1}]_{q^{-2}} \big), \\
\mathbb{R}(r_1, r_2 \mid s; i,j) =& \ R(r_1, r_2 \mid s; i,j) + R(r_2, r_1 \mid s; i,j). 
\end{align}

\Defi 
Let $\gDrOq$ be the $\kor$-algebra generated by $\Hg_{i,m}$ and $\Ag_{i,r}$, where $m\geq1$, $r\in\Z$, and invertible central elements $\ccc_i$ ($i \in \indx_0$), $\CCC$, subject to the following relations:
\begin{align}
%\K_1\K_1^{-1} =1, \quad C C^{-1} &=1, \quad \K_1, C \text{ are central}, \label{iDRo0} \\
[\Hg_{i,m},\Hg_{j,n}] &=0, \label{eq: grel1} \\
[\Hg_{i,m}, \Ag_{j,r}] &= \textstyle \frac{[m \cdot a_{ij}]}{m} (\Ag_{j,r+m}- \Ag_{j,r-m}\CCC^m), \label{eq: grel2} \\ 
[\Ag_{i,r}, \Ag_{j,s}] &= 0 \quad \mbox{if} \quad a_{ij} = 0, \label{eq: grel3} \\ 
[\Ag_{i,r}, \Ag_{j,s+1}]_{q^{-a_{ij}}}  -q^{-a_{ij}} [\Ag_{i,r+1}, \Ag_{j,s}]_{q^{a_{ij}}} &= 0  \quad \mbox{if} \quad i\neq j, \label{eq: grel4} \\ 
\label{eq: grel5}
[\Ag_{i,r}, \Ag_{i,s+1}]_{q^{-2}}  -q^{-2} [\Ag_{i,r+1}, \Ag_{i,s}]_{q^{2}}
&= \ccc_i\CCC^r\Thg_{i,s-r+1} - q^{-2}\ccc_i\CCC^{r+1}\Thg_{i,s-r-1}  \\ \notag
&\quad  + \ccc_i\CCC^s\Thg_{i,r-s+1}  -q^{-2}\ccc_i \CCC^{s+1}\Thg_{i,r-s-1}, \\
\mathbb{S}(r_1, r_2 \mid s; i,j) &= \mathbb{R}(r_1, r_2 \mid s; i,j) \quad \mbox{if} \quad a_{ij} = -1. \label{eq: grel6}
\end{align}
where $m,n\geq1$; $r,s, r_1, r_2\in \Z$; 
and 
\[
1+ \sum_{m=1}^\infty (q-q^{-1})\Thg_{i,m} z^m  =  \exp \left( (q-q^{-1}) \sum_{m=1}^\infty \Hg_{i,m} z^m \right).
\]
By convention, $\Thg_{i,0} = (q-q^{-1})^{-1}$ and $\Thg_{i,m} = 0$ for $m\leq -1$. 
\enDefi 

%Given $Y \in \{H_{i,r}, A_{i,r}, \Theta_{i,r}\}$, let $Y\gup{n}$ denote the sum of homogeneous components $y$ of $Y$ with $\degdr_i y = n$. 

\Rem
%The relations \eqref{eq: Ons rel} are the same as the defining relations in \cite[Example 7.6]{kolb-14}. 
To match our conventions with those of \cite{lu-wang-21}, set 
\eq \label{eq: central K}
v = q, \qquad \mathbb{K}_i = q^2 \ck_i.
\eneq
\enrem

By \cite[Theorem 2.16]{lu-wang-21}, there is an algebra isomorphism 
\eq
%\label{eq: droq->oq}
\gOq \isoto{} \gDrOq.
\eneq
This isomorphism and its inverse can be described precisely in terms of a braid group action (see \cite[\S 3.3, 3.4]{lu-wang-21} and \S \ref{subsec: qsp br act} below). We will not need the details in what follows. We note, however, that the isomorphism sends
\eq
\Bg_{i} \mapsto \Ag_{i,0}, \quad \ck_i \mapsto \ck_i, \quad  q^{2N+2}\ck_\delta\mapsto\Ck , 
\eneq
for $i \in \indx_0$. 

\subsection{Coideal structures}
%decide whether to use ci or Ki 

The following is a straightforward generalization of \S \ref{subsec: coideal str sl2}-- \S \ref{subsec: C-duality} to higher ranks.  
Given $\mathbf{c} = (c_0, \cdots, c_N) \in (\C^\times)^{N+1}$, we have a central character 
\[
\chi_{\mathbf{c}} \colon Z(\gOq) \to \C^\times, \qquad \ck_i \mapsto c_i. 
\]
Let $\gOqc$ be the corresponding central reduction, i.e., the quotient of $\gOq$ by the ideal generated by $\ck_i - c_i$. 

By \cite[Theorem 7.1]{kolb-14}, for any $\mathbf{s} = (s_0, \cdots, s_N) \in \C^{N+1}$, 
there exists an injective algebra homomorphism 
\eq
\label{eq: Kolb emb gen}
\eta_{\mathbf{c},\mathbf{s}} \colon \gOqc \monoto \qla, \quad \Bg_i \mapsto \Eg_i^- - \cc_i \Eg_i^+ \Kg_i^{-1} + \sss_i \Kg_i^{-1} \quad (i \in \indx),
\eneq
giving $\gOqc$ the structure of a right coideal subalgebra of $\qla$. Abbreviate $\mathcal{O}_q^{\mathbf{c}, \mathbf{s}}(\widehat{\g}) = \eta_{\mathbf{c},\mathbf{s}}(\gOqc)$. 
There is also a generic version of this embedding
\[
\eta_{\mathbf{s}} \colon \gOq \monoto \qlae, \quad \Bg_i \mapsto \Eg_i^- - q^{-2}\KK_i \Eg_i^+ \Kg_i^{-1} + \sss_i \Kg_i^{-1} \quad (i \in \indx). 
\]
If $\mathbf{s} = (0, \cdots, 0)$, we abbreviate $\eta = \eta_{\mathbf{s}}$. 

%%%the other chapters, too 

We remark that the other results from \S \ref{sec: q-ons}, concerning one-dimensional modules and $C$-duality, as well as the rationality properties from \S \ref{sec: rat} carry over to the $\gOq$ case in a natural way. 

\subsection{The QSP braid group action} \label{subsec: qsp br act}
%cite Kolb, too

%There is a braid group action given by %on...
By \cite[Lemma 3.4]{lu-wang-21}, for each $i \in \indx$, there exists an automorphism $\Tbr_i$ of $\gOq$ such that $\Tbr_i(\KK_\mu) = \KK_{s_i\mu}$ and 
\[
\Tbr_i(B_j) = 
\left\{ \begin{array}{l l}
\KK_j\mi B_j & \mbox{if } i=j, \\
B_j & \mbox{if } a_{ij} = 0, \\ 
B_jB_i - qB_iB_j  & \mbox{if } a_{ij} = -1, \\ 
{[2]^{-1}}(B_jB_i^2-q[2]B_iB_jB_i+q^2B_i^2B_j) + B_j\KK_i & \mbox{if } a_{ij} = -2, 
%[2]^{-1}(B_jB_i^2-q[2]B_iB_jB_i+q^2B_i^2B_j) + B_j\KK_i & \mbox{if } a_{ij} = -2, 
\end{array} \right. 
\] 
for $\mu \in \Z\indx$ and $j \in \indx$. 
For each $\lambda \in \Lambda$, there is also an automorphism $\Tbr_\lambda$ such that $\Tbr_\lambda(B_i) = B_{\lambda(i)}$ and $\Tbr_\lambda(\KK_i) = \KK_{\lambda(i)}$.   
These automorphisms define an action of the braid group $\widetilde{\mathcal{B}}$. 
We will refer to this action as the Lu--Wang, QSP or coideal subalgebra braid group action to distinguish it from Lusztig's braid group action from \S \ref{subsec: br action Lusztig}. 

We will need the following lemma. 

\Lem \label{lem: br act on wt}
We have $\Tbr_w(B_i) = B_{wi}$ for $i \in \indx$ and $w \in \widetilde{W}$ such that $wi \in \indx$. 
\enlem

\Proof
See \cite[Lemma 3.5]{lu-wang-21}. 
\enproof

We remark that the analogous statement for the usual braid group action holds as well. 

\subsection{Rank one subalgebras} \label{subsec: rk one subalg} 

\nc{\Ui}{U_i^{\mathbb{K}}} 
\nc{\sle}{U_q^{\mathbb{K}}L\mathfrak{sl}_2} 

Let us recall how the braid group actions can be used to obtain rank one subalgebras in $\qla$ and $\gOq$. 
For $i \in \indx_0$, let $\omega'_i = \omega_i s_i$, and let $U_i$ be the subalgebra of $\qla$ generated by 
\[
E_i, \ F_i, \ K_i^{\pm 1}, \ T_{\omega'_i}(E_i), \ T_{\omega'_i}(F_i), \ T_{\omega'_i}(K_i^{\pm 1}). 
\]
By \cite[Proposition 3.7]{beck-94}, there is an algebra isomorphism $\iota_i \colon \uLsl \to U_i$ sending 
\[
E_1 \mapsto E_i, \ F_1 \mapsto F_i, \ K_1^{\pm1} \mapsto K_i^{\pm1}, \quad 
E_0 \mapsto T_{\omega'_i}(E_i), \ F_0 \mapsto T_{\omega'_i}(F_i), \ K_0^{\pm1} \mapsto T_{\omega'_i}(K_i^{\pm1}). 
\] 
Let $\Ui$ be the subalgebra of $\qlae$ generated by $U_i$, $\mathbb{K}_i^{\pm1}$ and $(\KK_\delta \KK_i^{-1})^{\pm1}$. The algebra isomorphism $\iota_i$ extends to an isomorphism $\iota_i \colon \sle \to \Ui$ sending  
$\KK_1 \mapsto K_i$ and $\KK_0 \mapsto \KK_\delta \KK_i^{-1}$. For future use, we note that it follows easily from the definitions that 
\eq \label{eq: x- beck iso}
\iota_i(x^-_r) = x^-_{i,r}. 
\eneq

Moreover, for $i \in \indx_0$, let $\gOqi$ be the subalgebra of $\gOq$ generated by $B_i$, $\Tbr_{\omega'_i}(B_i)$, $\KK_i^{\pm1}$ and $(\KK_\delta \KK_i^{-1})^{\pm1}$. By \cite[Proposition 3.9]{lu-wang-21}, there is an algebra isomorphism $\iota \colon \Oq \to \gOqi$ sending 
\eq \label{eq: rank 1 q onsager}
B_1 \mapsto B_i, \quad B_0 \mapsto \Tbr_{\omega'_i}(B_i), \quad \KK_1 \mapsto \KK_i, \quad \KK_0 \mapsto \KK_\delta \KK_i^{-1}. 
\eneq 
We also note that $A_{-1}$ is sent to $\Tbr_{\omega_i}(B_i)$. 

%is the action on ordinary K also given like this? 

%%%%%%%%%%%%%%%%%%%%%%%%%%%%%
%%%%%%%%%%%%%%%%%%%%%%%%%%%%%

\section{Compatibility of braid group actions} 
\label{sec: comp br gr}

The braid group action on $\gOq$ is not compatible with the action on $\qla$, in the sense that the former is not the restriction of the latter. Similarly, the rank one subalgebras from \S \ref{subsec: rk one subalg} are not compatible, i.e., for each $i \in \indx_0$, the diagram 
\eq \label{eq: comm diagram compatible}
\begin{tikzcd}%[ row sep = 0.2cm]
\Oq \arrow[r, "\eta"] \arrow[d, "\iota_i", "\wr"'] & \sle \arrow[d, "\iota_i", "\wr"']   \\
\gOqi \arrow[d, hookrightarrow]  & \Ui \arrow[d, hookrightarrow]  \\
\gOq \arrow[r, "\eta"] & \qlae
\end{tikzcd}
\eneq
does \emph{not} commute. 
Nevertheless, 
the goal of this section is to show that, 
in an appropriate sense, 
diagram \eqref{eq: comm diagram compatible} does commute modulo the ``Drinfeld positive half" of~$\qlae$. 

\Rem
In rank 1, it was sufficient to work with $q$-Onsager algebra with fixed parameters $\mathbf{c}$ and $\mathbf{s}$. In higher rank, we first need to pass to the universal version of the generalized $q$-Onsager algebra. This is necessitated by our analysis of the braid group action, which is non-trivial on the central elements $\mathbb{K}_i$. However, in the main result of this section, Corollary \ref{cor: ultimate factorization theorem}, we can specialize the parameters again. 
\enrem

\subsection{Action on generators} 

Assume that $\g$ is not of type $\mathsf{A}_1$. Set 
\[ \tilde{E}_i = - q^{-2} \KK_i E_i K_i\mi.\]
We begin by investigating the image of elements such as $\Tbr_i(B_{i\pm1})$ inside $\qlae$. 

\Lem \label{lem: action difference 1} 
If $a_{ij} = -1$ then   
\[
\eta(\Tbr_i(B_j)) =  
%\left\{ \begin{array}{l l}
%\KK_j\mi(F_j - q^{-2} \KK_j E_j K_j\mi) & \mbox{if } i=j, \\
%(F_j - q^{-2} \KK_j E_j K_j\mi) & \mbox{if } a_{ij} = 0, \\ 
T_i(\eta(B_j)) +  (q-q\mi) F_i \tilde{E}_j= T_i(\eta(B_j)) + [\tilde{E}_j, F_i]_q. 
%\end{array} \right. 
\]
\enlem

\Proof 
If $a_{ij} = -1$ then  
\begin{align}
T_i(\eta(B_j)) = (F_jF_i - qF_iF_j) - q^{-2}\KK_i\KK_j (E_iE_j - q\mi E_jE_i) K_i\mi K_j\mi
\end{align}
and 
\begin{align}
\eta(\Tbr_i(B_j)) =& \ \eta(B_jB_i - qB_iB_j) \\
=& \ (F_j - q^{-2} \KK_j E_j K_j\mi) (F_i - q^{-2} \KK_i E_i K_i\mi) \\
& \ - q(F_i - q^{-2} \KK_i E_i K_i\mi) (F_j - q^{-2} \KK_j E_j K_j\mi) \label{eq: surv term} \\ 
=& \ (F_jF_i - qF_iF_j) - q^{-2}\KK_i \KK_j (E_iE_j - q\mi E_jE_i)K_i\mi K_j\mi \\
& \ + (q-q\mi)q^{-2} \KK_j F_iE_jK_j\mi, 
\end{align}
since
\eq \label{eq: vanishing term FE'}
F_jE_iK_i\mi - qE_iK_i\mi F_j = 0.
\eneq
Applying \eqref{eq: br act Fi Ei ord} concludes the proof. 
\enproof 

For future use, let us record \eqref{eq: vanishing term FE'} as 
\eq \label{eq: vanishing term FE}
[F_i, \tilde{E}_{i\pm1}] =0. 
\eneq

\subsection{Type $\mathsf{A}$ combinatorics}

From now on assume that $\g$ is of type $\mathsf{A}_N$ with $N \geq 2$. 
Then $\Lambda = \langle \pi \rangle$, where $\pi$ is the automorphism of order $N+1$ of the affine Dynkin diagram of type $\mathsf{A}_N^{(1)}$ sending $i \mapsto i+1 \mod N+1$. The highest root is $\theta = \alpha_1 + \alpha_2 + \cdots + \alpha_n$. 

Let $w_0$ be the longest element of the finite Weyl group $W_0$ of $\g$.
Given $i \in \indx_0$, let $w_i$ the longest element of the Weyl group corresponding to the Dynkin diagram of $\g$ with the node labelled by $i$ removed. 
Given $1 \leq k < l$, let $[k,l] = s_k s_{k+1} \cdots s_{l}$. 

\Pro \label{pro: type A wt}
Let $\widetilde{W}$ be the extended affine Weyl group of type $\mathsf{A}_N^{(1)}$, for $N \geq 1$. Then 
\[
\omega_i = \pi^i [N - i+1 ,N] \cdots[2,i+1][1,i]. 
\] 
This is a reduced expression. 
\enpro

\Proof
It follows from \cite[Ch.\ VI \S 2.3, Theorem 6]{bourbaki} that 
\[ \omega_i = \pi^i w_0w_i. \]
It is easy to check that $w_0w_i = \tau^{N+1-i}$, where $\tau$ is the inverse of the cyclic permutation $(12 \cdots N+1)$, if we identify $W$ with the symmetric group $S_{N+1}$. It is well known that such permutations have an expression of the form $[N - i+1 ,N] \cdots[2,i+1][1,i]$. To show that this expression is reduced, it suffices to check that $\ell(\omega_i) = i(N-i+1)$. By \cite[Proposition 1.23]{IM}, $\ell(\omega_i)$ is equal to the number of positive roots which contain $\alpha_i$ in their expression in terms of the basis of simple roots. Such roots are precisely $\alpha_k + \cdots + \alpha_l$, for $1 \leq k \leq i$ and $i \leq l \leq N$, completing the proof. 
\enproof

For example, let $N=5$. Then:
\[
\omega_1 \mapsto \pi s_5 s_4 s_3 s_2 s_1, \quad \omega_2 \mapsto \pi^2 s_4 s_5 s_3 s_4 s_2 s_3 s_1 s_2, \quad 
\omega_3 \mapsto \pi^3 s_3s_4s_5 s_2 s_3s_4 s_1 s_2 s_3,
\]
\[
\omega_4 \mapsto \pi^4 s_2s_3s_4s_5 s_1s_2s_3s_4, \quad \omega_5 \mapsto \pi^5 s_1s_2s_3s_4s_5. 
\] 

\subsection{Inductive arguments}

Define non-commutative polynomials $P_k(y_1, \cdots, y_k)$ over $\C$ by induction in the following way: %maybe different field in BC?
\begin{alignat}{3}
P_1(y_1) =& \ y_1, \qquad& P_{k+1}(y_1, \cdots, y_k) =& \ [P_k(y_1, \cdots, y_k), y_{k+1}]_q.
\end{alignat}
Clearly, we have 
\begin{align}
P_k(y_1, \cdots, y_k) =& \ P_{k-l+1}(P_l(y_1, \cdots, y_l), y_{l+1}, \cdots, y_k). 
\end{align}
for any $1 \leq l \leq k-1$. 
We will need the following lemma. 

\Lem \label{lem: EF ind}
Let $ 2 \leq m \leq N$.  
Then $\Tbr_{m} \cdots \Tbr_2(B_1) = P_{m}(B_1, \cdots, B_{m})$. Moreover,  
\[ 
\eta(\Tbr_m \cdots \Tbr_2 (B_1)) = \sum_{\substack{\varepsilon_k \geq \varepsilon_{k+1},\\ k=1,...,m-1}} P_m(\tilde{e}_{1}^{\varepsilon_1}, \cdots, \tilde{e}_{m}^{\varepsilon_m}), 
\]
where each $\varepsilon_l \in \{+,-\}$, and $\tilde{e}_i^+ = \tilde{E}_i$, $\tilde{e}_i^- = F_i$.
In other words, the RHS consists of polynomials in letters $\tilde{E}_1, \cdots, \tilde{E}_l, F_{l+1}, \cdots, F_m$ (with each letter occurring exactly once), for $0 \leq l \leq m$. In particular, these polynomials have degree $\sum_{j=1}^l \alpha_j - \sum_{j=l+1}^m \alpha_j$. 
\enlem

\Proof
We prove the statement by induction on $m$. The case $m=2$ reduces to Lemma~\ref{lem: action difference 1}. 
Otherwise, by induction, 
\[ 
\Tbr_{m-1} \cdots \Tbr_2(B_1) = P_{m-1}(B_1, \cdots, B_{m-1}). 
\]
Then 
\begin{align} \label{eq: TBP polynomial}
\Tbr_m \cdots \Tbr_2(B_1) =& \ P_{m-1}(B_1, \cdots, B_{m-2}, \Tbr_m(B_{m-1})) \\
=& \ [P_{m-1}(B_1, \cdots, B_{m-1}), B_m]_q \\ % - qB_m P_{m-1}(B_1, \cdots, B_{m-1}) \\
=& \ P_m(B_1, \cdots, B_m), 
\end{align}
since $B_m$ commutes with $B_1, \cdots, B_{m-2}$. 
This proves the first statement of the lemma. 

Moreover, by induction, 
\[
\eta(P_{m-1}(B_1, \cdots, B_{m-1})) = \sum_{\varepsilon_k \geq \varepsilon_{k+1}} P_{m-1} (\tilde{e}_{1}^{\varepsilon_1}, \cdots, \tilde{e}_{{m-1}}^{\varepsilon_{m-1}}). 
\]
Therefore, it suffices to show that 
\[
P_{m} (\tilde{e}_{1}^{\varepsilon_1}, \cdots, \tilde{e}_{{m-2}}^{\varepsilon_{m-2}}, F_{m-1}, \tilde{E}_m) = 0. 
\] 
But this is indeed the case since 
\[
P_{m} (\tilde{e}_{1}^{\varepsilon_1}, \cdots, \tilde{e}_{{m-2}}^{\varepsilon_{m-2}}, F_{m-1}, \tilde{E}_m) = 
 P_{m-1} (\tilde{e}_{1}^{\varepsilon_1}, \cdots, \tilde{e}_{{m-2}}^{\varepsilon_{m-2}}, [F_{m-1}, \tilde{E}_m]_q)
\] 
and $[F_{m-1}, \tilde{E}_m]_q = 0$ by \eqref{eq: vanishing term FE}.  
\enproof 

Lemma \ref{lem: EF ind} has the following ``dual". 

\Lem \label{lem: EF ind2}
Let $ 2 \leq m \leq N$.  
Then $\Tbr_{1} \cdots \Tbr_{m-1}(B_m) = P_{m}(B_m, \cdots, B_{1})$. Moreover,  
\[ 
\eta(\Tbr_1 \cdots \Tbr_{m-1} (B_m)) = \sum_{\substack{\varepsilon_k \leq \varepsilon_{k+1},\\ k=1,...,m-1}} P_m(\tilde{e}_{m}^{\varepsilon_m}, \cdots, \tilde{e}_{1}^{\varepsilon_1}), 
\]
In other words, the RHS can be expressed as a sum of polynomials consisting of letters $\tilde{E}_m, \cdots, \tilde{E}_l, F_{l-1}, \cdots, F_1$, for $1 \leq l \leq m+1$. In particular, these polynomials have degree $-\sum_{j=1}^{l-1} \alpha_j +  \sum_{j=l}^m \alpha_j$. 
\enlem 

\Proof
The proof is analogous to that of Lemma \ref{lem: EF ind}. 
\enproof

\subsection{Weak compatibility} 

We first recall the following fact. 

\Lem \label{lem: deg dr of E0} \label{pro: B-1 compat br gr action}
We have $\degdr E_0 = - \theta = - \sum_{j=1}^N \alpha_j$ and $\degdr F_0 = \theta = \sum_{j=1}^N \alpha_j$. 
\enlem

\Proof
This follows directly from the explicit description of the images of $E_0$ and $F_0$ in the Drinfeld realization, see, e.g., \cite[Remark 4.7]{beck-94} and \cite[Theorem 2.2]{chari-pressley-94}. 
\enproof 

Let us write $\eta(\Tbr_{\omega'_i}(B_i))$ as a sum of $Q$-homogeneous terms. We say that a term is \emph{irrelevant} if its degree is not $\pm \alpha_i$. 
The theorem below is our first result on the compatibility of the two braid group actions. It is the essential part of Theorem \ref{intro: br gr} from the Introduction. 

\Thm \label{pro: irrelevant terms}
Let $1 \leq i \leq N$. 
Every irrelevant term of $\eta(\Tbr_{\omega'_i}(B_i))$ lies in $U_{d_i = 1,+}$. 
Hence 
\begin{alignat}{3}
\eta(\mathbf{T}_{\omega_i'}(B_i)) &\equiv& \ T_{\omega_i'}(\eta(B_i)) &\mod U_{d_i = 1,+}. 
\end{alignat}
\enthm

\Proof 
By Lemma \ref{lem: EF ind2}, 
\[
\Tbr_1 \cdots \Tbr_{i-1} (B_i) = [P_{i-1}(B_i, \cdots, B_2),B_1]_q. 
\]
Next, observe that 
\begin{align}
(s_{N-i+1} \cdots s_{N}) \cdots (s_2 \cdots s_{i+1}) (\alpha_1) =& \ s_{N-i+1}\cdots s_2 (\alpha_{1}), \\ 
(s_{N-i+1} \cdots s_{N}) \cdots (s_2 \cdots s_{i+1}) (\alpha_k) =& \ \alpha_{k+N-i}, 
\end{align}
for $2 \leq k \leq i$. 
Therefore, Lemmas \ref{lem: br act on wt} and \ref{lem: EF ind}-\ref{lem: EF ind2} imply that 
\begin{align}
\Tbr_{\omega'_i}(B_i) =& \ \pi^i [P_{i-1}(B_{N}, \cdots, B_{N-i+2}), \Tbr_{N-i+1} \cdots \Tbr_2 (B_1)]_q \\
=& \ \pi^i [P_{i-1}(B_{N}, \cdots, B_{N-i+2}), P_{N-i+1}(B_1, \cdots, B_{N-i+1})]_q \\
=& \ \pi^i 
\sum_{\substack{\varepsilon_k \leq \varepsilon_{k+1},\\ k=N-i+2,...,N}} \sum_{\substack{\varepsilon_l \geq \varepsilon_{l+1},\\ l=1,...,N-i}} [ P_{i-1}(\tilde{e}_{N}^{\varepsilon_N}, \cdots, \tilde{e}_{N-i+2}^{\varepsilon_{N-i+2}})  , P_{N-i+1}(\tilde{e}_{1}^{\varepsilon_1}, \cdots, \tilde{e}_{N-i+1}^{\varepsilon_{N-i+1}})  ]_q \\
=& \  \label{eq: TBi q comm}
\sum_{\substack{\varepsilon_k \leq \varepsilon_{k+1},\\ k=1,...,i-1}} \sum_{\substack{\varepsilon_l \geq \varepsilon_{l+1} \geq \varepsilon_0,\\ l=i+1,...,N-1}} [ P_{i-1}(\tilde{e}_{i-1}^{\varepsilon_{i-1}}, \cdots, \tilde{e}_{1}^{\varepsilon_{1}})  , P_{N-i+1}(\tilde{e}_{i+1}^{\varepsilon_{i+1}}, \cdots, \tilde{e}_{0}^{\varepsilon_{0}})  ]_q. \\
\end{align} 

By Lemma \ref{lem: deg dr of E0}, each irrelevant term inside of the summation in \eqref{eq: TBi q comm} is already in $U_{d_i = 1,+}$ unless $\varepsilon_0 = +$. So suppose that $\varepsilon_0 = +$ and that at least one of the superscripts $\varepsilon_k = -$. This forces $\varepsilon_1 = -$ and $\varepsilon_{i+1} = \cdots = \varepsilon_0 = +$. We then have 
\begin{align}
 [ P_{i-1}(\tilde{e}_{i-1}^{\varepsilon_{i-1}}, \cdots, \tilde{e}_{1}^{\varepsilon_{1}}) ,& \ P_{N-i+1}(\tilde{e}_{i+1}^{\varepsilon_{i+1}}, \cdots, \tilde{e}_{0}^{\varepsilon_{0}})  ]_q = \\ 
=& \ [ [P_{i-2}(\tilde{e}_{i-1}^{\varepsilon_{i-1}}, \cdots, \tilde{e}_{2}^{\varepsilon_{2}}), F_1]_q, 
[ P_{N-i}(\tilde{E}_{i+1}, \cdots, \tilde{E}_N), \tilde{E}_{0}]_q ]_q \\ 
=& \ [ P_{N-i}(\tilde{E}_{i+1}, \cdots, \tilde{E}_N), [P_{i-2}(\tilde{e}_{i-1}^{\varepsilon_{i-1}}, \cdots, \tilde{e}_{2}^{\varepsilon_{2}}), [F_1, \tilde{E}_0]_q]_q]_q \\
=& \ 0
\end{align}
since, by \eqref{eq: vanishing term FE}, $[F_1, \tilde{E}_0]_q = 0$. 
We conclude that every irrelevant term in \eqref{eq: TBi q comm} is indeed in $U_{d_i = 1,+}$. 
Moreover, the relevant terms are precisely
\begin{align}
 [ P_{i-1}(\tilde{E}_{i-1}, \cdots, \tilde{E}_{1})  , P_{N-i+1}(\tilde{E}_{i+1}, \cdots, \tilde{E}_{0})  ]_q =& \ T_{\omega_i'}(\tilde{E}_i), \\ 
 [ P_{i-1}({F}_{i-1}, \cdots, {F}_{1})  , P_{N-i+1}({F}_{i+1}, \cdots, {F}_{0})  ]_q =& \ T_{\omega_i'}({F}_i). 
\end{align}
 This concludes the proof. 
\enproof

\Cor[Compatibility theorem] \label{thm: final br gr act compat}
The diagram \eqref{eq: comm diagram compatible} commutes modulo $U_{\neq i,+}$,  
for each $1 \leq i \leq N$. 
\encor

\Proof
Let $\mu$ denote the composition of the east and south arrows ($\reflectbox{\rotatebox[origin=c]{90}{$\Lsh$}}$), and $\nu$ the composition of the south and east arrows ($\reflectbox{\rotatebox[origin=c]{180}{$\Rsh$}}$) in \eqref{eq: comm diagram compatible}. We have $\mu(B_0) = \nu(B_0)$ and, by Theorem \ref{pro: irrelevant terms}, $\mu(B_{-1}) \equiv \nu(B_{-1})$ modulo $U_{d_i=1,+}$. Since $B_0$ and $B_{-1}$ generate $\Oq$, and $\mu(B_0)$ contains terms of degree $\pm \alpha_i$, it follows that the diagram commutes modulo~$U_{\neq i,+}$. 
\enproof

\subsection{Strong compatibility}

The goal of this section is to strengthen the results of the preceding subsection to pave the way for a higher rank generalization of the factorization theorem. %showing that diagram \eqref{eq: comm diagram compatible} commutes modulo $U_+$. 
Set $\CCC_i = \CCC\mi \mathbb{K}_i$. 
It is easily verified that %ind proof? 
\[
\Tbr_{\omega'_i}(\KK_i) = \CCC_i\mi, \qquad T_{\omega'_i}(K_i) = K_i\mi.
\]

In the calculations below, we will frequently make use of a remainder term $Q_i$. The statement of the following lemma contains its definition. 

\Lem 
For each $i \in \indx_0$: 
%\[
%x_{i,-1}^+ = - T_{\omega'_i}(F_i)K_i\mi, \qquad x^-_{i,1} = - K_i T_{\omega'_i}(E_i). 
%\] 
\[
\eta(A_{i,-1}) = x^-_{i,1} - \CCC_ix^+_{i,-1}K_i + Q_i, \qquad Q_i \in U_{d_i \geq 1,+}. 
\]
Hence 
\[
\eta(A_{i,-1}) \equiv \iota_i \eta(A_{-1}) \quad \mod  U_{d_i \geq 1,+}. 
\]
\enlem

%application to max comm Cartan
%Molev 
%stronger statement - mod d_i>1

\Proof
By \cite[Proposition 4.7]{beck-94}\footnote{More precisely, in the first equality below, an adjustment by a sign function $o(i)$  which alternates on adjacent vertices is necessary. An analogous adjustment occurs in the formula for $A_{i,-1}$. Since we work with a fixed vertex $i \in \indx_0$ most of the time, we may assume for simplicity that $o(i)=1$.}, 
\begin{align} 
\label{eq: x+beck}
x_{i,-1}^+ =& \ T_{\omega_i}(E_i) = T_{\omega'_i}(- F_iK_i) = - T_{\omega'_i}(F_i)K_i\mi, \\
\label{eq: x-beck}
x_{i,1}^- =& \ T_{\omega_i}(F_i) = T_{\omega'_i}(- K_i\mi E_i) = - K_iT_{\omega'_i}(E_i). 
\end{align} 
Moreover, by \cite[Theorem 3.13]{lu-wang-21}, 
\[
A_{i,-1} = \Tbr_{\omega_i}(B_i) = \Tbr_{\omega'_i}(\KK_i\mi B_i) = \CCC_i \Tbr_{\omega'_i}(B_i). 
\]
Hence, by Theorem \ref{pro: irrelevant terms} and \eqref{eq: x+beck}--\eqref{eq: x-beck}, 
\[
\eta(A_{i,-1}) = \CCC_i T_{\omega'_i}(F_i - K_i\mi\KK_iE_i) + Q_i = -\CCC_i x^+_{i,-1} K_i + x^-_{i,1} + Q_i, 
\]
completing the proof.  
\enproof 

%Tosidlo
%K-matrix note 
%Hall algebra inclusion/coideal structure geometrically 

%Prop: all Ak 
%zywia sie wspomnieniami 

\Lem \label{lem: Hi1 app}
We have 
\[
\eta(H_{i,1}) \equiv \iota_i\eta(H_1) + q^2 \CCC_i\mi [Q_i, F_i]_{q^{-2}} \quad \mod  U_{d_i \geq 2,+}.
\]
\enlem 

\Proof
We have 
\begin{align}
\eta(H_{i,1}) =& \ q^2 \CCC_i\mi \eta( [A_{i,-1}, A_{i,0}]_{q_i^{-2}}) \\ 
=& \ q^2 \CCC_i\mi [\iota_i \eta(A_{-1}) + Q_i, \iota_i \eta(A_{0})]_{q^{-2}} \\ 
=& \ \iota_i \eta(H_1) + q^2 \CCC_i\mi [Q_i, F_i + \tilde{E}_i]_{q^{-2}}. 
\end{align}
Since $[Q_i, \tilde{E}_i]_{q^{-2}} \in U_{d_i \geq 2,+}$, the lemma follows. 
\enproof

\Pro \label{pro: A_i general compat}
For all $r\in \Z$: 
\[
\eta(A_{i,r}) \equiv \iota_i \eta(A_{r}) \quad \mod  U_{d_i \geq 1,+}. 
\]
\enpro 

\Proof 
The result holds for $r = -1,0$ by Theorem \ref{pro: irrelevant terms}. The general case is proven using upward and downward induction. We explicitly do the former, the latter being analogous. 
Relation \eqref{eq: grel2} implies that 
\begin{align}
A_{i, r+1} = [2]\mi [H_{i,1}, A_{i, r}] + \CCC A_{i,r-1}. 
\end{align} 
Using Lemma \ref{lem: Hi1 app}, let us write 
$\eta(H_{i,1}) = \iota_i\eta(H_1) + q^2 \CCC_i\mi [Q_i, F_i]_{q^{-2}} + R$, with $R \in U_{d_i \geq 2,+}$. 
By induction, we may assume that 
\[ \eta(A_{i,r-1}) \equiv \iota_i \eta(A_{r-1}), \quad [ \iota_i \eta(H_{1}), \eta(A_{i, r})] \equiv \iota_i \eta([H_1, A_r])\] modulo $U_{d_i \geq 1,+}$. %Moreover, by ... %my paper
Therefore, it suffices to show that 
\[
[ q^2 \CCC_i\mi [Q_i, F_i]_{q^{-2}} + R, \eta(A_{i,r}) ] \in U_{d_i \geq 1,+}. 
\] 

By induction, we can write $\eta(A_{i,r}) = \iota_i \eta(A_r) + S$, with $S \in U_{d_i \geq 1,+}$. On the other hand, by Lemma \ref{lem: Ar-1} and \eqref{eq: x- beck iso},  
$\iota_i \eta(A_r) = x_{i,-r}^- + S'$ with $S' \in U_{[i],+}$, where $U_{[i],+}$ is the subalgebra of $U_{[i]}$ spanned by elements of positive degree (with respect to $\alpha_i$). We observe immediately that 
\[
[R, \eta(A_{i,r}) ] + [q^2 \CCC_i\mi [Q_i, F_i]_{q^{-2}}, S + S'] \in U_{d_i \geq 1,+}. 
\] 
Therefore, it suffices to show that 
\eq \label{eq: QF vanish}
[[Q_i, F_i]_{q^{-2}}, x_{i,-r}^-] \in U_{d_i \geq 1,+}. 
\eneq
In fact, we show that \eqref{eq: QF vanish} vanishes. 
This will be achieved in \S \ref{sec: completion of proof} below. 
\enproof

%corollary: lattice on all of U qsp. 

\subsection{Completion of the proof} \label{sec: completion of proof}

\nc{\tE}{\tilde{E}}

We first state some auxiliary lemmas. Define the non-commutative polynomial $P'_k(y_1, \cdots, y_k)$ over $\C$ by induction in the following way: %maybe different field in BC?
\begin{alignat}{3}
P'_1(y_1) =& \ y_1, \qquad& P'_{k+1}(y_1, \cdots, y_k) =& \ P'_k(y_1, \cdots, y_{k-1}, [y_k,y_{k+1}]_q), 
\end{alignat}
We have 
\begin{align} \label{eq: P' property}
P'_k(y_1, \cdots, y_k) =& \ P'_{l+1}(y_1, \cdots, y_l, P'_{k-l}(y_{l+1}, \cdots, y_k)), \\
\end{align}
for any $1 \leq l \leq k-1$. 

We say that a tuple $(y_1, \cdots, y_k)$ is \emph{almost commuting} if $y_my_n = y_ny_m$, for $|m-n| > 1$.

\Lem \label{lem: PnP'}
The polynomials $P_k$ and $P'_k$ are equal if $(y_1, \cdots, y_k)$ is almost commuting. 
\enlem

%We have 
%\begin{align}
%\Tbr_{\omega'_1}(B_1) =& \ P_{n-1}(B_2,\cdots,B_{n-2}, B_n, P_{n-1}(B_{n-1}, \cdots, B_{2},B_0)). 
%=& \ P(B_2, \cdots, B_k, P(B_{k+1}, \cdots, B_{))
%\end{align}
\Proof 
The first statement is proven by induction. By definition, $P_1 = P'_1$ and $P_2 = P'_2$. Moreover, 
\begin{align}
P'_{k+1}(y_1, \cdots, y_{k+1}) =& \ P'_2(y_1, P'_k(y_2, \cdots,y_{k+1})) \\ 
=& \  P'_2(y_1, P_k(y_2, \cdots,y_{k+1})) \\ 
=& \ P_k(P'_2(y_1, y_2), y_3, \cdots, y_{k+1}) \\
=& \  P_k(P_2(y_1, y_2), y_3, \cdots, y_{k+1}) 
= P'_{k+1}(y_1, \cdots, y_{k+1}). \qedhere
\end{align}
\enproof

\Lem \label{lem: aux lemma str 1}
Let $2 \leq i \leq n-1$. Then 
\be
\item $[\tE_i, [\tE_{i-1}, Y]_q]_q = q [E_{i}, [\tE_{i-1} K_i\mi, Y]_q]_{q^2}$ if $[K_i\mi, Y]=0$, 
\item $[[\tE_{i-1}, Y]_q, F_i]_{q^{-2}} = [\tE_{i-1}, [Y, F_i]_{q^{-1}}]$, 
\item $[[\tE_{i+1}, Y]_q, F_i]_{q^{-1}} = [\tE_{i+1}, [Y,F_i]]$, 
\ee
for any $Y \in \qlae$. 
If $[Y,F_i] = 0$ then 
\eq \label{eq: EEY van} 
[[\tE_{i-1}, [\tE_{i+1}, Y]_q]_q, F_i]_{q^{-2}} = 0. 
\eneq
\enlem 

\Proof
The first three statements follow by a direct calculation. The fourth statement follows from the second and third if we substitute $[\tE_{i+1}, Y]_q$ for the $Y$ in (2). 
\enproof

\Lem \label{lem: Z x van}
Let $i-1 \leq k \leq 1$. 
The element 
\[
Z = [[P_{i-1}(\tE_{i-1}, \cdots, \tE_k, F_{k-1}, \cdots, F_1), P_{N-i+1}(F_{i+1}, \cdots, F_N, F_0)]_q, F_i]_{q^{-2}}
\]
can be expressed as a polynomial in $x_{j,r}^+$ $(j \neq i)$, with coefficients in $\C$. Hence 
\eq \label{eq: Z x van}
[Z, x_{i,s}^-] = 0
\eneq
for all $s \in \Z$. 

Similarly, if $i+1 \leq k \leq n$ or $k=0$, then 
\[
Z' = [[P_{i-1}(F_{i-1}, \cdots, F_1), P_{N-i+1}(\tE_{i+1}, \cdots, \tE_k, F_{k-1}, \cdots F_N, F_0)]_q, F_i]_{q^{-2}}
\]
can also be expressed as a polynomial in $x_{j,r}^+$ $(j \neq i)$, with coefficients in $\C$. Hence 
\eq \label{eq: Z x van2}
[Z', x_{i,s}^-] = 0.
\eneq
\enlem

\Proof
We will only prove the statement for $Z$, the one for $Z'$ being analogous. 
The equality \eqref{eq: Z x van} follows directly from relation \eqref{eq: g x+-} and $Z$ being a polynomial in $x_{j,r}^+$ $(j \neq i)$. Therefore, it suffices to prove the latter. Using \eqref{eq: P' property} and Lemma \ref{lem: PnP'}, we have 
\begin{align}
[P_{i-1}(\tE_{i-1}, \cdots, \tE_k,& \ F_{k-1}, \cdots, F_1), P_{N-i+1}(F_{i+1}, \cdots, F_N, F_0)]_q = \\ 
=& \ P'_{i}(\tE_{i-1}, \cdots, \tE_k, F_{k-1}, \cdots, F_1, P'_{N-i+1}(F_{i+1}, \cdots, F_N, F_0)) \\
=& \ P'_{N}(\tE_{i-1}, \cdots, \tE_k, F_{k-1}, \cdots, F_1, F_{i+1}, \cdots, F_N, F_0).
\end{align}

In the following, let us abbreviate 
\[
\text{``$\vec{F}$\ '' $=$ ``$F_{k{\shortminus}1}, \cdots, F_1, F_i, F_{i+1}, \cdots, F_N, F_0$''.} 
\] 
Direct calculation using Lemma \ref{lem: aux lemma str 1} shows that 
\begin{align}
Z =& 
\ {\shortminus}q[\tE_{i{\shortminus}1}, P'_N(F_i, \tE_{i{\shortminus}2}, \cdots, \tE_k, F_{k{\shortminus}1}, \cdots, F_1, F_{i+1}, \cdots, F_N, F_0)] \\ 
=& \ {\shortminus}q [\tE_{i{\shortminus}1}, P'_N(\tE_{i{\shortminus}2}, \cdots, \tE_k, F_{k{\shortminus}1}, \cdots, F_1, F_i, F_{i+1}, \cdots, F_N, F_0)] \\ 
=& \ {\shortminus}q [\tE_{i{\shortminus}1}, [E_{i{\shortminus}2}, P'_{N{\shortminus}1}(E_{i{\shortminus}3}, \cdots, E_k K_{i{\shortminus}2}\mi\cdots K_{k}\mi, \vec{F})] ] \\ 
=& \ {\shortminus}q^2 [E_{i{\shortminus}1}, [E_{i{\shortminus}2}K_{i{\shortminus}1}\mi, P'_{N{\shortminus}1}(E_{i{\shortminus}3}, \cdots, E_k K_{i{\shortminus}2}\mi\cdots K_{k}\mi, \vec{F})]_{q\mi} ]_{q^2} \\ 
=& \ {\shortminus}q^2 [E_{i{\shortminus}1}, [E_{i{\shortminus}2}, [E_{i{\shortminus}3}K_{i{\shortminus}1}\mi, P'_{N{\shortminus}2}(E_{i{\shortminus}4}, \cdots, E_k K_{i{\shortminus}2}\mi\cdots K_{k}\mi, \vec{F})]]_{q\mi} ]_{q^2} \\ 
=& \ {\shortminus}q^2 [E_{i{\shortminus}1}, [E_{i{\shortminus}2}, P'_{N{\shortminus}1}(E_{i{\shortminus}3}, E_{i{\shortminus}4}, \cdots, E_k K_{i{\shortminus}1}\mi\cdots K_{k}\mi,\vec{F})]_{q\mi} ]_{q^2} \\ 
=& \ q^3(q{\shortminus}q\mi) 
[E_{i{\shortminus}1}, [E_{i{\shortminus}2}, P'_{i{\shortminus}k{\shortminus}1}(E_{i{\shortminus}3}, E_{i{\shortminus}4}, \cdots, \\
& \ E_{k+1}, E_k  P'_{N{\shortminus}i+k+2}(\vec{F})K_{i{\shortminus}1}\mi\cdots K_{k}\mi)]_{q\mi} ]_{q^2}. 
\end{align} 

Therefore, it suffices to show that $P'_{N{\shortminus}i+k+2}(\vec{F})K_{i{\shortminus}1}\mi\cdots K_{k}\mi$ can be expressed as a polynomial in $x_{j,r}^+$ $(j \neq i)$. 
First assume that $i > 1$. 
We claim that 
\[
P'_{N{\shortminus}i+k+2}(\vec{F})K_{i{\shortminus}1}\mi\cdots K_{k}\mi = 
P'_{i-k}(x_{i-1,-1}^+, E_{i-2}, \cdots, E_{k}) =: Y. 
\] 
Let 
\[
\zeta = \pi^{i-1} [N-i+2, N]\cdots [2, i]. 
\]
Indeed, 
\begin{align}
Y =& \ P'_{i-k}(x_{i-1,-1}^+, E_{i-2}, \cdots, E_{k}) 
= T_{\omega_{i-1}}P'_{i-k}(E_{i-1}, E_{i-2}, \cdots, E_{k}) \\ 
=& \ T_{\omega_{i-1}}T_k \cdots T_{i-2}(E_{i-1}) 
= T_{\zeta} T_{[1,k]}(E_k) \\ 
=& \ -T_{\zeta} T_{[1,k-1]}(F_kK_k) 
= -T_{\zeta} (P'_k(F_k, \cdots, F_1) K_1 \cdots K_k) \\ 
=& \ -T_{\pi^{i-1}} (P'_k(F_{N+k-i+1}, \cdots, F_{N-i+3}, F_{1}, \cdots, F_{N-i+2})\cdot \\
& \ \cdot K_{N+k-i+1} \cdots K_{N-i+3} K_1 \cdots K_{N-i+2}) \\ 
=& \ - P'_k(F_{k-1}, \cdots, F_{1}, F_{i}, \cdots, F_{N}, F_0) K_{k-1} \cdots K_{1} K_i \cdots K_{N}K_0 \\ 
=& \ - P'_{N{\shortminus}i+k+2}(\vec{F})K_{i{\shortminus}1}\mi\cdots K_{k}\mi, %explain more
\end{align} 
proving the claim. The case $i=1$ is similar but easier, so we leave it to the reader. 
\enproof

\Pro %take notation out. We need P plus and minus, depending on sign og epsilon i+1 - another letter? 
We have 
\[
[ [Q_i, F_i]_{q^{-2}} , x^-_{i,r}] = 0. 
\]
\enpro

\Proof 
%First let $2 \leq i \leq n-1$. relab 
By the definition of $Q_i$ and \eqref{eq: TBi q comm}, we have 
\eq \label{eq: sum to vanish}
Q_i = \CCC_i \sum_{\substack{\varepsilon_k \leq \varepsilon_{k+1},\\ k=1,...,i-1}} \sum_{\substack{\varepsilon_l \geq \varepsilon_{l+1} \geq \varepsilon_0,\\ l=i+1,...,N-1}} [ P_{i-1}(\tilde{e}_{i-1}^{\varepsilon_{i-1}}, \cdots, \tilde{e}_{1}^{\varepsilon_{1}})  , P_{N-i+1}(\tilde{e}_{i+1}^{\varepsilon_{i+1}}, \cdots, \tilde{e}_{0}^{\varepsilon_{0}})  ]_q, 
\eneq
where $(\varepsilon_{i-1}, \varepsilon_{i+1}) \in \{(+,-), (-,+), (+,+)\}$. 
Let us denote each summand on the RHS of \eqref{eq: sum to vanish}, corresponding to a sequence $\vec{\varepsilon} = (\varepsilon_{i-1}, \cdots \varepsilon_1, \varepsilon_{i+1}, \cdots \varepsilon_N, \varepsilon_0)$, by $P_{\vec{\varepsilon}}$. 

It follows from Lemma \ref{lem: PnP'} and \eqref{eq: P' property} that 
\begin{align}
 P_{i-1}(\tilde{e}_{i-1}^{\varepsilon_{i-1}}, \tilde{e}_{i-2}^{\varepsilon_{i-2}}, \cdots, \tilde{e}_1^{\varepsilon_1}) =& \ [ \tilde{e}_{i-1}^{\varepsilon_{i-1}}, P_{i-2}(\tilde{e}_{i-2}^{\varepsilon_{i-2}}, \cdots, \tilde{e}_1^{\varepsilon_1})]_q, \\
  P_{N-i+1}(\tilde{e}_{i+1}^{\varepsilon_{i+1}}, \tilde{e}_{i+2}^{\varepsilon_{i+2}}, \cdots, \tilde{e}_0^{\varepsilon_0}) =& \ [ \tilde{e}_{i+1}^{\varepsilon_{i+1}}, P_{N-i}(\tilde{e}_{i+2}^{\varepsilon_{i+2}}, \cdots, \tilde{e}_0^{\varepsilon_0})]_q. 
\end{align}
Let us abbreviate $Y_1 = P_{i-2}(\tilde{e}_{i-2}^{\varepsilon_{i-2}}, \cdots, \tilde{e}_1^{\varepsilon_1})$ and $Y_2 = 
P_{N-i}(\tilde{e}_{i+2}^{\varepsilon_{i+2}}, \cdots, \tilde{e}_0^{\varepsilon_0})$. Then 
\begin{align}
P_{\vec{\varepsilon}} =& \ [ [\tilde{e}_{i-1}^{\varepsilon_{i-1}} , Y_1]_q, [\tilde{e}_{i+1}^{\varepsilon_{i+1}} , Y_2]_q ]_q \\
=& \ (1-q) [\tilde{e}_{i-1}^{\varepsilon_{i-1}}, [\tilde{e}_{i+1}^{\varepsilon_{i+1}}, Y_1Y_2]_q]_q. 
\end{align} 

If $(\varepsilon_{i-1}, \varepsilon_{i+1}) = (+,+)$, then \eqref{eq: EEY van} in Lemma \ref{lem: aux lemma str 1} now implies that $[P_{\vec{\varepsilon}}, F_i]_{q^{-2}} = 0$. 
Otherwise, $[ [P_{\vec{\varepsilon}}, F_i]_{q^{-2}} , x^-_{i,r}] = 0$ follows directly from \eqref{eq: Z x van}--\eqref{eq: Z x van2} in Lemma \ref{lem: Z x van}. 
\enproof

\subsection{Factorization theorem in higher rank} 
%and coproduct

We can now generalize the factorization and coproduct theorems to higher ranks. 
Let us choose parameters $\mathbf{c}, \mathbf{s}$, which determines the coideal structure through 
\eqref{eq: Kolb emb gen}. 

\Cor \label{cor: ultimate factorization theorem} 
Let $\g$ be of type $\mathsf{A}$. Then: 
\begin{align}
\eta_{\mathbf{s}}(\Thgsr_i(z)) \equiv& \ \xi_{\mathbf{s}}(\Thgsr_i(z)) \pmb{\phi}_i^-(z\mi)\pmb{\phi}_i^+(\CC z) \fext{\mod U_{+}}{z}, \\
\Delta_{\mathbf{s}}(\Thgsr_i(z)) \equiv& \ \eta_{\mathbf{s}}(\Thgsr_i(z)) \otimes \eta(\Thgsr_i(z)) \quad  \mod \fext{(\Oqp{\mathbf{c},\mathbf{s}}(\widehat{\g}) \otimes U_{+})}{z}.  
\end{align} 
\encor 

\Proof
If $\mathbf{s} = 0$, the corollary follows from Corollary \ref{cor: fact} and Theorem \ref{pro: A_i general compat}. 
Let us explain this more precisely. In analogy to \eqref{eq:Thr}, the defining relations of $\gOq$ imply that
\begin{align}
\label{eq: gThr}
\Thgsr_i(z)= 1 + \frac{\cc_i^{-1}\CC z^2\left(z\mi[A_{i,-1}, \Apsi]_{q^{-2}}-q^{-2}[A_{i,0},\Apsi]_{q^2}\right)}{1-\CC z^2}. 
\end{align}
By Proposition \ref{pro: A_i general compat}, we can write $\eta(\Apsi) = \iota_i \eta(\Aps) + R(z)$ with $R(z) \in \fext{U_{d_i \geq 1, +}}{z}$. Moreover, $\eta(A_{i,-1}) = \iota_i \eta (A_{-1}) + R_{-1}$, with $R_{-1} \in U_{d_i \geq 1, +}$, and $\eta(A_{i,0}) = \iota_i \eta(A_0)$. We call $R(z)$ and $R_{-1}$ ``remainder terms''. 
Since the elements $\iota_i \eta (A_{-1})$ and $\iota_i \eta(A_0)$ contain precisely homogeneous terms of degree $\pm \alpha_i$, multiplying them with remainder terms yields an element of $U_+$. Similarly, a product of remainder terms lies in $U_+$. 
It follows that 
\begin{align}
\eta(\Thgsr_i(z)) \equiv& \ \iota_i \eta \left(1 + \frac{\cc_i^{-1}\CC z^2\left(z\mi[A_{-1}, \Aps]_{q^{-2}}-q^{-2}[A_{0},\Aps]_{q^2}\right)}{1-\CC z^2}\right) \\
=& \ \iota_i \eta (\Thgsr(z)) \equiv \iota_i(\pmb{\phi}^-(z\mi)\pmb{\phi}^+(\CC z)) = \pmb{\phi}_i^-(z\mi)\pmb{\phi}_i^+(\CC z)
\end{align} 
modulo $\fext{U_{+}}{z}$. In the third equivalence, Corollary \ref{cor: fact} was used. 

The case of non-standard parameters $\mathbf{s} \neq 0$ follows from a higher rank analogue of \eqref{eq: eta cs delta} and Theorem \ref{cor: group like Theta}, using the fact that $\overline{\Delta}_{\mathbf{s}}(U_{+}) \subseteq \C \otimes U_{+}$. 
The group-like property of the coproduct is deduced in the same way as in the proof of Theorem~\ref{cor: group like Theta}. 
\enproof

%\Rem
%Although not necessary for our applications, it would be interesting to find out whether Corollary \ref{cor: ultimate factorization theorem} can be strengthened so that $U_{+}^i$ and $U^i_{-}$ are replaced by $U_+$ and $U_-$, respectively. Preliminary calculations suggest that this is the case at least in rank two. 
%\enrem

Thanks to Corollary \ref{cor: ultimate factorization theorem}, 
Corollary \ref{cor: FR thm Oq} and Theorem \ref{cor: Drinfeld rational functions} carry over to higher rank, as formulated in the introduction (Theorem \ref{intro thm2}). 

%analogue of compatibility theorem with s 

%%%%

\providecommand{\bysame}{\leavevmode\hbox to3em{\hrulefill}\thinspace}
\providecommand{\MR}{\relax\ifhmode\unskip\space\fi MR }
\providecommand{\MRhref}[2]{%
  \href{http://www.ams.org/mathscinet-getitem?mr=#1}{#2}
}
\providecommand{\href}[2]{#2}

\end{document}